\documentclass[11pt,reqno]{article}

\usepackage[usenames]{color}
\usepackage{amssymb}
\usepackage{amsmath}
\usepackage{amsthm}
\usepackage{amsfonts}
\usepackage{amscd}
\usepackage{enumerate}
\usepackage{amsfonts}
\usepackage{graphicx}
\usepackage[toc,page]{appendix}
\usepackage[margin=0pt,font+=small,labelformat=parens,labelsep=space, skip=6pt,list=false,hypcap=false] {caption}
\usepackage{subfig}
\usepackage{relsize}
\usepackage{url}
\usepackage{amsthm}
\usepackage[colorlinks=true,
linkcolor=webgreen,
filecolor=webbrown,
citecolor=webgreen]{hyperref}

\definecolor{webgreen}{rgb}{0,.5,0}
\definecolor{webbrown}{rgb}{.6,0,0}

\usepackage{color}
\usepackage{fullpage}
\usepackage{float}
\usepackage{graphics}
\usepackage{latexsym}
\usepackage{epsf}
\usepackage{breakurl}

\numberwithin{equation}{section}

\usepackage[top=.9in,bottom=.9in]{geometry}
 \numberwithin{equation}{section}
 
\begin{document}

\parskip 8pt
\parindent 20pt
\theoremstyle{plain}
\newtheorem{theorem}{Theorem}
\newtheorem{Coroly}[theorem]{Corollary}
\newtheorem{lemma}[theorem]{Lemma}
\newtheorem{Def}[theorem]{Definition}
\newtheorem{prop}[theorem]{Proposition}

\theoremstyle{definition}
\newtheorem{definition}[theorem]{Definition}
\newtheorem{example}[theorem]{Example}
\newtheorem{conjecture}[theorem]{Conjecture}

\theoremstyle{remark}
\newtheorem{rem}[theorem]{Remark}

\pagestyle{plain}
\numberwithin{equation}{section}
\iftrue
\title{Determining the Indeterminate: On the Evaluation of Integrals that connect Riemann's, Hurwitz' and Dirichlet's Zeta, Eta and Beta functions.}
\author{ Michael Milgram\footnote{mike@geometrics-unlimited.com}\\{Consulting Physicist, Geometrics Unlimited, Ltd.}}
\date{March 13, 2023}
\maketitle
\begin{flushleft} \vskip 0.3 in 
\centerline{Box 1484, Deep River, Ont. Canada. K0J 1P0}
\vskip .2in
\centerline{Author's manuscript, Orcid:0000-0002-7987-0820}
\vskip 0.2 in
\centerline{\bf REVISED:  {\date{\today}}}
\centerline{}
\vskip .1in
\fi
\centerline{}
\vskip .1in
MSC classes: 	11M06, 11M26, 11M35, 11M41, 26A09, 26A15, 26A30, 30B40, 30E20, 32D15, 32D20, 33B99
\vskip 0.1in
Keywords: Riemann Zeta Function, Dirichlet Eta function, alternating Zeta function, Hurwitz Zeta function, infinite series, evaluation of integrals, analytic continuation, regularization, holomorphic continuation, edge-of-the-wedge, indeterminate, essential singularity
\vskip 0.1in 

\centerline{\bf Abstract}\vskip .3in
By applying the inverse Mellin transform to some simple closed form identities, a number of relationships are established that connect integrals containing Riemann's and Hurwitz' zeta functions ($\zeta(s)$ and $\zeta(s,a)$) and their alternating equivalents $\eta(s)$ and $\eta(s,a)$. In particular, special cases involving improper integrals containing $\zeta(\sigma+it)$ and few other functions in the integrand are identified. Many of these integrals that do not appear in the literature, can be, and were, verified numerically. In one limit, the use of analytic continuation generates a family of improper integrals containing only the real and imaginary parts of $\zeta(\sigma+it)$ with and without simple trigonometric factors; the associated closed form contains an (unclassified) entity that has many of the attributes of an essential singularity. Consequently, this means that the associated integrals are indeterminate (i.e. non-single valued), so a new symbol is introduced to label the indeterminism. Much of this paper examines this singularity from several angles in order to resolve the associated ambiguities, before eventually showing how it blends into the classical study of functions of real and complex variables in an unusual manner. This is done by establishing a self-consistent way to remove the singularity and thereby evaluate new members of a family of integrals of general interest that contain $\zeta(s,a)$ and $\eta(s,a)$. Some implications are proposed.

\section{Introduction} \label{sec:Intro}

\subsection{Basics} \label{sec:Basics}
The use of the Mellin transform (\cite[Chapter 3]{Paris&Kaminski}) in the analysis of Riemann's function $\zeta(s)$ and its relatives, those being (with $s\equiv \sigma+it$) Hurwitz' function $\zeta(s,a)$ 

\begin{equation}
\zeta(s,a)\equiv \sum_{k=0}^{\infty}1/(k+a)^s\,,\hspace{15pt} \sigma>1,
\label{Hzeta}
\end{equation}

and its alternating counterpart $\eta(s,a)$, defined by
 
\begin{equation}
\eta(s,a)\equiv \sum_{k=0}^{\infty}(-1)^{k}/(k+a)^s,\hspace{15pt} \sigma>0\,\,,
\label{HzetaAlt}
\end{equation}

alongside special cases corresponding to Riemann's Zeta function itself
\begin{equation}
\zeta(s)\equiv\zeta(s,1)=\sum_{k=0}^{\infty}1/(k+1)^s,\hspace{15pt} \sigma>0,
\label{ZetaDef}
\end{equation}
Dirichlet's alternating equivalent
\begin{equation}
\eta(s)\equiv\eta(s,1)=\sum_{k=0}^{\infty}(-1)^{k}/(k+1)^s=(1-2^{1-s})\,\zeta(s),\hspace{15pt} \sigma>0,
\label{etLim1}
\end{equation}
plus the elementary identities
\begin{equation}
\displaystyle \zeta(s,1/2)=\sum _{k=0}^{\infty }   1/\left( k+1/2 \right) ^{s} = \left( {2}^{s}-1 \right) \zeta \left( s \right)\,, 
\label{Sumh}
\end{equation}
and
\begin{equation}
\eta(1)=\ln(2)\,
\label{eta11}
\end{equation}
have been long known and fruitful (\cite[Chapter 4]{Paris&Kaminski}, \cite[chapter 9]{IvicZ},\cite[Section 2.15]{Titch2}). Additionally, we note Nielsen's Beta function \cite{Nielsen} defined by

\begin{equation}
B(s)\equiv \eta(1,s)\,,
\label{NBeta}
\end{equation}

as well as the special case Dirichlet's function $\beta(s)$ 
\begin{equation}
\beta(s)\equiv \eta(s,1/2)/2^s\,
\label{Beta(s)}
\end{equation}

(which here differs from some definitions by the factor $2^s$). Therefore it is somewhat surprising to find that an elementary application of this transform unearths some significant and important properties pertaining to the evaluation of integrals containing these functions that appear to have been overlooked. 

\subsection{Motivation}
In Section (\ref{sec:Basic}), two simple, some might say ``trivial", theorems are derived, from which three corollaries are obtained by the simple application of the Mellin inverse and convolution transforms \cite{Paris&Kaminski}. At this point one might expect the author to isolate a few special cases of minor interest from these very simple identities and move on to the study of weightier matters elsewhere. However, great oaks from tiny acorns grow, and in this case it turns out that the extraction of special cases from simple identities is not always simple, and many of the identities that arise appear to be unknown in the literature, where many related evaluations are limited to integrals with finite limits - e.g. \cite{HuKimKim}, \cite{Esp&Moll2002} and \cite{Shpot&Paris}. In particular, one case uncovers a pathology that stands out from the usual classification of singularities found in many textbooks, lying somewhere between a functional discontinuity when viewed in the context of real variables, and an essential singularity when examined in the complex plane. It was the discovery of this pathology (cf. Section \ref{sec:PhiPiw0HR}) that diverted this paper from the original goal of listing a number of new and straightforward improper integrals involving the $\zeta(s),~\eta(s)$ and $\beta(s)$ family of functions and then moving on.

\subsection{Notation} \label{sec:Notation}
With respect to notation, in the entirety of this work, $j,k,m$ and $n$ are always non-negative integers except where noted; all other symbols are complex but are usually treated as real, except if noted. $M_{h}\equiv m+1/2$ and $M_{u}\equiv m+1$ are used for typographical brevity. $\psi^{(k)}(x)$ (i.e. polygamma) is the $k^{th}$ derivative of the digamma function $\psi(x)$, except for the first derivative written as $\psi^{\prime}(x)$. $\gamma$ is the Euler-Mascheroni constant and $\gamma(1)=-0.0728158...$ is the first Stieltjes constant in the expansion of $\zeta(s)|_{s=1}$ (see \cite[Eq. (25.2.5)]{NIST}). Throughout, the symbol $s$ may be referenced by its real and/or imaginary components, that is $s=\sigma+it$ where $\sigma,~t\in\Re$. Other symbols are defined either as they appear, or in Appendix \ref{sec:Defs}. Any sum whose lower limit exceeds its upper limit is zero and $:=$ indicates symbolic replacement. The real ($\Re$) and imaginary ($\Im$) components of a complex function $f(s)$ are denoted by subscripts $f_{R}(s)$ and $f_{I}(s)$ respectively. Complicated and lengthy identities are frequently quoted; these were obtained by either of the computer codes Mathematica \cite{Math} or Maple \cite{Maple}, which are cited by name at the corresponding points in the text.

\subsection{Outline} \label{sec:Outline}
The plan of the paper follows. Each of three corollaries derived in Section \ref{sec:Basic} is studied, and special cases are extracted in an individual section devoted to each corollary. In Section \ref{sec:Cor3a}, special cases are studied that involve an integrand containing the product of binary combinations of Hurwitz' function $\zeta(s,a)$ and/or its alternating equivalent $\eta(s,a)$ . Section \ref{sec:GenCon1} is devoted to a detailed study of Corollary (\ref{sec:Cor1}) that focusses on integrands containing the alternating function $\eta(s,a)$. Section \ref{sec:GenCor2} is devoted to special cases related to Corollary (\ref{sec:Cor2}) that deals with integrands containing the function $\zeta(s,a)$ directly. Section \ref{sec:Except_eta} examines an unusual pathology arising in Section \ref{sec:GenCon1} that manifests itself in the form of an ambiguous limit when some special cases of the the more general forms are evaluated. This means that there exist exceptional cases where the integrals are multi-valued and are therefore indeterminate at a point. In this Section it is shown how the pathology relates to a similar integral examined from the point of view of a real variable. Section \ref{sec:Except_zeta}  examines a similar pathology that arises in Section \ref{sec:GenCor2}. In both these Sections, it is shown how to bypass the pathology in different ways and arrive at a self-consistent method to evaluate the relevant integrals. Section \ref{sec:Pv} explores a different method of circumventing the pathology by evaluating related integrals using an analogue of the well-known contour deformation technique of mathematical and physical analysis.

The presentation of cases in each of these Sections is categorized according to the value of a real parameter $0\leq w\leq 1/2$ introduced in Section \ref{sec:Basic} that interpolates between $\zeta(s,a)$ and $\eta(s,a)$.  In a summary Section (\ref{sec:Summary}), I suggest possible implications and applications that could be fruitfully pursued. In Appendix \ref{sec:Defs} will be found accessory definitions, identities and lemmas used throughout this work. Appendix \ref{sec:LongEqs} records a number of lengthy identities that are the most general forms of the shorter, special-case identities employed in the main text, in order that the reader can reproduce the special results presented. That Appendix also lists a few numerical constants. For readers who may not be well-versed in the  evaluation of a singular integral by analytic continuation and/or regularization, the final Appendices \ref{sec:Analogue} and \ref{sec:Reg} review these analysis methods, by means of simple analogues extracted from the pages of mathematical history and modern quantum  physics.  

\subsection{Highlights} \label{Highlights}
Because many of the special cases are buried in a mass of derivation identities, listed below are a dozen of the more interesting results. See below for an explanation of the cases labelled by (*).

\item{From Section \ref{sec:Cor3a}}, Equation \eqref{Wb0}:
\begin{align} \nonumber
\int_{-\infty}^{\infty}&\frac{{| \zeta \! \left(\frac{1}{2}+i\,v \right)|}^{2} \left(\sqrt{2}\,\cos \! \left(v\,\ln \! \left(2\right)\right)-1\right)}{\cosh ^{2}\! \left(\pi \,v \right)}d v\\ \nonumber
&=-\frac{1}{\pi}\int_{1/2}^{1}\psi \! \left(v \right)^{2}d v+\frac{5\,\pi}{24}+\frac{1}{\pi}\left(\frac{1}{2}\ln^{2}\left(\pi \right)-\ln^{2} \left(2\right)-\frac{\gamma^{2}}{2}-\gamma \! \left(1\right)\right)\,. \tag{*}
\end{align}

\item{From Subsection (\ref{sec:SpCasesHR}}), Equations \eqref{Hiworh}, \eqref{Diff1}, \eqref{s01ma} and \eqref{Az01m} respectively
\begin{equation} \nonumber
\int_{0}^{\infty}\frac{\beta_{I} \! \left(\frac{1}{2}+i\,v \right) }{\cosh \! \left(\pi \,v \right)}\,v\,d v
 = \frac{\pi^{2}}{48}-\frac{\ln \! \left(2\right)}{4}\,,\tag{*}
\end{equation}

\item{}
\begin{equation} \nonumber
\int_{0}^{\infty}\eta_{I} \! \left(1/2+i\,v\right) v d v  = 0\,,
\end{equation}

\item{}
\begin{equation} \nonumber
\int_{0}^{\infty}\left(3-
{2  \sqrt{2}} \cos \! \left(v\ln \! \left(2\right)\right)\right)\,\zeta_{R} \! \left(1/2+i\,v \right)d v =0\,
\end{equation}
and
\item{}
\begin{equation} \nonumber
\int_{0}^{\infty}\left(3-2\,\sqrt{2}\,\cos \! \left(v\,\ln \! \left(2\right)\right)\right) \zeta_{I} \! \left(1/2+i\,v \right)\,\tanh \! \left(\pi \,v \right)d v
 = 3\,\ln \! \left(2\right)-2\,.
\end{equation}

\item{ From Subsection \ref{sec:Pv}, Equation \eqref{HAB4}}
\begin{equation} \nonumber
\int_{0}^{\infty}v^{4}\,\eta_{I} \! \left(1/2+i\,v \right)\,\tanh \! \left(\pi \,v \right)d v
 = 
-\frac{\ln \! \left(2\right)}{16}-\frac{5\,\pi^{2}}{12}-\frac{87\,\zeta \! \left(3\right)}{4}-\frac{7\,\pi^{4}}{15}-\frac{45\,\zeta \! \left(5\right)}{2}\,.
\end{equation}

\item{From Section \ref{sec:GenCor2}, Equations \eqref{Rz}}, \eqref{Ip4a}, \eqref{CR3}, \eqref{CR5}, \eqref{Wdm1} and \eqref{Oa4b} respectively:
\item{}
\begin{equation} \nonumber
\int_{0}^{\infty}\frac{\zeta_{R} \! \left(\frac{1}{2}+i\,v \right)}{\cosh \! \left(\pi \,v \right)}d v
 = \gamma-1\,,\tag{*}
\end{equation}

\item{}
\begin{equation} \nonumber
\int_{0}^{\infty}\frac{ \zeta_{I} \! \left(\frac{1}{2}+i\,v \right)}{\cosh \! \left(\pi \,v \right)}\,v\,d v
 = -\frac{3}{2}+\frac{\pi^{2}}{6}-\frac{\gamma}{2}\,,\tag{*}
\end{equation}

\item{}
\begin{equation} \nonumber
\int_{0}^{\infty}\sin \! \left(v\,\ln \! \left(2\right)\right) \zeta_{I} \! \left(1/2+i\,v \right)d v
 = \frac{\pi}{2\,\sqrt{2}} \,,
\end{equation}

\item{}
\begin{equation} \nonumber
\int_{0}^{\infty}\zeta_{R} \! \left(1/2+iv \right)d v =\, -\pi\,,
\end{equation}

\item{}
\begin{equation} \nonumber
\int_{0}^{\infty} \zeta_{I} \! \left(1/2+i\,v \right)\tanh \! \left(\pi \,v \right)d v
 = -\gamma\,
\end{equation}
and
\item{}
\begin{equation} \nonumber
\int_{0}^{\infty}\left(2\,\cos \! \left(v\,\ln \! \left(2\right)\right)-\sqrt{2}\,\right) \zeta_{I} \! \left(1/2+i\,v \right)\,\tanh \! \left(\pi \,v \right)d v
 \,{\underset{\looparrowleft}=} 
\sqrt{2} \left(1-\frac{3\,\sqrt{2}\,\ln \! \left(2\right)}{2}-\frac{\gamma}{2}\right)\,.
\end{equation}

The integrals in the identities listed above that are tagged by (*) are obviously convergent because of the exponentially increasing denominator and have all been verified numerically; others are true in the sense of analytic continuation (i.e. regularization), recognizing that all the zeta-family of  functions are known to both fluctuate in sign and possess unknown (or conjectured) asymptotic behaviour in the argument range of interest here (see Lindel{\"o}f's Hypothesis, Definition \ref{sec:Lind} in Appendix \ref{sec:Defs}) and so the convergence of those not so-tagged is a matter of conjecture. For readers unfamiliar with regularization by analytic continuation, a review is given in Appendix C. In Subsection (\ref{sec:PhiPiwh}) an attempt is made to numerically evaluate a representative candidate integral with limited success using Ces{\`a}ro summation \cite[Appendix A]{MilHughMMM2025}. The symbol ${\underset{\looparrowleft}=}$ is introduced in Subsection (\ref{sec:PhiPiw0HR}).

\section{Basic Identities} \label{sec:Basic}

It is convenient here to introduce the real variable $w>0$ and the complex variables $a,~b$ and $z$.

\begin{theorem}{\bf The Mellin transform of Hurwitz' alternating function $\eta(s,a)$}:\newline \label{sec:Thm1}
If $\Re(z)>0$, $w>0$, $w\in\Re$ and $\Re(b+1-z)>0$ then
\begin{equation}
\int_{0}^{\infty}v^{z -1}\,\eta \! \left(b +1, v +w +1/2\right)d v
 = 
\frac{\Gamma \! \left(z \right)  \Gamma \! \left(b +1-z \right)}{\Gamma \! \left(b +1\right)}\eta \! \left(b +1-z , w +{1}/{2}\right).
\label{Cry1a}
\end{equation}
\begin{proof}
\begin{align} \nonumber
&\int_{0}^{\infty}v^{z -1}\,\eta \! \left(b +1, v +w +1/2\right)d v
 = 
 \overset{\infty}{\underset{j =0}{\sum}}\! \left(-1\right)^{j}\int_{0}^{\infty}\frac{v^{z -1}}{\left(j +v +w +\frac{1}{2}\right)^{b +1}}d v \\
& = 
\frac{ \Gamma \! \left(b +1-z \right) \Gamma \! \left(z \right)}{\Gamma \! \left(b +1\right)}\overset{\infty}{\underset{j =0}{\sum}}\! \frac{\left(-1\right)^{j}}{\left(j +w +\frac{1}{2}\right)^{b +1-z}}
=\frac{\Gamma \! \left(z \right)  \Gamma \! \left(b +1-z \right)}{\Gamma \! \left(b +1\right)}\eta \! \left(b +1-z , w +{1}/{2}\right).
\label{Proof1}
\end{align}

The first equality is obtained by applying \eqref{HzetaAlt} followed by the transposition of the summation and integration operators, validated because both are convergent under the conditions specified; the second equality is the result of an elementary integration \cite[Eq. (3.251(11)]{G&R}, and the third again reflects the definition \eqref{HzetaAlt}, giving the final result Theorem \eqref{sec:Thm1}.  {\bf QED}

\end{proof}
\end{theorem}

\begin{theorem} {\bf The Mellin transform of Hurwitz' function $\zeta(s,a)$}:\newline \label{sec:Thm2}
If $\Re(z)>0$, $w>0$, $w\in\Re$, $\Re(b)>0$ and $\Re(b-z)>0$ then
\begin{equation}
\int_{0}^{\infty}v^{z -1}\,\zeta\! \left(b +1, v +w +{1}/{2}\right)d v
 = 
\frac{\Gamma \! \left(z \right)  \Gamma \! \left(b +1-z \right)}{\Gamma \! \left(b +1\right)}\zeta\! \left(b +1-z , w +{1}/{2}\right).
\label{Cry1A}
\end{equation}
\begin{proof}
The proof is exactly the same as for Theorem \ref{sec:Thm1}, except that the bound on the variable $b$ differs.
\end{proof}
\end{theorem}

\begin{rem}
\leavevmode
\begin{itemize}
\item \eqref{Cry1A} is listed in \cite[Eq. 2.3.1(1)]{prudnikov};
\item 
in the above, ${\Gamma \! \left(z \right) \Gamma \! \left(b +1-z \right)}/{\Gamma \! \left(b +1\right)}=\beta(z,b+1-z)$ where $\beta(x_{1},x_{2})$ is Euler's beta integral \cite[Eq. (5.12.1)]{NIST}.
\end{itemize}
\end{rem}

\begin{Coroly} \label{sec:Cor1}
\begin{align} \nonumber
&\frac{1}{2\,\pi\,i \,\Gamma \! \left(b +1\right)} \int_{\sigma -i\,\infty }^{\sigma +i\,\infty }a^{-v}\,\Gamma \! \left(v \right) \eta \! \left(b +1-v , w +1/2\right) \Gamma \! \left(b +1-v \right)d v \\ \nonumber \\
& = \eta \! \left(b +1, a +w +1/2\right).
\label{M1ng1}
\end{align}
\begin{proof}
This is the inverse Mellin transform of Theorem \ref{sec:Thm1} provided that $0<\sigma<b+1$ and $-\pi<\arg(a)<\pi$ (see \cite{Paris&Kaminski}).
\end{proof}

\end{Coroly}

\begin{Coroly} \label{sec:Cor2}
\begin{align} \nonumber
&\frac{1}{2\pi\,i\,\Gamma \! \left(b +1\right) } \int_{\sigma -i\,\infty}^{\sigma +i\,\infty }a^{-v}\,\Gamma \! \left(v \right) \zeta\! \left(b +1-v , w +1/2\right) \Gamma \! \left(b +1-v \right)d v \\\nonumber\\
&=\zeta\! \left(b +1, a +w +1/2\right). 
\label{M1ng2}
\end{align}
\begin{proof}
This is the inverse Mellin transform of Theorem \ref{sec:Thm2} provided that $0<\sigma<b$ and $-\pi<\arg(a)<\pi$ (see \cite{Paris&Kaminski}).
\end{proof}
\end{Coroly}


\begin{Coroly} \label{Cor3} 
If any of $Z_{1}(s,a)$ and $Z_{2}(s,a)$ represent one of $\zeta(s,a)$ or $\eta(s,a)$, then
\begin{align} \nonumber
&\int_{0}^{\infty}v^{z -1}\,{Z_{1}\! \left(b +1, v +w +1/2\right)}\,Z_{2}\! \left(b +1, v +w +1/2\right)d v\\ \nonumber &= 
\frac{1}{2\,\pi\,i\,\Gamma \left(b +1\right)^{2}}  \int_{\sigma -i\,\infty }^{\sigma +i\,\infty }\Gamma \left(v \right) \Gamma \left(z -v \right)\Gamma \left(b +1-v \right) \Gamma \left(b +1-z +v \right)
\\&\times
  Z_{1}\left(b +1-v , w+1/2\right) Z_{2}\left(b +1-z +v , w +1/2\right) d v.
\label{Ev5}
\end{align}
\begin{proof}
This is the Mellin convolution transform of Theorem \ref{sec:Thm2} (see \cite[Eq. 3.1.14]{Paris&Kaminski}) provided that $~0<w\in\Re$, $0<\sigma<1,~\Re(z)>0$ and $\sigma<\Re(b)$.
\end{proof}
\end{Coroly}

\begin{rem}
\leavevmode
\begin{itemize}

\item{For those identities amenable to such analysis, all of the above can be (and were) numerically verified in the range of parameters specified, utilizing the identity \eqref{EtaZeta} as the basis for numerical verification. This assumes that the numerical methods for $\zeta(s,a)$ embedded in the two computer codes used - Mathematica \cite{Math} and Maple \cite{Maple} - are numerically robust and accurate.}

\item{If $w=0$ any of the functions $\eta\,(s,w+1/2)$ in the above reduce to Dirichlet's Beta function - see \eqref{Beta(s)} and any of the functions $\zeta\,(s,w+1/2)$ are proportional to $\zeta(s)$ - see \eqref{Sumh}.}

\item{If $w=1/2$, any of the functions $\eta\,(s,w+1/2)$ in the above reduce to Hurwitz' alternating zeta function (see \eqref{etLim1}) and hence are proportional to Riemann's zeta function, and any of the functions $\zeta(s,w+1/2)$ reduce to Riemann's zeta function itself (see \eqref{ZetaDef}).}

\item{The variable $w$ in the range $0\leq w \leq 1/2$ therefore interpolates between $\eta(s,a)$ and $\zeta(s,a)$.}

\end{itemize}
\end{rem}

\section{Special Cases related to  Corollary \ref{Cor3}} \label{sec:Cor3a}
\subsection{$Z_{1}(s,a)=Z_{2}(s,a)=\zeta(s,a)$}
Since the left-hand side of \eqref{Ev5} is independent of the parameter $\sigma$ in the range $0<\sigma<1$, let $\sigma=1/2$ along with $z=1,~w=0, ~ \Re(b)>\sigma$, so that \eqref{Ev5} becomes

\begin{equation}
\int_{0}^{\infty}{\zeta\! \left(b +1, v +1/2\right)}^{\,2}d v
 = 
\frac{1}{2\,\Gamma \! \left(b +1\right)^{2}}\int_{-\infty}^{\infty}\frac{{| \zeta \left( b +\frac{1}{2}+i\,v , \frac{1}{2}\right)|}^{2}\,{| \Gamma \left(b +\frac{1}{2}+i\,v \right)|}^{2}}{\cosh \left(\pi \,v \right)}d v\,. 
\label{Cv2}
\end{equation}
Because the right-hand side of \eqref{Cv2} is effectively the contour integral \eqref{Ev5} following a change of variables, as the value of $b$ varies, the integrand on the right-hand side of the contour integral \eqref{Ev5} develops poles whose residues must be added or subtracted as they cross the original contour of integration. Thus, with the help of \eqref{Sumh} along with the more restrictive requirement that $ b\in \Re$, \eqref{Cv2} becomes
\begin{align} \nonumber
&\int_{0}^{\infty}{\zeta\! \left(b +1, v +1/2\right)}^{\,2}d v
 = \mathfrak{X}(b,1/2,0)+
\frac{1}{2\,\Gamma \! \left(b +1\right)^{2}}\\
&\times \int_{-\infty}^{\infty}\frac{\left(1+2^{2\,b +1}-2^{\frac{3}{2}+b}\,\cos \left(v\,\ln \left(2\right)\right)\right) {| \zeta \left(b +\frac{1}{2}+i\,v \right)|}^{2}\,{| \Gamma \left(b +\frac{1}{2}+i\,v \right)|}^{2}}{\cosh \left(\pi \,v \right)}d v\,.
\label{Cv2a}
\end{align}
The following summarizes those residue terms that must be included in \eqref{Cv2a}, for general values of $0<\sigma<1$ and $-1<b<1$ with $b\in\Re$:

\begin{align} \label{ResX}
\mathfrak{X}(b,\sigma,w)\equiv\left\{\begin{array}{cc}
{\displaystyle \frac{ \pi \,\zeta \left(2\,b , w +\frac{1}{2}\right) \Gamma \! \left(2\,b \right)}{\Gamma(b+1)^{2}\,\sin(\pi\,b)}} & b <\sigma \mbox{ and/or } \sigma<1-b  
\\\\
{\displaystyle \frac{ \pi\,w \,\zeta \left(2\,b +1, w +\frac{1}{2}\right) \Gamma \! \left(2\,b +1\right)}{\Gamma(b+1)^{2}\,\sin(\pi\,b)} }  & b +1<\sigma \mbox{ and/or } \sigma<-b 
\\\\
{\displaystyle \frac{\pi\,\zeta\left(2\,b +2, w +\frac{1}{2}\right) \left(12\,w^{2}-1\right) \Gamma \left(2\,b +2\right)}{24\,\Gamma(b+1)^{2}\,\sin(\pi\,b)}} & \sigma <-b -1 
\\\\
 0 & \mathit{otherwise.}  
\end{array}\right.
\end{align}

With reference to \eqref{Zasyv}, the left-hand side of \eqref{Cv2a} converges if $b>1/2$ because the integrand asymptotically approaches $v^{-2b\,}$; the right-hand side converges because $|\Gamma(b+1/2+iv)|<|\Gamma(b+1/2)|$ (see \cite[Eq. (5.6.6)]{NIST}) and the denominator increases exponentially. Similarly, let $w=1/2$ and, again with $b>1/2\,$, we arrive at
\begin{align} \nonumber
\int_{0}^{\infty}{\zeta\! \left(b +1, v +1\right)}^{2}&d v
 = \mathfrak{X}(b,1/2,1/2)\\
 &+
\frac{1}{2\,\Gamma \! \left(b +1\right)^{2}}\int_{-\infty}^{\infty}\frac{{| \zeta \left(b +\frac{1}{2}+i\,v \right)|}^{2}\,{| \Gamma \left(b +\frac{1}{2}+i\,v \right)|}^{2}}{\cosh \left(\pi \,v \right)}d v
\label{Cv1b1}
\end{align}
both sides of which again clearly converge for the same reasons. Subtracting \eqref{Cv1b1} from \eqref{Cv2a} allows us to extend the range of validity with respect to the variable $b$. Application of a simple change of variables on the left-hand side of each, leads to the following identity: 

\begin{align} \nonumber
\int_{1/2}^{1}&{\zeta\! \left(b +1, v\right)}^{2}d v =\mathfrak{X}(b,1/2,0)-\mathfrak{X}(b,1/2,1/2)\\
&+ 
\frac{2^{2\,b}}{\Gamma \! \left(b +1\right)^{2}}\int_{-\infty}^{\infty}\frac{{| \zeta \left(b +\frac{1}{2}+i\,v \right)|}^{2}\,{| \Gamma \left(b +\frac{1}{2}+i\,v \right)|}^{2} \left(1-2^{1/2-b}\,\cos \left(v\,\ln \left(2\right)\right)\right)}{\cosh \left(\pi \,v \right)}d v\,.
\label{W1m2}
\end{align}

Notice that the left-hand side converges for $b\neq 0$, whereas the left-side integrals in \eqref{Cv2a} and \eqref{Cv1b1} formally diverge at their upper limit if $b<1/2$. As a simple example, since the $v=0$ limit of the integrand of the right-hand side integral with $b=1/2$ equals $\ln^2(2)/4/\pi$, both sides converge, giving the special case

\begin{equation}
\int_{1/2}^{1}{\zeta \left(3/2, v\right)}^{2}d v
 = 
16 \int_{-\infty}^{\infty}\frac{{| \zeta \! \left(1+i\,v \right)|}^{2}\, \left(1-\cos \! \left(v\,\ln \! \left(2\right)\right)\right)}{\sinh \! \left(2\,\pi \,v \right)}\,v\,d v +8\,\ln \! \left(2\right)
\label{W1m2b}
\end{equation}
using \eqref{Nist5p4p3} in the simplification process, a result that was numerically verified. The limiting case $b=0$ is more interesting. From \eqref{W1m2} and \eqref{Nist5p4p4} we have
\begin{align} \nonumber
\int_{-\infty}^{\infty}&\frac{{| \zeta \! \left(\frac{1}{2}+i\,v \right)|}^{2} \left(\sqrt{2}\,\cos \! \left(v\,\ln \! \left(2\right)\right)-1\right)}{\cosh ^{2}\! \left(\pi \,v \right)}d v\\ \nonumber
& = 
\underset{b \rightarrow 0}{\mathrm{lim}}\! \left(-\frac{1}{\pi}\int_{\frac{1}{2}}^{1}{\zeta\! \left(b +1, v\right)}^{2}d v+\frac{4\,\zeta \! \left(2\,b \right) \left(2^{2\,b-1}-1\right) \Gamma \! \left(2\,b \right)}{\Gamma \! \left(b +1\right)^{2}\,\sin \! \left(\pi \,b \right)}\right)\\ \nonumber
&=
\underset{b \rightarrow 0}{\mathrm{lim}}\! \left(-\frac{1}{\pi}\int_{1/2}^{1}\left(\zeta\! \left(b +1, v\right)-\frac{1}{b}\right)^{2}d v+\frac{1}{\pi}\int_{1/2}^{1}\left(\frac{1}{b^{2}}-\frac{2}{b}\,\zeta\left(b +1, v\right)\right)d v\right. \\  \nonumber
&\left. \hspace{30pt}+\frac{1}{\pi}\left(\frac{\pi^{2}}{4}-\gamma^{2}+\ln^{2} \left(\pi \right)-2\,\ln^{2} \left(2\right)-2\,\gamma \! \left(1\right)\right)+\frac{\ln \! \left(\pi \right)}{\pi \,b}+\frac{1}{2\,\pi \,b^{2}}\right)\\ \nonumber
&=\underset{b \rightarrow 0}{\mathrm{lim}}\! \left(-\frac{1}{\pi}\int_{1/2}^{1}\psi \! \left(v \right)^{2}d v+\frac{\frac{1}{2}-2 \left(2^{b}-1\right) \zeta \! \left(b \right)+2\,\zeta \! \left(b \right)}{\pi \,b^{2}} \right. \\ \nonumber
&\left. \hspace{30pt}+\frac{1}{\pi}\left(+\frac{1}{2\,b^{2}}+\frac{\ln \left(\pi \right)}{b}+\frac{\pi^{2}}{4}-\gamma^{2}+\ln \! \left(\pi \right)^{2}-2\,\ln \! \left(2\right)^{2}-2\,\gamma \! \left(1\right)\right)\right)\\
&=-\frac{1}{\pi}H+\frac{5\,\pi}{24}+\frac{1}{\pi}\left(\frac{1}{2}\ln^{2}\left(\pi \right)-\ln^{2} \left(2\right)-\frac{\gamma^{2}}{2}-\gamma \! \left(1\right)\right)=0.447498183476330\dots\,,
\label{Wb0}
\end{align}
where
\begin{equation}
H\equiv\int_{1/2}^{1}\psi^{2} \! \left(v \right)d v\,.
\label{Hdef}
\end{equation}
In the above, the first equality reflects the addition of two different residue pairs from \eqref{ResX}, the second equality presents the evaluation of the residue terms in the desired limit as well as the addition and subtraction of diverging terms, the third equality makes use of \eqref{Zlim0} as well as the evaluation of the second integral (see \eqref{Zident}) leading to a final result at the fourth equality, an identity verified numerically. The constant $H$ is given numerically in Table \ref{sec:Consts}; it does not appear to be reducible to fundamental constants, although some literature revolves about the function $\psi^{2}(v)$ and its appearance in an integral (see \cite{connon2012integral}, \cite[Section 4]{CoppoAndCandel} and \cite[Section 20]{Nielsen}). Further, an integral of the form defined by the left-hand side of \eqref{W1m2} between the limits $(0,1)$ is evaluated in \cite[Eq. (3.1)]{Esp&Moll2002}.

\subsection{$Z_{1}(s,a)=Z_{2}(s,a)=\eta(s,a)$} \label{sec:Cor3b}


Similar to the special cases presented previously, the general case here with $z=1,~b>-1/2$ is
\begin{align} \nonumber
&\hspace{-60pt}\int_{0}^{\infty}\eta \! \left(b +1, v +w +1/2\right)^{\,2}d v
 = \mathfrak{Y}(b,\sigma,w)\\ + \frac{1 }{2\,i\,\Gamma \! \left(b +1\right)^{2}} 
&\int_{\sigma -i\,\infty }^{\sigma +i\,\infty }\frac{\eta \left(b +1-v , w +\frac{1}{2}\right) \Gamma \left(b +1-v \right) \eta \left(b +v , w +\frac{1}{2}\right) \Gamma \left(b +v \right)}{\sin \left(\pi \,v \right)}\,d v 
\label{Ev1}
\end{align}
where certain residues need to be added as $b$ varies and the contour moves, giving:
\begin{align}
\mathfrak{Y}(b,\sigma,w)\equiv
\left\{\begin{array}{cc}
{\displaystyle -\frac{\pi \,\eta \left(2\,b +1, w +\frac{1}{2}\right) \Gamma \left(2\,b +1\right)}{2\,\sin \left(\pi \,b \right) \Gamma \left(b +1\right)^{2}}} & \sigma <-b \mbox{ and/or } \sigma>b+1 
\\\\
 0 & \mathit{otherwise,}  
\end{array}\right.
\label{Res}
\end{align}
accordingly to the value of $b$ relative to  $\sigma$.
\subsubsection{Special case $w=0$} \label{sec:w0Cor3p2}
In the case that $b=0,~w=0~\mbox{and}~\sigma=1/2$ from \eqref{Ev1} we obtain
\begin{equation}
\int_{-\infty}^{\infty}\frac{{| \beta \! \left(\frac{1}{2}+i\,v \right)|}^{2}}{\cosh ^{2}\! \left(\pi \,v \right)}d v
 = 
\frac{1}{\pi}\int_{1/2}^{\infty}\eta \! \left(1, v\right)^{2}d v=\frac{1}{\pi}\int_{1/2}^{\infty}B(v)^{2}d v
\label{Xs2}
\end{equation}
a relation between a weighted integration of Dirichlet's Beta function on the critical line and a special case integral of Hurwitz' alternating function on the real line, otherwise referred to as Nielsen's beta function - see \eqref{NBeta}, \cite[Eq. 8.372]{G&R}, \cite[Section 8]{connon2012integral} . 

\subsubsection{Special case $w=1/2$} \label{sec:whCor3p2}
Similarly, with $w=1/2$ we find
\begin{equation}
\int_{-\infty}^{\infty}\frac{\left(3-2\,\sqrt{2}\,\cos \! \left(v\,\ln \! \left(2\right)\right)\right) {| \zeta \! \left(\frac{1}{2}+i\,v \right)|}^{2}}{\cosh^{2} \! \left(\pi \,v \right)}d v
 = 
\frac{2 }{\pi}\int_{1}^{\infty}\eta \! \left(1, v\right)^{2}d v \,.
\label{Ev2a}
\end{equation}\newline
The above results yield scope for interesting variations. For example, combining \eqref{Wb0} and \eqref{Ev2a} produces the transformation 
\begin{align} 
\int_{-\infty}^{\infty}\frac{{| \zeta \! \left(\frac{1}{2}+i\,v \right)|}^{2}}{\cosh^{2} \! \left(\pi \,v \right)}d v
 = \frac{2}{\pi}\left(\int_{1}^{\infty}\eta \! \left(1, v\right)^{2}d v
-H\right) +\frac{5\,\pi}{12}
-\frac{2\,\ln \! \left(2\right)^{2}-\ln \! \left(\pi \right)^{2}+\gamma^{2}+2\,\gamma \! \left(1\right)}{\pi}\,.
\label{Xs0}
\end{align}
Alternatively, combining the same two identities with different weighting yields
\begin{align} \nonumber
\int_{-\infty}^{\infty}\frac{{| \zeta \! \left(\frac{1}{2}+i\,v \right)|}^{2}\,\cos \! \left(v\,\ln \! \left(2\right)\right)}{\cosh^{2} \! \left(\pi \,v \right)}&d v
 = 
-\frac{3\,\sqrt{2} }{2\,\pi}\,H +\frac{\sqrt{2}}{\pi} \int_{1}^{\infty}\eta \! \left(1, v\right)^{2}d v -\frac{3\,\sqrt{2}\,\ln^{2} \! \left(2\right)}{2\,\pi}\\
&+\frac{3\,\sqrt{2}\,\ln^{2} \! \left(\pi \right)}{4\,\pi}-\frac{3\,\sqrt{2}\,\gamma \! \left(1\right)}{2\,\pi}+\frac{5\,\sqrt{2}\,\pi}{16}-\frac{3\,\sqrt{2}\,\gamma^{2}}{4\,\pi}\,.
\label{Xs1}
\end{align}
All the above were verified numerically.
\subsection{$Z_{1}(s,a)=\zeta(s,a),~Z_{2}(s,a)=\eta(s,a)$} \label{sec:Cor3c}


For $\sigma<z$ and $0<\sigma<1$, the general case \eqref{Ev5} gives
\begin{align} \nonumber
&\int_{0}^{\infty}v^{z -1}\,\zeta \left(b +1, v +w +\frac{1}{2}\right) \eta \! \left(b +1, v +w +\frac{1}{2}\right)d v  
 = \mathfrak{Z}(z,b,\sigma,w)+\frac{1}{2\,\pi \,\Gamma \! \left(b +1\right)^{2}} \\ \nonumber 
& \times\int_{-\infty}^{\infty}\Gamma \! \left(\sigma+i\,t  \right) \zeta \left(b-\sigma-i\,t   +1, w +1/2\right) \Gamma \! \left(b-\sigma-i\,t   +1\right) \Gamma \! \left(z-\sigma -i\,t  \right)\\
&\times \eta \! \left(\sigma+i\,t +b  -z +1, w +\frac{1}{2}\right) \Gamma \! \left( b+\sigma +i\,t -z +1\right)d t
\label{CbA1}
\end{align}

where, as in the previous cases, the following general residues need to be included according to the value of the parameter $b$ relative to $\sigma$. With $z>\sigma$ we have
\begin{align}
\mathfrak{Z}(z,b,\sigma,w)=\left\{\begin{array}{cc}
{\displaystyle \frac{\Gamma \! \left(b \right) \Gamma \! \left(z -b \right) \eta \! \left(2\,b +1-z , w +\frac{1}{2}\right) \Gamma \! \left(2\,b +1-z \right)}{\Gamma \! \left(b +1\right)^{2}}}  & b <\sigma \\\\
%
%
{\displaystyle \frac{\Gamma \! \left(-b -1+z \right) \zeta\left(2\,b +2-z , w +\frac{1}{2}\right) \Gamma \! \left(2\,b +2-z \right)}{2\,\Gamma \! \left(b +1\right)}} & \sigma <-b \\\\
%
%
{\displaystyle -\frac{w\,\Gamma \! \left(-b -1+z \right) \eta \! \left(2\,b +2-z , w +\frac{1}{2}\right) \Gamma \! \left(2\,b +2-z \right)}{\Gamma \! \left(b +1\right)}}& b +1<\sigma  \\\\
 0 & \mathit{otherwise}  \,.
\end{array}\right.
\label{Mres}
\end{align}

\subsubsection{Special case $w=0$}  \label{sec:w0Cor3p3}

There do not appear to be any interesting cases. 

\subsubsection{Special case $w=1/2$}  \label{sec:whChr3p2}

Of special interest is the case $b=0,~z=1$ and $w=1/2$, so, first set $\sigma=b+1/2$, and consider the limit $b\rightarrow 0$ by subtracting a term $1/b$ from the integrand of the left-hand side of \eqref{CbA1}, then adding and identifying the equivalent integral using \eqref{Cry1a}, to obtain
\begin{align} \nonumber
&\underset{b \rightarrow 0}{\mathrm{lim}}\! \left(\int_{0}^{\infty}v^{z -1}\,\psi \! \left(v +1\right) \eta \! \left(b +1, v +1\right)d v \right.\\& \left. \nonumber +\frac{\Gamma \! \left(b \right) \Gamma \! \left(z -b \right) \eta \! \left(2\,b +1-z , 1\right) \Gamma \! \left(2\,b +1-z \right)}{\Gamma \! \left(b +1\right)^{2}}
-\frac{\Gamma \! \left(z \right) \eta \! \left(b +1-z , 1\right) \Gamma \! \left(b +1-z \right)}{\Gamma \! \left(b +1\right) b}\right)
\\& = 
\frac{\pi}{2}\int_{-\infty}^{\infty}\frac{\,\zeta \left(\frac{1}{2}-i\,t \right) \eta \left(\frac{3}{2}-z +i\,t \right)d t}{\cosh \left(\pi \,t \right) \left(i\,\sin \left(\pi \,z \right) \sinh \left(\pi \,t \right)+\cos \left(\pi z \right) \cosh \left(\pi \,t \right)\right)},
\label{Sc2L}
\end{align}
where \eqref{Zlim0} has also been employed and the appropriate residue from \eqref{Mres} subtracted and evaluated using \eqref{Cry1a}. Evaluating the limit $b\rightarrow 0$ is now straightforward, as is the evaluation of the limit $z\rightarrow 1$ if a term $\ln(v)$ is subtracted from the integrand so that it will converge at its upper limit, and an equivalent term added and identified using \eqref{T2p2a}, to eventually arrive at
\begin{equation}
\frac{2}{\pi} \int_{0}^{\infty}\left(\psi \! \left(v +1\right)-\ln \! \left(v \right)\right) \eta \! \left(1, v +1\right)d v 
 = 
\int_{-\infty}^{\infty}\frac{\left(\sqrt{2}\,\cos \! \left(t\,\ln \! \left(2\right)\right)-1\right) {| \zeta \! \left(\frac{1}{2}+i\,t \right)|}^{2}}{\cosh^{2} \! \left(\pi \,t \right)}d t\,.
\label{Sc2Ld}
\end{equation}
\begin{rem} Comparing \eqref{Sc2Ld} and \eqref{Wb0} identifies the left-hand side of \eqref{Sc2Ld} in terms of fundamental constants.
\end{rem}
\section{Special Cases: Corollary \ref{sec:Cor1}} \label{sec:GenCon1}

{\bf For most of the remainder of this work, let $a=re^{i\phi}$, with $0<\{r,w\} \in \Re$ and $-\pi<\phi<\pi,\mbox{ } \phi \in \Re$  and $m=0,1,2\dots$. In Sections \ref{sec:Except_eta} and \ref{sec:Except_zeta}, we will consider a subset consisting of exceptional cases corresponding to the limit $\phi\rightarrow \pm\pi.$}

\subsection{General Identities - A Summary}

In the following Subsection (\ref{sec:Analytic}), we first consider \eqref{M1ng1} (i.e. Corollary \ref{sec:Cor1}) in general, by decomposing it into its real and imaginary parts, and determine the analytic structure of each, at which time special cases are identified. This is followed by Subsection (\ref{sec:SpCasesHR}) that considers special case identities corresponding to $\phi=0$. Another Subsection (\ref{sec:PhiPiallw}) follows where, for general values of the parameter $w$, the case $\phi=\pm \pi$  is considered. This leads to the discovery of exceptional cases identified by a special combination of parameters; the analysis of these cases is deferred to later Sections \ref{sec:Except_eta} and \ref{sec:Pv} that are devoted to the development of methods required to deal with an associated pathology.
\subsection{Analytic structure and exceptional cases} \label{sec:Analytic}


Since $0<\sigma<1$, the identity \eqref{M1ng1} with $|\arg(a)|<\pi$ can be cast into a more familiar form


\begin{equation}
\frac{1}{2\,\Gamma \! \left(b +1\right)}\int_{-\infty}^{\infty}\frac{a^{\sigma-i\,v  -1}\,\eta \left( b +\sigma-i\,v , w +\frac{1}{2}\right) \Gamma \left( b +\sigma-i\,v \right)}{\sin \left(\pi  \left(\sigma-i\,v  \right)\right) \Gamma \left( \sigma -i\,v\right)}d v
 = \eta \! \left(b +1, a +w +1/2\right)
\label{C1}
\end{equation}

by replacing $\sigma:=1-\sigma$ and employing the Gamma function reflection property. Notice that the right-hand side of \eqref{C1} is formally independent of the variable $\sigma$ in the range specified. Limiting our interest to the case $b=0$, which satisfies the condition \eqref{M1ng1}, \eqref{C1} can be rewritten

\begin{equation}
 {{\int}}_{\!\!\!-\infty}^{\infty}\frac{r^{{-i} v} {\mathrm e}^{\phi  v} \eta \! \left(\sigma -{i} v ,w +\frac{1}{2}\right)}{\sin \! \left(\pi  \left(\mathrm{i} v -\sigma \right)\right)}{d}v =-2 \eta \! \left(1,1/2+r \,{\mathrm e}^{i \phi}+w \right) {\mathrm e}^{{i} \phi  \left(1-\sigma \right)} r^{1-\sigma}\,.
\label{C1c}
\end{equation}

It is now convenient to split \eqref{C1c} into its real and imaginary parts, yielding, after converting the integration range to $(0,\infty)$ and utilizing the convergent series representation \eqref{HzetaAlt}, a substantial pair of equations, which, for the record, are presented in their full form as \eqref{C2R} and \eqref{C2I} in Appendix \ref{sec:LongEqs}. For simplicity, let $\sigma=1/2$ giving the (slightly) shorter forms


\begin{align} \nonumber
&\int_{0}^{\infty}\Re \! \left(r^{i\,v}\,\eta \! \left(\frac{1}{2}+i\,v , w +\frac{1}{2}\right)\right) \frac{ \cosh \! \left(\phi \,v \right)}{\cosh \! \left(\pi \,v \right)}d v \\ \nonumber
& = 
\int_{0}^{\mathit{\infty}}\left({\sin \! \left(v\,\ln \! \left(r \right)\right) \eta_{I} \! \left(\frac{1}{2}+i\,v , w +\frac{1}{2}\right)}{}-{\cos \! \left(v\,\ln \! \left(r \right)\right) \eta_{R} \! \left(\frac{1}{2}+i\,v , w +\frac{1}{2}\right)}{}\right) \frac{\cosh \! \left(\phi \,v \right)}{\cosh \! \left(\pi \,v \right)}\,d v \\ \nonumber
 &= 
\sqrt{r} \left(-r\,\sin \! \left(\phi \right)\,\sin \! \left(\frac{\phi}{2}\right) \overset{\infty}{\underset{j =0}{\sum}}\, \frac{\left(-1\right)^{j}}{2\,r \left(j +w +\frac{1}{2}\right) \cos \! \left(\phi \right)+r^{2}+\left(j +w +\frac{1}{2}\right)^{2}} \right.\\
&\left. \hspace{20pt} -\cos \! \left(\frac{\phi}{2}\right)\overset{\infty}{\underset{j =0}{\sum}}\, \frac{\left(-1\right)^{j} \left(r\,\cos \! \left(\phi \right)+j +w +\frac{1}{2}\right)}{2\,r \left(j +w +\frac{1}{2}\right) \cos \! \left(\phi \right)+r^{2}+\left(j +w +\frac{1}{2}\right)^{2}}\right) 
\label{HR}
\end{align}

and
\begin{align} \nonumber
\int_{0}^{\infty}&\Im \! \left(r^{i\,v}\,\eta \! \left(\frac{1}{2}+i\,v , w +\frac{1}{2}\right)\right)\frac{ \sinh \! \left(\phi \,v \right)}{\cosh \! \left(\pi \,v \right)}d v\\ \nonumber
=\int_{0}^{\infty}&\left({\cos \! \left(v\,\ln \! \left(r \right)\right) \eta_{I} \! \left(\frac{1}{2}+i\,v , w +\frac{1}{2}\right)}{}+{\sin \! \left(v\,\ln \! \left(r \right)\right) \eta_{R} \! \left(\frac{1}{2}+i\,v , w +\frac{1}{2}\right)}{}\right) \frac{ \sinh \! \left(\phi \,v \right)}{\cosh \! \left(\pi \,v \right)}\,d v \\ \nonumber
& = 
\sqrt{r} \left(r\,\cos \! \left(\frac{ \phi}{2}\right)\sin \! \left(\phi \right) \overset{\infty}{\underset{j =0}{\sum}}\, \frac{\left(-1\right)^{j}}{2\,r \left(j +w +\frac{1}{2}\right) \cos \! \left(\phi \right)+r^{2}+\left(j +w +\frac{1}{2}\right)^{2}}\right. \\ 
&\left. \hspace{20pt}-\sin \! \left(\frac{\phi}{2}\right) \overset{\infty}{\underset{j =0}{\sum}}\, \frac{\left(-1\right)^{j} \left(r\,\cos \! \left(\phi \right)+j +w +\frac{1}{2}\right)}{2\,r \left(j +w +\frac{1}{2}\right) \cos \! \left(\phi \right)+r^{2}+\left(j +w +\frac{1}{2}\right)^{2}}\right)\,.
\label{HI}
\end{align}

At this point, it is important to determine the analytic properties of these identities. To begin, both \eqref{HR} and \eqref{HI} are invariant under the interchange $\phi\leftrightarrow\,-\phi$. Are there any poles? Since it is known that the (alternating) series themselves are convergent, being simply the representation of a special case of the alternating function $\eta(s,a)$, we must ask if any of the individual terms of the sums are singular? Such a singularity would require that at least one term in the (shared) denominator of the sums vanish, in which case the index $j$ must satisfy
\begin{equation}
j=-r\,\cos \! \left(\phi \right)-w -\frac{1}{2}\pm r\,\sqrt{-1+\cos^{2}\left(\phi \right)}\,,
\label{jpm}
\end{equation}
which is impossible unless $|\cos(\phi)|=1$, that is $\phi=0$ or $\phi=\pm \pi$. In the case that $\phi=0$, a singularity will only occur if
\begin{equation}
j=-r-w-1/2
\label{Jzero}
\end{equation}
which cannot happen because here we demand both $0<\{r,w\}$. However, if $\phi\rightarrow\pm \pi$, then one term in the sum satisfying \eqref{Jzero} will diverge. Thus the case
\begin{equation}
r-w-1/2=m,\hspace{1cm} \phi=\pm\pi
\label{SpCaseDef}
\end{equation}
where $m$ is a non-negative integer must be treated separately when $\phi\rightarrow\pm \pi$, indicating that the right-hand sides sides of \eqref{HR} and \eqref{HI} are singular only if the condition \eqref{SpCaseDef} is satisfied. 

\begin{rem} 
\leavevmode
\begin{itemize}
\item{For reasonable values of $\phi\approx 0$, \eqref{HR} and \eqref{HI} were verified numerically.}
\item{Notice that the left-hand sides of \eqref{HR} and \eqref{HI} respectively contain factors $\frac{\cosh \, \left(\phi \,v \right)}{\cosh \, \left(\pi \,v \right)}=1$ and $ \frac{ \sinh \, \left(\phi \,v \right)}{\cosh \, \left(\pi \,v \right)}=\tanh(\pi v)\sim \pm 1$ when $\phi=\pm \pi$. Therefore, in this limiting case, the integrals may not converge, so their value and meaning will arise from analytic continuation in the variable $\phi$, the right-hand side being well-defined when $|\phi|<\pi$.} 
\end{itemize}
\end{rem}
\subsection{Special cases of \eqref{HR} and \eqref{HI}:  $\phi=0$} \label{sec:SpCasesHR}
\subsubsection{The Real part: Eq. \eqref{HR}:}
Let $\phi=0$ in \eqref{HR} to find, after a minor amount of simplification and identification of the right-hand side series,

\begin{equation}
\int_{0}^{{\infty}}\frac{\Re \! \left(r^{i\,v}\,\eta \! \left(1/2+i\,v , w +\frac{1}{2}\right)\right)}{\cosh \! \left(\pi \,v \right)}d v
 = \sqrt{r}\,\eta \! \left(1, 1/2+r +w \right)\,.
\label{CR2a}
\end{equation}
\begin{itemize}
\item{With $w=0$, \eqref{CR2a} reduces to }
\begin{equation}
\int_{0}^{\infty}\frac{\Re \! \left(\left(2\,r \right)^{i\,v}\,\beta \! \left(\frac{1}{2}+i\,v \right)\right)}{\cosh \! \left(\pi \,v \right)}d v
 = 
\sqrt{\frac{r}{2}}\,\,\eta \! \left(1, \frac{1}{2}+r \right)
\label{W0}
\end{equation}
with reference to \eqref{Beta(s)}, and further, if $r=1/2$ we find
\begin{equation}
\int_{0}^{\infty}\frac{\beta_{R}\! \left(\frac{1}{2}+i\,v \right)}{\cosh \! \left(\pi \,v \right)}d v
 = \frac{\ln \! \left(2\right)}{2}
\label{W0h}
\end{equation}
using \eqref{eta11}, a result that was numerically verified.

\item{Alternatively, in \eqref{CR2a} let $w=1/2$, to obtain with \eqref{etLim1}}

\begin{equation}
\int_{0}^{\infty}\frac{\Re \! \left(r^{i\,v}\,\eta \! \left(\frac{1}{2}+i\,v \right)\right)}{\cosh \! \left(\pi \,v \right)}d v
 = \sqrt{r}\,\eta \! \left(1, 1+r \right)\,,
\label{Wh1}
\end{equation}
which, with \eqref{etLim1}, can be written
\begin{equation}
\int_{0}^{\infty}\frac{\Re \! \left(r^{i\,v}\,\zeta \! \left(\frac{1}{2}+i\,v \right)\right)}{\cosh \! \left(\pi \,v \right)}d v -\sqrt{2} \int_{0}^{\infty}\frac{\Re \! \left(\left(\frac{r}{2}\right)^{i\,v}\,\zeta \! \left(\frac{1}{2}+i\,v \right)\right)}{\cosh \! \left(\pi \,v \right)}d v 
 = \sqrt{r}\,\eta \! \left(1, 1+r \right)\,.
\label{Wh}
\end{equation}.
Now, if $r=1$ and $r=2$, we respectively have
\begin{align} \nonumber
\int_{0}^{\infty}\left(\frac{\left(1-\sqrt{2}\,\cos \! \left(v\,\ln \! \left(2\right)\right)\right) \zeta_{R} \! \left(\frac{1}{2}+i\,v \right)}{\cosh \! \left(\pi \,v \right)}\right.&\left.-\frac{\sin \! \left(v\,\ln \! \left(2\right)\right) \sqrt{2}\,\zeta_{I} \! \left(\frac{1}{2}+i\,v \right)}{\cosh \! \left(\pi \,v \right)}\right)d v \\&
 =\,1 -\ln \! \left(2\right)
\label{dg1}
\end{align}
and
\begin{align} \nonumber
\int_{0}^{\infty}\left(\frac{\left(-\sqrt{2}+\cos \! \left(v\,\ln \! \left(2\right)\right)\right) \zeta_{R} \! \left(\frac{1}{2}+i\,v \right)}{\cosh \! \left(\pi \,v \right)}\right.&\left.-\frac{\sin \! \left(v\,\ln \! \left(2\right)\right) \zeta_{I} \! \left(\frac{1}{2}+i\,v \right)}{\cosh \! \left(\pi \,v \right)}\right)d v\\&
 = \sqrt{2}\,\ln \! \left(2\right)-\frac{\sqrt{2}}{2}\,,
\label{dg2}
\end{align}

which, by adding and subtracting (with a $\sqrt {2}\,$ multiplier) gives the identities

\begin{equation}
\int_{0}^{\mathit{\infty}}\frac{ \left(3-2\,\sqrt{2}\,\cos \! \left(v\,\ln \! \left(2\right)\right)\right)}{\cosh \! \left(\pi \,v \right)}\zeta_{R} \! \left({1}/{2}+i\,v \right)\,d v
 = 2\,-3\,\ln \! \left(2\right)
\label{dg3}
\end{equation}
and
\begin{equation}
\int_{0}^{\infty}\frac{3\,\sin \! \left(v\,\ln \! \left(2\right)\right) \zeta_{I} \! \left(\frac{1}{2}+i\,v \right)+\cos \! \left(v\,\ln \! \left(2\right)\right) \zeta_{R} \! \left({1}/{2}+i\,v \right)}{\cosh \! \left(\pi \,v \right)}d v
 = -\frac{\sqrt{2}}{2}
\label{dg1a}
\end{equation}
as well as several obvious variations (by changing the multipliers when adding and subtracting), such as

\begin{equation}
\displaystyle 2\, \sqrt{2}\int_{0}^{\infty }\!{\frac {\sin \left(v\ln  \left( 2 \right)  \right) \zeta_{I} \left( 1/2+iv \right) }{\cosh
\mbox{} \left( \pi\,v \right) }}\,{\rm d}v+\int_{0}^{\infty }\!{\frac {\zeta_{R} \left( 1/2+iv \right)   }{\cosh
\mbox{} \left( \pi\,v \right) }}\,{\rm d}v=\,-\ln  \left( 2 \right) \,
\label{R1b12pA}
\end{equation}
and
\begin{equation}
\displaystyle 6\,\int_{0}^{\infty }\!{\frac {\sin \left( v\ln  \left( 2 \right)  \right)  \zeta_{I} \left( 1/2+iv \right)   
\mbox{}}{\cosh \left( \pi\,v \right) }}\,{\rm d}v+2\,\int_{0}^{\infty }\!{\frac {\cos \left( v\ln  \left( 2 \right)  \right)  \zeta_{R} \left( 1/2+iv \right)   }{\cosh \left( \pi\,v \right) }}\,{\rm d}v
\mbox{}= \,- \sqrt{2}\,.
\label{R12b12pm}
\end{equation}

All of the above were numerically verified.
\end{itemize}

\subsubsection{The Imaginary part:  Eq. \eqref{HI}} \label{sec:SpCasesHI}



First, notice that if $\phi=0$, the identity \eqref{HI} reduces to the trivial statement $0=0$. However, because the integrations are convergent for $\phi\approx 0$, expanding \eqref{HI} about $\phi=0$ to first order, equating the coefficients of these terms (effectively differentiating with respect to $\phi$), and following the same procedures as in the previous subsection, gives the first moment equivalent of \eqref{CR2a}, that is
\begin{equation}
\int_{0}^{\infty}\frac{\Im \! \left(r^{-i\,v}\,\eta \! \left(\frac{1}{2}-i\,v , w +\frac{1}{2}\right)\right)}{\cosh \! \left(\pi \,v \right)}v\,d v
 = 
\sqrt{r} \left(-r\,\eta \! \left(2, \frac{1}{2}+r +w \right)+\frac{1}{2}\,\,{\eta \! \left(1, \frac{1}{2}+r +w \right)}\right)\,.
\label{Hi0}
\end{equation}
\begin{itemize}
\item{With $w=0$, we find the equivalent of \eqref{W0}, that is}

\begin{equation}
\int_{0}^{\infty}\frac{\Im \! \left(\left(2\,r \right)^{-i\,v}\,\beta \! \left(\frac{1}{2}-i\,v \right)\right) }{\cosh \! \left(\pi \,v \right)}\,v\,d v
 = 
\frac{\sqrt{2\,r} \left(-r\,\eta \! \left(2, \frac{1}{2}+r \right)+{\eta \left(1, \frac{1}{2}+r \right)}/2\right)}{2}\,.
\label{Hiw0}
\end{equation}
If $r=1/2$ this becomes

\begin{equation}
\int_{0}^{\infty}\frac{\beta_{I} \! \left(\frac{1}{2}+i\,v \right) }{\cosh \! \left(\pi \,v \right)}\,v\,d v
 = \frac{\pi^{2}}{48}-\frac{\ln \! \left(2\right)}{4}\,.
\label{Hiworh}
\end{equation}

\item{Alternatively, if $w=1/2$ we obtain}
\begin{align} \nonumber
\int_{0}^{\infty}\frac{\Im \! \left(r^{-iv}\,\zeta \! \left(\frac{1}{2}-iv \right) \right)}{\cosh \! \left(\pi \,v \right)}&v\,d v -\sqrt{2} \int_{0}^{\infty}\frac{\Im \! \left(\left(\frac{2}{r}\right)^{i\,v}\,\zeta \! \left(\frac{1}{2}-i\,v \right) \right)}{\cosh \! \left(\pi \,v \right)}vd v \\&
 = 
\sqrt{r} \left(\eta \! \left(1, 1+r \right)/2-r\,\eta \! \left(2, 1+r \right)\right).
\label{Hi0h}
\end{align}
Let $r=1$ then $r=2$ to respectively obtain the coupled pair of identities
\begin{align} \nonumber
\int_{0}^{\infty}&\frac{ \left(-\sqrt{2}\,\zeta_{R} \! \left(\frac{1}{2}+i\,v \right)\,\sin \! \left(v\,\ln \! \left(2\right)\right) +\zeta_{I} \! \left(\frac{1}{2}+i\,v \right) \left(-1+\sqrt{2}\,\cos \! \left(v\,\ln \! \left(2\right)\right)\right)\right)}{\cosh \! \left(\pi \,v \right)}v\,d v\\&
 \hspace{4.5cm}= \frac{\pi^{2}}{12}-\frac{1}{2}-\frac{\ln \! \left(2\right)}{2}
\label{Hi0ha}
\end{align}

and

\begin{align} \nonumber
\int_{0}^{\infty}&\frac{\left(-\zeta_{R} \! \left(\frac{1}{2}+i\,v \right)\,\sin \! \left(v\,\ln \! \left(2\right)\right)+\zeta_{I} \! \left(\frac{1}{2}+i\,v \right) \left(\sqrt{2}-\cos \! \left(v\,\ln \! \left(2\right)\right)\right)\right) }{\cosh \! \left(\pi \,v \right)}v\,d v\\&
 \hspace{4.5cm}= 
\frac{\sqrt{2}}{2} \left(\frac{5}{2}-\frac{\pi^{2}}{3}+\ln \! \left(2\right)\right).
\label{Hi0hb}
\end{align}
Adding and subtracting as before gives

\begin{equation}
\int_{0}^{\infty}\frac{ \left(3-2\,\sqrt{2}\,\cos \! \left(v\,\ln \! \left(2\right)\right)\right) }{\cosh \! \left(\pi \,v \right)}\zeta_{I} \! \left(1/2+i\,v \right)\,v\,d v
 = -\frac{5\,\pi^{2}}{12}+3+\frac{3\,\ln \! \left(2\right)}{2}
\label{Him}
\end{equation}
and
\begin{equation}
\int_{0}^{\infty}\frac{ 2\, \sqrt{2}\,\zeta_{R} \! \left(\frac{1}{2}+i\,v \right)\,\sin \! \left(v\,\ln \! \left(2\right)\right)-\zeta_{I} \! \left(\frac{1}{2}+i\,v \right)}{\cosh \! \left(\pi \,v \right)}v\,d v
 = \frac{\pi^{2}}{4}-2-\frac{\ln \! \left(2\right)}{2}\,.
\label{Hip}
\end{equation}


As before, with different multipliers, other variations are possible. For example


\begin{align} \nonumber
\displaystyle 
\int_{0}^{\infty }\,{\frac { \zeta_{I} \left( 1/2+iv \right)   \cos \left( v\ln  \left( 2 \right)  \right) 
}{\cosh \left( \pi\,v \right) }}\,v\,{\rm d}v&-3\,
\int_{0}^{\infty }\,{\frac { 
\zeta_{R} \left( 1/2+iv \right)  \sin \left( v\ln  \left( 2 \right)  \right)
\mbox{}}{\cosh \left( \pi\,v \right) }}\,v\,{\rm d}v\\&=
 \frac{\sqrt{2}}{4} \left( 3-\frac{\pi^2}{3} \right) \,.
\label{R2b12pmD}
\end{align}

\item{Since the integrals are convergent due to the exponentially increasing denominator term, higher moments are similarly available by multiple differentiation of \eqref{HI} with respect to $\phi$, yielding for example}
\begin{equation}
\displaystyle \int_{0}^{\infty }\,{\frac { \left( 3-2\, \sqrt{2}\cos \left( v\ln  \left( 2 \right)  \right)  \right) 
\mbox{} \zeta_{R} \left( 1/2+iv \right)   }{\cosh \left( \pi\,v \right) }}\,v^2\,{\rm d}v=\,-\frac{17}{2}+\frac{3\,\ln  \left( 2 \right) }{4}
\mbox{}-\frac{5\,{\pi}^{2}}{6}+{\frac {27\,\zeta \left( 3 \right) }{2}}\,.
\label{R1bmD}
\end{equation}
\item{Because $\sin(\pi\sigma)=\pm\cos(\pi\sigma)$ when $\sigma=1/4 \mbox{ or } 3/4$, possibilities for simpler results occur.}
Let $\sigma=1/4$ and follow the procedures outlined above by evaluating the difference between cases $r=1$ and $r=2$ with the multiplier $2^{3/4}$ to yield the following identity

\begin{align} \nonumber
\displaystyle\int_{0}^{\infty }\,&{ \left(1 - 2^{1/4}\cos \left( v\ln  \left( 2 \right)  \right)+ \sqrt{2}/4 \right)\frac {   \cosh \left( \pi\,v \right)\zeta_{R} \left( 1/4+iv \right)   +\sinh \left( \pi\,v \right) \zeta_{I} \left( 1/4+iv \right) 
\mbox{}  }{\cosh \left( 2\,\pi\,v \right) }}\,{\rm d}v\\
&=\,-\frac{\sqrt{2}}{2}  \left(\ln  \left( 2 \right)  -1/2 \right) +\frac{1}{4}
- \frac{\ln  \left( 2 \right)}{4}\,.
\label{Pr12a}
\end{align}
%
Similarly, for the case $\sigma=3/4$ we obtain
\begin{align} \nonumber
\int_{0}^{\infty} &\left(1-2^{\frac{3}{4}}\,\cos \! \left(v\,\ln \! \left(2\right)\right)+\frac{\sqrt{2}}{2}\right)\frac{\sinh \! \left(\pi \,v \right)\,\zeta_{I} \! \left(\frac{3}{4}+i\,v \right)-\cosh \! \left(\pi \,v \right)\,\zeta_{R} \! \left(\frac{3}{4}+i\,v \right)}{\cosh \! \left(2\,\pi \,v \right)}d v \\
& = 
\left(\frac{1}{2}+\frac{\sqrt{2}}{2}\right) \ln \! \left(2\right)-\frac{1}{2}-\frac{\sqrt{2}}{4}\,.
\label{Pr12b}
\end{align}
All the above were verified numerically.
\end{itemize}
\subsection{non-Exceptional Cases of \eqref{HR}  and \eqref{HI} where $\phi=\pm\pi$, general $w$} \label{sec:PhiPiallw}
%
%
%
As will be shown in the following subsections, to evaluate the limit $\phi\rightarrow\,\pm\,\pi$ is not straightforward. In preparation, we consider first some simple special cases. Rewriting  \eqref{HR} and \eqref{HI} and recalling the definition $a=r\,e^{i\phi}$, we have

\begin{equation}
\int_{0}^{\infty}\frac{\Re \! \left(r^{i\,v}\,\eta \! \left(\frac{1}{2}+i\,v , w +\frac{1}{2}\right)\right) \cosh \! \left(\phi \,v \right)}{\cosh \! \left(\pi \,v \right)}dv
 = 
\Re \! \left(r^{1/2}e^{i\phi/2}\overset{\infty}{\underset{j =0}{\sum}}\, \frac{\left(-1\right)^{j}}{j +\frac{1}{2}+r\,e^{i\phi} +w}\right)
\label{HR1}
\end{equation}
and
\begin{equation}
\int_{0}^{\infty}\frac{\Im \! \left(r^{i\,v}\,\eta \! \left(\frac{1}{2}+i\,v , w +\frac{1}{2}\right)\right) \sinh \! \left(\phi \,v \right)}{\cosh \! \left(\pi \,v \right)}d v
 = 
-\Im \! \left(r^{1/2} e^{i\phi /2} \overset{\infty}{\underset{j =0}{\sum}}\, \frac{\left(-1\right)^{j}}{j +\frac{1}{2}+r\,e^{i\phi} +w}\right)\,.
\label{HI1}
\end{equation}
\subsubsection{The Real part: \eqref{HR1}}
By applying the substitution $\phi=\pm \pi$ in \eqref{HR1}, {\bf if \boldmath $r-w-1/2\neq m$}, we find
\begin{equation}
\int_{0}^{\infty}\Re \! \left(r^{i\,v}\,\eta \! \left(1/2+i\,v , w +\frac{1}{2}\right)\right)d v = 0,\hspace{10pt}r-w-1/2\neq m
\,,
\label{IntREta}
\end{equation}
a result that might be difficult (impossible?) to obtain by studying the properties of the integral in isolation, or by numerical means. Yet, by the principle of analytic continuation in the variable $\phi$, the identity \eqref{IntREta} is true. 

Of particular interest are special cases of the above.
\begin{itemize}
\item{If $w=0,~r=m+1$, we find}
%
\begin{equation}
\int_{0}^{\infty}\left(\cos \! \left(v\,\ln \! \left(m +1\right)\right) \eta_{R} \! \left(\frac{1}{2}+i\,v , \frac{1}{2}\right)-\sin \! \left(v\,\ln \! \left(m +1\right)\right) \eta_{I} \! \left(\frac{1}{2}+i\,v ,\frac{1}{2}\right)\right)d v=0\,,
\label{Wc0}
\end{equation}
reducing to
\begin{equation}
\int_{0}^{\infty}\eta_{R} \! \left(1/2+i\,v , 1/2\right)d v = 0
\label{W0cm0}
\end{equation}
if $m=0$.  In terms of $\beta(s)$, \eqref{Wc0} becomes
\begin{equation}
\int_{0}^{\infty}\left(\beta_{I} \! \left({1}/{2}+i\,v \right)\,\sin \! \left(v\,\ln \! \left(2\,m +2\right)\right)-\beta_{R} \! \left(1/2+i\,v \right)\,\cos \! \left(v\,\ln \! \left(2\,m +2\right)\right)\right)d v = 0\,.
\label{Wc0b}
\end{equation}
\item{If $w=1/2,~r=M_{h}\equiv m+1/2$}, we obtain
\begin{equation}
\int_{0}^{\infty}\Re \! \left(M_{h}^{i\,v}\,\eta \! \left(1/2+i\,v \right)\right)d v=0\,\,,
\label{Whalfb}
\end{equation}
which, with reference to \eqref{etLim1} becomes

\begin{align} \nonumber
\int_{0}^{\infty}&\left(\left( \sqrt{2}\cos \! \left(v\,\ln \! \left({M_{h}}/{2}\right)\right)-\cos \! \left(v\,\ln \! \left(M_{h} \right)\right)\right) \zeta_{R} \! \left(1/2+i\,v \right)\right. \\ & \left. 
-\left(\sqrt{2}\,\sin \! \left(v\,\ln \! \left({M_{h}}/{2}\right)\right) -\sin \! \left(v\,\ln \! \left(M_{h} \right)\right)\right) \zeta_{I} \! \left(1/2+i\,v \right)\right)d v = 0\,.
\label{Whalfc}
\end{align}

\end{itemize}
\subsubsection{ The Imaginary part: \eqref{HI1}}
\begin{itemize}

\item{When $w=0$ and $r=M_{u}\equiv m+1$ \eqref{HI1} yields}
\begin{equation}
\int_{0}^{\infty}\Im \! \left(M_{u}^{i\,v}\,\eta \! \left(1/2+i\,v , 1/2\right)\right) \tanh \! \left(\pi \,v \right)d v
 = 
\left(\frac{1}{M_{h}}-\eta \! \left(1, M_{h}\right)+  \left(-1\right)^{m}\pi\right) \sqrt{M_{u}},
\label{CI1}
\end{equation}
reducing to
\begin{equation}
\int_{0}^{\infty}\eta_{I} \! \left(1/2+i\,v ,1/2\right)\tanh \! \left(\pi \,v \right)d v
 = 2+{\pi}/{2}\,
\label{CI1a}
\end{equation}

if $m=0$. 

\item{In the case that $w=1/2$ and $r=M_{h}\equiv m+1/2$, we find}

\begin{equation}
\int_{0}^{\infty}\Im \! \left(M_{h}^{i\,v}\,\eta \! \left(1/2+i\,v \right)\right) \tanh \! \left(\pi \,v \right)d v
 = 
\left(\eta \! \left(1,M_{h}\right)-\pi  \left(-1\right)^{m}\right) \sqrt{M_{h}}.
\label{CI2}
\end{equation}
If $m=0$, in terms of $\zeta(1/2+iv)$, \eqref{CI2} becomes
\begin{align} \nonumber
\int_{0}^{\infty}&\left(\left(\sqrt{2}\,\sin \! \left(2\,v\,\ln \! \left(2\right)\right)-\sin \! \left(v\,\ln \! \left(2\right)\right)\right) \zeta_{R} \! \left(1/2+i\,v \right) \right. \\ 
& \left.+\left(-\sqrt{2}\,\cos \! \left(2\,v\,\ln \! \left(2\right)\right)+\cos \! \left(v\,\ln \! \left(2\right)\right)\right) \zeta_{I} \! \left(1/2+i\,v \right)\right) \tanh \! \left(\pi \,v \right)d v
 = -\frac{\sqrt{2}}{4}\,\pi\,.
\label{CI2a}
\end{align}
\end{itemize}
\begin{rem}
It is interesting to consider the calculation of higher order moments by differentiating the basic equations with respect to $\phi$ prior to the evaluation of \eqref{HI1} with $\phi=\pm\pi$. Because all integrals are convergent when $|\phi|\leq\pi$ with the specified value of $w$, there are no issues with respect to the transposition of differentiation and integration operators, and any superficial limit operation reduces to simple substitution in both sides of \eqref{HI1}.  Differentiating once and twice respectively yields

\begin{equation}
\int_{0}^{\infty}\eta_{I} \! \left(1/2+i\,v\right) v d v
 = 0
\label{Diff1}
\end{equation}
and
\begin{equation}
\int_{0}^{\infty}\tanh \! \left(\pi \,v \right) \eta_{I} \! \left(1/2+i\,v \right)v^{2}\,d v
 = 
\frac{\ln \! \left(2\right)}{4}+\frac{\pi^{2}}{6}+\frac{3\,\zeta \! \left(3\right)}{2}\,.
\label{Diff2}
\end{equation}

\begin{rem}
Further differentiation here will reproduce the fourth order moment identity \eqref{HAB4} obtained in a different way in Section (\ref{sec:PvalExcpt}) below; for now, we observe by direct computation, that the third order moment vanishes, a property that is possibly true for all odd moments.
\end{rem}

\begin{rem}
Integrands appearing in the left-hand sides of the above are highly oscillatory. Thus convergence properties of the integrals themselves are difficult to deduce since they almost certainly diverge; however, divergent or not, the right-hand sides represent, and therefore assign value to, the left-hand sides by the principle of analytic continuation in the variable $\phi$ since both sides of the identities agree when $|\phi|<\pi$.  The correspondence is true, exact and does not involve limits.   At this point, readers unfamiliar with the subtleties of analytic continuation and regularization might wish to consult Definition \eqref{sec:AnalCont} and Appendices \ref{sec:Analogue} and \ref{sec:Reg} for simpler analogues.
\end{rem}
\section{Special Cases: Corollary \ref{sec:Cor2}} \label{sec:GenCor2}
\subsection{General Results}
%
Following the methods developed in the previous section, \eqref{M1ng2} gives us the general form
\begin{align} \nonumber
\frac{1}{2\pi\,\Gamma \! \left(b +1\right) }&\int_{-\infty}^{\infty}r^{-i\,v}\,{\mathrm e}^{v\,\phi}\,\Gamma \! \left(i\,v +1 -\sigma\right) \zeta\! \left(\sigma-i\,v +b  , w +1/2\right) \Gamma \! \left(\sigma-i\,v +b  \right)d v\\
& = 
\zeta \left(b +1, r\,{\mathrm e}^{i\,\phi}+w +1/2\right) r^{1-\sigma}\,{\mathrm e}^{i\,\phi  \left(1-\sigma \right)}\,,
\label{Ming2a}
\end{align}
with the replacement $\sigma:=1-\sigma$. Of particular interest is the case $\sigma>b$ and $b=0$, which requires that the residue of the pole at $\sigma=b$ in the original contour integral \eqref{M1ng2} be added to the left-hand side to eventually yield the general form
\begin{align} \nonumber
\frac{1}{2}&\int_{-\infty}^{\infty}\frac{r^{-i\,v}\,{\mathrm e}^{v\,\phi}\,\zeta\left(\sigma-i\,v  , w +\frac{1}{2}\right)}{\sin \left(\pi \,\sigma \right) \cosh \left(\pi \,v \right)-i\,\cos \left(\pi \,\sigma \right) \sinh \left(\pi \,v \right)}d v
\\ \nonumber 
&= \,
\underset{b \rightarrow 0}{\mathrm{lim}}\! \left(\left(\frac{1}{b}-\psi \! \left(r\,{\mathrm e}^{i\,\phi}+w +1/2\right)\right) r^{1-\sigma}\,{\mathrm e}^{i\,\phi  \left(1-\sigma \right)}-\frac{1}{b}{\mathrm e}^{-i \left(b +\sigma -1\right) \phi}\,r^{-b -\sigma +1}\right)\\
&=-\,r^{1-\sigma}e^{-i(\sigma-1)\phi}\left(\psi \! \left(r\,\cos \! \left(\phi \right)+w +1/2+i\,r\,\sin \! \left(\phi \right)\right)-\ln \! \left(r \right)-i\,\phi \right)  
\,.
\label{Ming2ca}
\end{align}
In \eqref{Ming2ca}, the limit \eqref{Zlim0} and the Gamma function reflection property \cite[Eq. (5.5.3)]{NIST} have been employed; notice that the left-hand side was simply evaluated with $b=0$ - any limits were applied only to the right-hand side because we have limited $0<\sigma<1$. Now consider the more restrictive case $\sigma=1/2$, and after some simple changes of integration variables we obtain
\begin{align} \nonumber
&\int_{0}^{\infty}\left(\frac{\left(\cos \! \left(v\,\ln \! \left(r \right)\right) \zeta_{R} \! \left(\frac{1}{2}+i\,v,w+1/2 \right)-\sin \! \left(v\,\ln \! \left(r \right)\right) \zeta_{I} \! \left(\frac{1}{2}+i\,v,w+1/2 \right)\right) \cosh \! \left(v\,\phi \right)}{\cosh \! \left(\pi \,v \right)}\right. \\ \nonumber
 &\left. -i\,\frac{ \left(\cos \! \left(v\,\ln \! \left(r \right)\right) \zeta_{I} \! \left(\frac{1}{2}+i\,v \right)+\sin \! \left(v\,\ln \! \left(r \right)\right) \zeta_{R} \! \left(\frac{1}{2}+i\,v,w+1/2 \right)\right) \sinh \! \left(v\,\phi\right)}{\cosh \! \left(\pi \,v \right)}\right)d v\\ \nonumber
 &= 
-\sqrt{r}\left({\mathrm e}^{i\phi/2} \psi \! \left(r\,\cos \! \left(\phi \right)+w +\frac{1}{2}+ir\,\sin \! \left(\phi \right)\right)-\left(i\,\phi+\ln \! \left(r \right)\right) \cos \! \left(\frac{\phi}{2}\right)\right.\\&\left.\hspace{5cm}-\sin \! \left(\frac{\phi}{2}\right) \left(i\,\ln \! \left(r \right)-\phi \right)\right) \,.
\label{Ming2e}
\end{align}
This result was easily evaluated numerically for small values of $\phi,\mbox{ with } w,r>0$. We now consider some special cases.

\subsection{The case $\phi=0$ using \eqref{Ming2e}} \label{sec:Phi0a}


\subsubsection{The Real part} \label{sec:Phi0Real}

\begin{itemize}
\item{With $w=0$ and \eqref{Sumh}, it is a simple matter to show that \eqref{Ming2e} becomes}
\begin{align} \nonumber
\int_{0}^{\infty}&\frac{\left(\sqrt{2}\,\cos \! \left(v\,\ln \! \left(2\,r \right)\right)-\cos \! \left(v\,\ln \! \left(r \right)\right)\right) \zeta_{R} \! \left(\frac{1}{2}+i\,v \right)}{\cosh \! \left(\pi \,v \right)}\\ \nonumber
&-\frac{\left(\sqrt{2}\,\sin \! \left(v\,\ln \! \left(2\,r \right)\right)-\sin \! \left(v\,\ln \! \left(r \right)\right)\right) \zeta_{I} \! \left(\frac{1}{2}+i\,v \right)}{\cosh \! \left(\pi \,v \right)}d v \\
& = 
-\sqrt{r}\left(\psi \! \left(r +1/2\right)-\ln \! \left(r \right)\right) \,,
\label{Ming2f}
\end{align}
which, with $r=1$ gives
\begin{align} \nonumber
\int_{0}^{\infty}\frac{\left(\sqrt{2}\,\cos \! \left(v\,\ln \! \left(2\right)\right)-1\right) \zeta_{R} \! \left(\frac{1}{2}+i\,v \right)}{\cosh \! \left(\pi \,v \right)}&-\frac{\sqrt{2}\,\sin \! \left(v\,\ln \! \left(2\right)\right) \zeta_{I} \! \left(\frac{1}{2}+i\,v \right)}{\cosh \! \left(\pi \,v \right)}d v
 \\
 &=\, -2+\gamma +2\,\ln \! \left(2\right)\,.
\label{Ming2fa}
\end{align}
Now, respectively adding and subtracting \eqref{Ming2fa} and \eqref{dg1} yields
\begin{equation}
2\,\sqrt{2} \int_{0}^{\infty}\frac{\sin \! \left(v\,\ln \! \left(2\right)\right) \zeta_{I} \! \left(\frac{1}{2}+i\,v \right)}{\cosh \! \left(\pi \,v \right)}d v
 = 1-\gamma -\ln \! \left(2\right)\,,
\label{X1}
\end{equation}
and
\begin{equation}
\int_{0}^{\infty}\frac{\left(\sqrt{2}\,\cos \! \left(v\,\ln \! \left(2\right)\right)-1\right) \zeta_{R} \! \left(\frac{1}{2}+i\,v \right)}{\cosh \! \left(\pi \,v \right)}d v
 =\, -\frac{3}{2}+\frac{\gamma}{2}+\frac{3\,\ln \! \left(2\right)}{2}\,.
\label{X2}
\end{equation}
Further, combining \eqref{X2} with \eqref{dg3} we find

\begin{equation}
\int_{0}^{\infty}\frac{\zeta_{R} \! \left(\frac{1}{2}+i\,v \right)}{\cosh \! \left(\pi \,v \right)}d v
 = \gamma-1
\label{Rz}
\end{equation}
and
\begin{equation}
\int_{0}^{\infty}\frac{\cos \! \left(v\,\ln \! \left(2\right)\right) \zeta_{R} \! \left(\frac{1}{2}+i\,v \right)}{\cosh \! \left(\pi \,v \right)}d v
 = \frac{-5+3\,\ln \! \left(2\right)+3\,\gamma}{2\,\sqrt{2}}\,.
\label{Rz2}
\end{equation}

\item{Now consider the case $w=1/2$, producing the identity}

\begin{equation}
\int_{0}^{\infty}\frac{\cos \! \left(v\,\ln \! \left(r \right)\right)\zeta_{R} \! \left(\frac{1}{2}+i\,v \right)-\sin \! \left(v\,\ln \! \left(r \right)\right) \zeta_{I} \! \left(\frac{1}{2}+i\,v \right)}{\cosh \! \left(\pi \,v \right)}d v
 = 
\sqrt{r}\left(\ln \! \left(r \right)-\psi \! \left(r +1\right)\right) \,.
\label{Rz3}
\end{equation}
If $r=1$ in \eqref{Rz3}, we immediately recover \eqref{Rz}. Let $r=2$ in \eqref{Rz3}, and combine with \eqref{dg1} to find
\begin{equation}
\int_{0}^{\infty}\frac{\zeta_{R} \! \left(\frac{1}{2}+i\,v \right) \left(2\,\sqrt{2}\,\cos \! \left(v\,\ln \! \left(2\right)\right)-1\right)}{\cosh \! \left(\pi \,v \right)}d v
 = -4+2\,\gamma +3\,\ln \! \left(2\right)
\label{Rz4}
\end{equation} 
and
\begin{equation}
\int_{0}^{\infty}\frac{2\,\sqrt{2}\,\sin \! \left(v\,\ln \! \left(2\right)\right) \zeta_{I} \! \left(\frac{1}{2}+i\,v \right)-\zeta_{R} \! \left(\frac{1}{2}+i\,v \right)}{\cosh \! \left(\pi \,v \right)}d v
 = 2-2\,\gamma -\ln \! \left(2\right)\,,
\label{Rz5}
\end{equation}
which, when combined with \eqref{Rz} reproduces \eqref{X1}. All of the above were validated numerically. 
\end{itemize}
\subsubsection{The Imaginary part}

As in Section (\ref{sec:SpCasesHI}), since the imaginary part of \eqref{Ming2e} vanishes identically when $\phi=0$, but the integral converges for other values of $\phi$, we transpose the integration and differentiation operator and differentiate once with respect to $\phi$, then set $\phi=0$, to arrive at the equivalent of \eqref{Hi0}

\begin{align} \nonumber
\int_{0}^{\infty}&\frac{\cos \! \left(v\,\ln \! \left(r \right)\right) \zeta_{I} \! \left( \frac{1}{2}+i\,v ,w+ \frac{1}{2}\right)+\sin \! \left(v\,\ln \! \left(r \right)\right) \zeta_{R} \! \left( \frac{1}{2}+i\,v ,w+ \frac{1}{2}\right) }{\cosh \! \left(\pi \,v \right)}\,v\,d v
  \\ = & 
\frac{\sqrt{r} \left(2\,r\,\psi^{\prime}\! \left(r +w +\frac{1}{2}\right)-\ln \! \left(r \right)+\psi \! \left(r +w +\frac{1}{2}\right)-2\right)}{2}\,.
\label{Ipart2}
\end{align}
\begin{itemize}
\item{If $w=0$, with \eqref{Sumh}, we arrive at }

\begin{align} \nonumber
\int_{0}^{\infty}&\left(\frac{\left(\sqrt{2}\,\sin \! \left(v\,\ln \! \left(2\,r \right)\right)-\sin \! \left(v\,\ln \! \left(r \right)\right)\right) \zeta_{R} \! \left(\frac{1}{2}+i\,v \right)}{\cosh \! \left(\pi \,v \right)} 
\right.\\&\left. \nonumber
+\frac{\left(\sqrt{2}\,\cos \! \left(v\,\ln \! \left(2\,r \right)\right)-\cos \! \left(v\,\ln \! \left(r \right)\right)\right) \zeta_{I} \! \left(\frac{1}{2}+i\,v \right)}{\cosh \! \left(\pi \,v \right)}\right) v d v \\
& = 
\frac{\left(2\,r\,\psi^{\prime}\! \left(r +\frac{1}{2}\right)-\ln \! \left(r \right)+\psi \! \left(r +\frac{1}{2}\right)-2\right) \sqrt{r}}{2}\,,
\label{Ip3}
\end{align} reducing to
\begin{align} \nonumber
\int_{0}^{\infty}\left(\frac{\sqrt{2}\,\sin \! \left(v\,\ln \! \left(2\right)\right) \zeta_{R} \! \left(\frac{1}{2}+i\,v \right)}{\cosh \! \left(\pi \,v \right)}\right.&\left.+\frac{\left(\sqrt{2}\,\cos \! \left(v\,\ln \! \left(2\right)\right)-1\right) \zeta_{I} \! \left(\frac{1}{2}+i\,v \right)}{\cosh \! \left(\pi \,v \right)}\right) v d v
 \\&= -4+\frac{\pi^{2}}{2}-\frac{\gamma}{2}-\ln \! \left(2\right)
\label{Ip3a}
\end{align}
if $r=1$.

\item{In the case that $w=1/2$, \eqref{Ipart2} becomes}

\begin{align} \nonumber
\int_{0}^{\infty}&\frac{\cos \! \left(v\,\ln \! \left(r \right)\right) \zeta_{I} \! \left(\frac{1}{2}+i\,v \right)+\sin \! \left(v\,\ln \! \left(r \right)\right) \zeta_{R} \! \left(\frac{1}{2}+i\,v \right) }  {\cosh \! \left(\pi \,v \right)}\,v\,d v \\
&\hspace{2cm} = 
\frac{\left(2\,r\,\psi^{\prime}\! \left(r +1\right)-\ln \! \left(r \right)+\psi \! \left(r +1\right)-2\right) \sqrt{r}}{2}
\label{Ip4}
\end{align}
and, if $r=1$, we find
\begin{equation}
\int_{0}^{\infty}\frac{ \zeta_{I} \! \left(\frac{1}{2}+i\,v \right)}{\cosh \! \left(\pi \,v \right)}\,v\,d v
 = -\frac{3}{2}+\frac{\pi^{2}}{6}-\frac{\gamma}{2}\,.
\label{Ip4a}
\end{equation}

Combining \eqref{Ip3a} and \eqref{Hi0hb} yields

\begin{align} \nonumber
\int_{0}^{\infty}&\left(-\frac{2\,\sqrt{2}\,\sin \! \left(v\ln \! \left(2\right)\right) \zeta_{R} \! \left(\frac{1}{2}+i\,v \right)}{\cosh \! \left(\pi \,v \right)}+\frac{\left(3-2\,\sqrt{2}\,\cos \! \left(v\ln \! \left(2\right)\right)\right) \zeta_{I} \! \left(\frac{1}{2}+i\,v \right)}{\cosh \! \left(\pi \,v \right)}\right) v d v \\
 & \hspace{3cm}=
\frac{13}{2}-\frac{5\,\pi^{2}}{6}+2\,\ln \! \left(2\right)+\frac{\gamma}{2};
\label{Ip5}
\end{align}
combining \eqref{Ip5} with \eqref{Him} gives the identity
\begin{equation}
\int_{0}^{\infty}\frac{\sin \! \left(v\,\ln \! \left(2\right)\right) \,\zeta_{R} \! \left(\frac{1}{2}+i\,v \right)}{\cosh \! \left(\pi \,v \right)}\,v\,d v
 = 
\frac{ \sqrt{2}}{4}\left(\frac{5\,\pi^{2}}{12}-\frac{7}{2}-\frac{\ln \left(2\right)}{2}-\frac{\gamma}{2}\right)\,.
\label{Ip5a}
\end{equation}
\end{itemize}

All of the above were verified numerically.

\subsection{The case $\phi=\pi$, the non-Exceptional cases}  \label{sec:phiPi}


In order to evaluate the limit $\phi\rightarrow \pi$ we shall need the real and imaginary parts of the right-hand side of \eqref{Ming2ca}. The real and imaginary parts of the relevant digamma (i.e. $\psi$) function, based on \eqref{PsiZ}, are listed in Appendix \ref{sec:Defs} (see \eqref{PsR} and \eqref{PsI}). As found previously, when $\phi=\pm\pi$, one term of the infinite series in \eqref{PsR} and/or \eqref{PsI} is singular if $r-w-1/2=m$. These are the exceptional cases. If $r-w-1/2 \neq m,~m\geq0$, by evaluating both sides of \eqref{Ming2ca} and setting $\phi=\pi$, we get the most general identities for the real and imaginary parts respectively, given in full in Appendix \ref{sec:LongEqs} (see \eqref{LhsRb} and \eqref{LhsIa}). 

\subsubsection{The real part}

\begin{itemize}
\item{In the case $w=0,~r\neq m+1/2$, the real part (i.e. \eqref{LhsRb}) becomes}
\begin{align} \nonumber
&\int_{0}^{\infty}\left(\frac{\left(C(v,\sigma,r)\,\cosh \! \left(\pi \,v \right) \sin \! \left(\pi \,\sigma \right)+S(v,\sigma,r)\,\sinh \! \left(\pi \,v \right) \cos \! \left(\pi \,\sigma \right)\right) \zeta_{R}(\sigma+iv)\,\cosh \! \left(\pi \,v \right)}{\cos \! \left(2\,\pi \,\sigma \right)-\cosh \! \left(2\,\pi \,v \right)} \right. \\& \nonumber \left.
+\frac{\left(S(v,\sigma,r)\,\cosh \! \left(\pi \,v \right) \sin \! \left(\pi \,\sigma \right)-C(v,\sigma,r)\,\sinh \! \left(\pi \,v \right) \cos \! \left(\pi \,\sigma \right)\right)  \zeta_{I}(\sigma+iv)\,\cosh \! \left(\pi \,v \right)}{-\cos \! \left(2\,\pi \,\sigma \right)+\cosh \! \left(2\,\pi \,v \right)}\right)d v \\
 &= 
\frac{\pi \,r^{1-\sigma}\,\sin \! \left(\pi \,\sigma \right)}{2}-\frac{r^{1-\sigma} \left(\psi \! \left(-r +\frac{1}{2}\right)-\ln \! \left(r \right)\right) \cos \! \left(\pi \,\sigma \right)}{2}
\label{Rparta}
\end{align}
where
\begin{equation}
C(v,\sigma,r)\equiv 2^{\sigma}\,\cos \! \left(v\,\ln \! \left(2\,r \right)\right)-\cos \! \left(v\,\ln \! \left(r \right)\right)
\label{Cdef}
\end{equation}
and
\begin{equation}
S(v,\sigma,r)\equiv 2^{\sigma}\,\sin \! \left(v\,\ln \! \left(2\,r \right)\right)-\sin \! \left(v\,\ln \! \left(r \right)\right)\,,
\label{Sdef}
\end{equation}
generating the special case
\begin{equation}
\int_{0}^{\infty}\left(\left(1-\sqrt{2}\,\cos \! \left(v\,\ln \! \left(2\right)\right)\right) \zeta_{R} \! \left(1/2+i\,v \right)+\sqrt{2}\,\sin \! \left(v\,\ln \! \left(2\right)\right)\zeta_{I} \! \left(1/2+i\,v \right)\right)d v
 = \pi,
\label{CR1}
\end{equation}
when $r=1$ and $\sigma=1/2$. Adding and subtracting \eqref{CR1} and \eqref{S0} gives 

\begin{equation}
\int_{0}^{\infty}\left(\sqrt{2}\,\cos \! \left(v\,\ln \! \left(2\right)\right)-1\right) \zeta_{R} \! \left(1/2+i\,v \right)d v
 = -\frac{\pi}{2}
\label{CR2}
\end{equation}
and
\begin{equation}
\int_{0}^{\infty}\sin \! \left(v\,\ln \! \left(2\right)\right) \zeta_{I} \! \left(1/2+i\,v \right)d v
 = \frac{\pi}{2\,\sqrt{2}} \,.
\label{CR3}
\end{equation}
Comparison of \eqref{CR3} with \eqref{s01pa} and \eqref{spboth} in the following section produces
\begin{equation}
\int_{0}^{\infty}\zeta_{R} \! \left(1/2+iv \right)d v =\, -\pi
\label{CR5}
\end{equation}
and
\begin{equation}
\int_{0}^{\infty}\cos \! \left(v\,\ln \! \left(2\right)\right) \zeta_{R} \! \left(1/2+i\,v \right)d v
 =\,-\frac{3\,\sqrt{2}}{4}\,\pi\,.
\label{CR4}
\end{equation}
\begin{rem} The general form of the identity corresponding to the real part when $w=1/2$ 
can be expressed symmetrically with respect to \eqref{Rparta}, and so is relegated to Appendix \ref{sec:LongEqs} as a matter of record and labelled \eqref{Rpartb}.
\end{rem} 

\item{When $w=1/2$,~ $r\neq m$ and $ \sigma=1/2$, from \eqref{Rpartb} we find the general result}
\begin{equation}
\int_{0}^{\infty}\Re \! \left(\zeta \! \left(1/2+i\,v \right) r^{i\,v}\right)\,d v
 = -\pi \,\sqrt{r}\,,
\label{Rx}
\end{equation}
which, if $r=1/2$, gives
\begin{equation}
\int_{0}^{\infty}\left(\sin \! \left(v\,\ln \! \left(2\right)\right)\zeta_{I} \! \left(1/2+i\,v \right)+\cos \! \left(v\,\ln \! \left(2\right)\right) \zeta_{R} \! \left(1/2+i\,v \right)\right)d v
 = -\frac{\pi}{\sqrt{2}}\,  \,,
\label{Rpc1}
\end{equation}
reproducing the sum of \eqref{CR3} and \eqref{CR4}.

\subsubsection{The imaginary part}

We now consider the imaginary part corresponding to \eqref{LhsIa} again with $\phi=\pi$, which, in the case:
\item{ $w=0$ and $r\neq m+1/2$, becomes}
\begin{align} \nonumber
&\int_{0}^{\infty}\left(\frac{\left(C \! \left(v , \sigma , r\right)\,\zeta_{R} \! \left(\sigma+i\,v \right)-S \! \left(v , \sigma , r\right)\,\zeta_{I} \! \left(\sigma+i\,v \right)\right) \cos \! \left(\pi \,\sigma \right)\sinh^{2}\left(\pi \,v \right)}{-\cos \! \left(2\,\pi \,\sigma \right)+\cosh \! \left(2\,\pi \,v \right)} \right.\\ \nonumber
& \left. +\frac{\left(-C \! \left(v , \sigma , r\right)\,\zeta_{I} \! \left(\sigma+i\,v \right)-S \! \left(v , \sigma , r\right)\,\zeta_{R} \! \left(\sigma+i\,v \right)\right) \sin \! \left(\pi \,\sigma \right) \cosh \! \left(\pi \,v \right) \sinh \! \left(\pi \,v \right)}{-\cos \! \left(2\,\pi \,\sigma \right)+\cosh \! \left(2\,\pi \,v \right)}\right)d v \\
& = 
r^{1-\sigma}\left(-\frac{\pi \,\cos \! \left(\pi \,\sigma \right)}{2}+\left(-\frac{\psi \! \left( \frac{1}{2}-r\right)}{2}+\frac{\ln \! \left(r \right)}{2}\right) \sin \! \left(\pi \,\sigma \right)\right) 
\label{Iparta}
\end{align} 
analogous to \eqref{Rparta}. If $r=1$ and $\sigma=1/2$, we obtain
\begin{align} \nonumber
\int_{0}^{\infty}&\left(\left(\sqrt{2}\,\cos \! \left(v\,\ln \! \left(2\right)\right)-1\right) \zeta_{I} \! \left(1/2+i\,v \right)+\sqrt{2}\,\sin \! \left(v\,\ln \! \left(2\right)\right) \zeta_{R} \! \left(1/2+i\,v \right)\right) \tanh \! \left(\pi \,v \right)d v \\
& = 2-\gamma -2\,\ln \! \left(2\right)\,.
\label{DR1}
\end{align}

Comparing \eqref{DR1} and \eqref{Wm0} in the following section, yields the identities

\begin{equation}
\int_{0}^{\infty}\zeta_{I} \! \left(1/2+i\,v \right) \left(\sqrt{2}\,\cos \! \left(v\,\ln \! \left(2\right)\right)-1\right)\tanh \! \left(\pi \,v \right) d v
 = \,-\frac{3\,\ln \! \left(2\right)}{2}+1-\frac{\gamma}{2}
\label{Wdm}
\end{equation}
and
\begin{equation}
\int_{0}^{\infty} \zeta_{R} \! \left(1/2+i\,v \right) \sin \! \left(v\,\ln \! \left(2\right)\right)\,\tanh \! \left(\pi \,v \right)d v
 =\, -\frac{\sqrt{2} \left(\ln \! \left(2\right)-2+\gamma \right)}{4}\,.
\label{Wdp}
\end{equation}
Subtracting \eqref{DR1} from \eqref{Wm1} gives
\begin{equation}
\int_{0}^{\infty} \zeta_{I} \! \left(1/2+i\,v \right)\tanh \! \left(\pi \,v \right)d v
 = -\gamma\,.
\label{Wdm1}
\end{equation}
Finally, adding\eqref{Wdm} and \eqref{Wdm1} yields 

\begin{equation}
\int_{0}^{\infty}\zeta_{I} \! \left(1/2+i\,v \right)\,\cos \! \left(v\,\ln \! \left(2\right)\right) \tanh \! \left(\pi \,v \right)d v
 =\, 
-\frac{3\,\sqrt{2}}{4}\,\gamma-\frac{3\,\sqrt{2}\,\ln \! \left(2\right)}{4}+\frac{\sqrt{2}}{2}\,.
\label{Wdm2}
\end{equation}


\item{In the case that both $w=1/2$ and $\sigma=1/2$, if $r\neq m+1$, \eqref{LhsIa} becomes}

\begin{align} \nonumber
\int_{0}^{\infty}{\left(\sin \! \left(v\,\ln \! \left(r \right)\right) \zeta_{R} \! \left({1}/{2}+i\,v \right)+\cos \! \left(v\,\ln \! \left(r \right)\right) \zeta_{I} \! \left({1}/{2}+i\,v \right)\right) \tanh \! \left(\pi \,v \right)}d v \\
 = 
\left(\pi \,\cot \! \left(\pi \,r \right)+\psi \! \left(r \right)-\ln \! \left(r \right)\right) \sqrt{r}\,.
\label{WhI1}
\end{align}

As expected, if $r=1/2$, adding \eqref{WhI1} and \eqref{Wdp} reproduces \eqref{Wdm2}.

\end{itemize}


\section{Exceptional Cases: Corollary \ref{sec:Cor1}, $\phi=\pm \pi$} \label{sec:Except_eta}
\subsection{The Pathology } \label{sec:PhiPiw0HR}


It will now be shown that analysis of the exceptional cases identified above is anything but simple. However, integrals of the type being discussed in this Section with $\phi=\pm\pi$, generate a template for many similar integrals of widespread interest in the literature. Therefore this case must be studied in great detail, by considering the identity \eqref{HR} more closely. With $\sigma=1/2$, we have

\begin{equation}
\int_{-\infty}^{\infty}a^{-i\,v}\,\eta \! \left(\frac{1}{2}-i\,v , w +\frac{1}{2}\right)d v
 = 
2\,\sqrt{a} \left(\overset{\infty}{\underset{\substack{j =0\\j\neq m}}{\sum}}\, \frac{\left(-1\right)^{j}}{a +w +\frac{1}{2}+j}+\frac{\left(-1\right)^{m}}{a +w +\frac{1}{2}+m}\right)
\label{C1a}
\end{equation}
indicating, on the right-hand side, the presence of a cut ($\sqrt{a}$) sitting along the negative real axis $a=\,-r$ and a single, simple, pole at the point corresponding to $a=-m-w-1/2$ lying atop that cut. {\bf Since it sits on the cut,} the pole is non-isolated and so the resulting analysis differs from convention that almost exclusively deals with isolated poles. We first investigate the nature of the combined  singularities. Consider the singular term $S$, defined  for some integer $j=m$, by 

\begin{equation}
S(r,\phi)\equiv S(a)=\frac{2\,\sqrt{a} \left(-1\right)^{m}}{a +w +\frac{1}{2}+m}\,.
\label{pole}
\end{equation}
In the immediate vicinity of the cut ($\epsilon\sim 0$) as it is approached from above or below

\begin{align} \nonumber
S \! \left(r ,\pm \pi \mp\epsilon\right)& = 
\frac{2\,\sqrt{r}\,\sqrt{-{\mathrm e}^{\mp i\,\epsilon}} \left(-1\right)^{m}}{-r\,{\mathrm e}^{\mp i\,\epsilon}+w +
\frac{1}{2}+m}\sim
\,\pm\frac{2\,i\,\sqrt{r} \left(-1\right)^{m}}{w-r +\frac{1}{2}+m}+\frac{\sqrt{r} \left(-1\right)^{m} \left(r +w +\frac{1}{2}+m \right) \epsilon}{\left(w-r +\frac{1}{2}+m \right)^{2}}\\
&\mp\frac{i\,\sqrt{r} \left(\frac{\left(2\,m +6\,r +2\,w +1\right)^{2}}{4}-8\,r^{2}\right) \epsilon^{2} \left(-1\right)^{m}}{4 \left(w-r  +\frac{1}{2}+m \right)^{3}}+\frac{O\! \left(\epsilon^{3}\right)}{\left(w-r  +\frac{1}{2}+m \right)^{4}}+\dots
\label{Ep2p}
\end{align}
indicating that {\bf if $r\neq w+1/2+m$}, the limiting value on the cut is given by the first term  of \eqref{Ep2p} where $\epsilon=0$. This term is imaginary, its sign changes with the variation of sign on the left-hand side of \eqref{C1a} (or equivalently \eqref{HI}) and so we may deduce that the real part of the non-exceptional cases corresponding to $\phi=\pm\pi$ will vanish, and the imaginary part will be finite and dependent on the direction of approach (from above or below). Special cases demonstrating this inference were given in previous sections.

\begin{figure}[h] 
\centering
\includegraphics[width=0.5\textwidth,height=.5\textwidth]{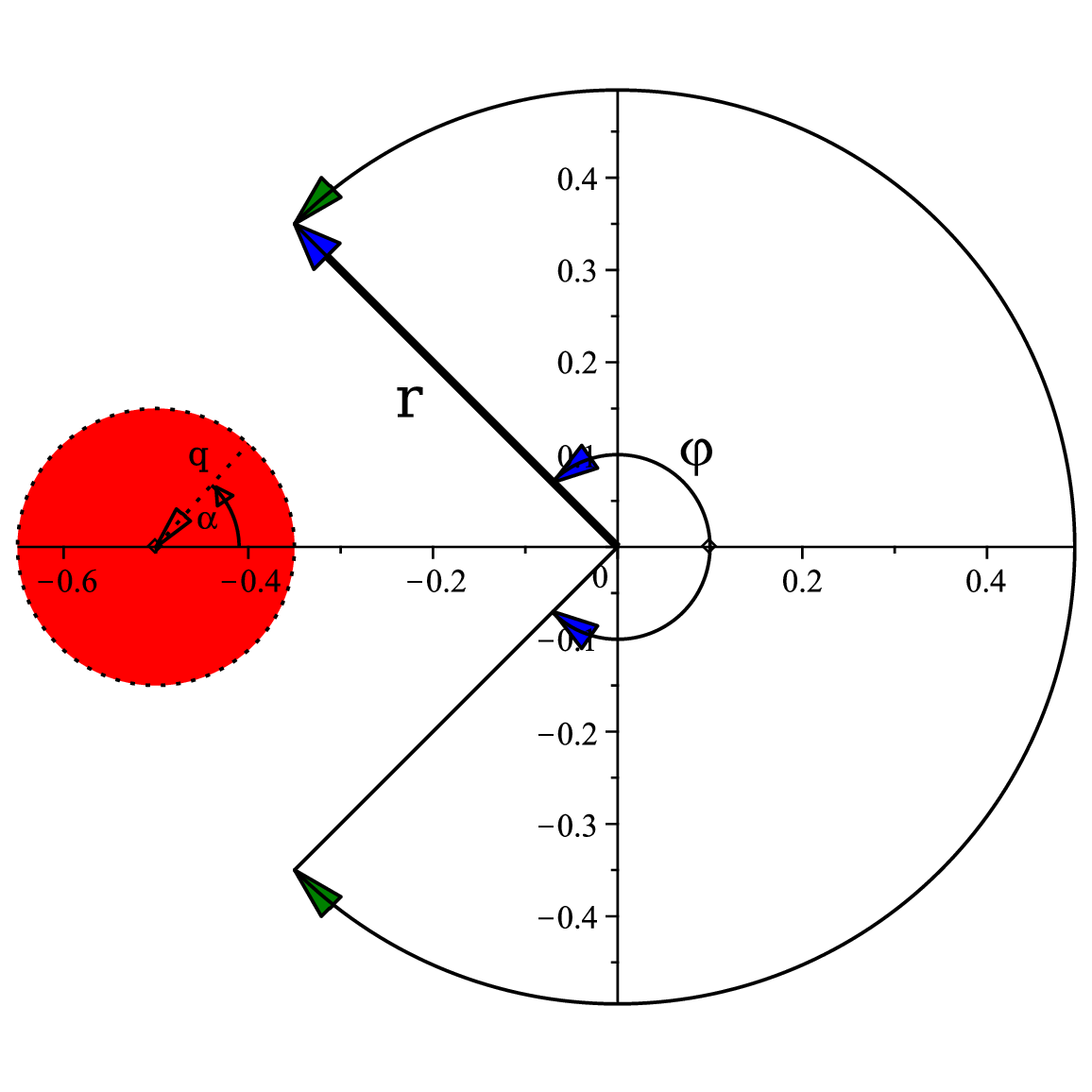}
\caption{With $m=0$, this shows two of the many different ways in which the exceptional limit point $(-1/2,i0)$ can be approached - either along the arc of an origin-centred circle of radius $r=\frac{1}{2}$ as $\phi\rightarrow\pm\pi$ - marked solid- or along a ray at an angle $\alpha$ inside a shrinking disk of radius $q$ as $q\rightarrow0$ marked red and dotted -- see \eqref{Sc1} and \eqref{Sc2} below. The limiting case corresponds to $\phi=\pm\pi,~r= 1/2,~\alpha=\{0,\pi\}$. In the case $r\rightarrow 1/2,~ \phi=\pm\pi$ and $\alpha=\{0,\pi\}$, both the singularity and the integral vanish as the exceptional point is approached along the real line $r\rightarrow 1/2^{\pm}$. Alternatively, the same point can be reached with $\phi\rightarrow\pm\pi$ when $r= 1/2$ or, equivalently $\alpha\rightarrow\pm\pi/2$ as $q\rightarrow 0$, in which case the integral diverges. For other values of $\alpha,$ the value of the integral varies depending on the path chosen -- see Subsection (\ref{sec:Discuss}). }
\label{fig:Diagram}
\end{figure}

Now, we ask what happens exactly {\bf at} the pole, since this corresponds to (symbolically) well-defined integrals of significant interest on the left-hand side of \eqref{HR} and \eqref{HI}. Continuing from \eqref{Ep2p}, let $r=m+w+1/2\pm\rho$, and examine the joint limits $\rho\sim 0$ and $\epsilon\sim 0$. This could be done by approaching the pole term by moving along the negative real $a=-r$ axis starting from a point far from the singularity ($\epsilon=0$) to land on the pole, either from the right ($\rho \sim 0^{-}$) or the left ($\rho\sim 0^{+}$) by considering the limit $r=m+w+1/2\pm\rho$ as $\rho\rightarrow 0$ and $\epsilon=0$. See Figure \ref{fig:Diagram}. We find
\begin{align} \nonumber
S \! \left(w+1/2+m\pm\rho ,\pi\right)&\sim\,-\frac{2\,i \left(-1\right)^{m}}{\sqrt{4\,m +4\,w +2}}\mp\frac{i\,\sqrt{4\,m +4\,w +2} \left(-1\right)^{m}}{\rho}
\\&\pm\frac{i \left(-1\right)^{m}\,\rho \,\sqrt{2}}{2 \left(4\,m +4\,w +2\right)^{\frac{3}{2}}}+\dots
\label{Ep4pX}
\end{align}

This result indicates that, by approaching the singularity along the real axis, the real part \eqref{HR1} will vanish at the singularity, whereas the contribution to the imaginary part \eqref{HI1} will be indefinite .  Alternatively, we may approach the limit point by choosing $r=m+w+1/2$ and consider the limit $\phi \rightarrow \pm\pi\mp\epsilon$, effectively approaching the singularity along a circle of radius $r=m+w+1/2$, (again see Figure \ref{fig:Diagram}) in which case we find
\begin{equation}
S \! \left(w+1/2+m ,\pm\pi\mp\epsilon\right)\sim\,\pm \frac{1}{\epsilon}\,\frac{4 \left(-1\right)^{m}}{\sqrt{4\,m +4\,w +2}}\pm\frac{\epsilon \left(-1\right)^{m}}{6\,\sqrt{4\,m +4\,w +2}}+\dots
\label{Ep5p}
\end{equation}
indicating that the real part is indefinite and the imaginary part vanishes. So, we have an ambiguity. Notice that the left-hand sides of \eqref{HR1} and \eqref{HI1} are immune to the ordering of limits applied to the right-hand side -- in fact, insofar as the left-hand sides are concerned, no ordered limits are required because $r>0$; it is a case of simple substitution: $r=m+w+1/2$ and $\phi=\,\pm\pi$ -- leaving an entity with a clear symbolic meaning, but a foggy value.

The important point here is that \eqref{pole} defines a singularity parametrized by two independent real variables, and as noted by Vladimirov who discusses a similar example, \cite[page 30]{Vlad}, Hartog's theorem: ``If a function is holomorphic with respect to each variable individually, it is also holomorphic with respect to the entire set of variables" is not valid for functions of real variables.
%

\subsubsection{Discussion} \label{sec:Discuss}
In short, we have arrived at a situation where an entity with a well-defined symbolic meaning has at least two, equally valid values. Which shall we choose? For the remainder of this section, I will explore various possibilities, in order to resolve the ambiguity inherent in \eqref{Ep4pX} and \eqref{Ep5p}. The goal is to obtain a preferred and consistent choice, but that is not to say it is the only one, given that the direction of approach to the singularity is arbitrary.

At this point it is worth asking: can the singularity associated with this ambiguity be identified as a classified singularity-- possibly an essential singularity because the ambiguous limits displayed by \eqref{Ep4pX} and \eqref{Ep5p} share some, but thus far, not all of the classification requirements of an essential singularity, one property of which requires that a function take on all values, except possibly one, as it approaches a limit point from all directions (Picard's Second theorem - see Definition \eqref{Picard} in Appendix \ref{sec:Defs}). As we have seen, the right-hand sides of \eqref{HR1} and \eqref{HI1} can be regularly expanded in a neighbourhood of the limit point, and the ratio (of independent variables) $\epsilon/\rho$ takes on infinitely many possible values as the limit point is approached along curves defined by those variables - see Figure \ref{fig:Diagram}. This suggests that Picard's Second theorem is satisfied.

For guidance, note that ambiguous limit points are known in the literature, but rarely discussed. In their tome dedicated to counterexamples, Gelbaum and Olmsted \cite[Chapter 9, Ex. 7]{CounterX} present just such an example. In a footnote, Hardy \cite[page 202] {Hardy} also gives an example but (perhaps wisely) defers further comment. Vladimirov makes reference to non-single-valued functions \cite[Section I.8.5]{Vlad} and gives an example \cite[Section 1.4.2]{Vlad} of a non-continuous function that is holomorphic for each variable individually in a neighbourhood of a limit point, but also eschews further discussion. Whittaker and Watson \cite[page 100 ff]{WW} discuss the properties of an essential singularity in terms of a Laurent expansion, but never focus on how to deal with an ambiguous singularity at the limit point. Uniquely among authorities, Morse and Feshbach \cite[page 380, ff]{M&F} distinguish between the value of an essential singularity {\it near} the singularity and {\it at} the singularity and comment that ``the behaviour of a function near an essential singularity is most complicated... their properties must therefore be carefully investigated in the neighbourhood of such points", but offer no further guidance. All of the above limit their comments to an isolated singularity.
\subsubsection{Analysis} \label{sec:Analysis}
The remainder of this Section is devoted to an examination of this phenomenon, as it applies to \eqref{HR1} annd \eqref{HI1}, particularly to determine what happens when the limit point is approached from all directions - since Picard's Second theorem appears to be satisfied, what is the nature of the Laurent series? 


To study the latter question, reconsider \eqref{Ep5p} expanded to more terms with less detail in order to survey the structure of the expansion. Applying the sequence of limits $\rho\rightarrow 0$ first, we find
\begin{align} \nonumber
\underset{\underset{\displaystyle \epsilon\rightarrow 0}{\displaystyle\rho\rightarrow 0}} {\lim} \,&S(r=m+w+1/2+\rho,\,\phi=\pi-\epsilon)\sim\left(\frac{73}{360}\,R_{1}\,\rho^{3}-\frac{1}{6}\,R_{2}\,\rho^{2}+\frac{1}{6}\,R_{3}\,\rho \right) i 
+\frac{2\,R_{4}}{\epsilon}\\& +\frac{\left(-\frac{16}{3}\,R_{1}\,\rho^{3}+4\,R_{2}\,\rho^{2}-4\,R_{3}\,\rho \right) i}{\epsilon^{2}}+\frac{16\,\rho^{3}\,R_{5}-8\,\rho^{2}\,R_{6}}{\epsilon^{3}}+\frac{16\,i\,R_{1}\,\rho^{3}}{\epsilon^{4}}+\dots
\label{P7Px}
\end{align}
where each of the coefficients $R_{j}$ is a function only of the parameters $m,w$ and contains no $\epsilon$ or $\rho$ dependence. The first terms in parentheses on the right-hand side indicates the existence of a regular part, the second term proves the existence of a single pole, and the other terms can be interpreted to be the first few terms of the Principal part of a Laurent series expansion about $\epsilon=0$. Notice that the higher order terms of the expansion vanish if $\rho=0$, leaving a single isolated simple pole at the limit point, but a Laurent expansion with a recognizable principal part nonetheless. Now, reconsidering \eqref{Ep4pX}, repeat the limiting exercise in the opposite order, that is apply the limit $\epsilon\rightarrow 0$ first, producing 
\begin{align} \nonumber
&\underset{\underset{\displaystyle \rho\rightarrow 0}{\displaystyle\epsilon\rightarrow 0}} {\lim} \,S(r=m+w+1/2+\rho,\,\phi=\pi-\epsilon)\sim\frac{\epsilon \,E_{15}}{8}-i\,E_{16}-\frac{i\,E_{1}}{\rho}\\ \nonumber
&+\frac{\left(-\frac{1}{24}\,\epsilon^{4}+\frac{1}{2}\,\epsilon^{2}\right) i\,E_{1}+\left(-\frac{1}{6}\,\epsilon^{3}+\epsilon \right) E_{6}}{\rho}
+\frac{\left(-\frac{7}{24}\,\epsilon^{3}+\frac{1}{4}\,\epsilon \right) E_{7}+\left(-\frac{5}{32}\,\epsilon^{4}\,E_{2}+\frac{3}{16}\,\epsilon^{2}\,E_{10}\right) i}{\rho^{2}}\\
&+\frac{\left(-\frac{25}{96}\,\epsilon^{4}\,E_{3}+\frac{1}{16}\,\epsilon^{2}\,E_{11}\right) i -\frac{\epsilon^{3}\,E_{8}}{4}}{\rho^{3}}
+\frac{-\frac{5}{32}\,\epsilon^{4}\,i\,E_{4}-\frac{1}{16}\,\epsilon^{3}\,E_{9}}{\rho^{4}}+\dots
\label{Ep2px}
\end{align}
where $E_{j}$ are expansion coefficients independent of both $\rho$ and $\epsilon$. Here we again see a regular term as well as a sequence of terms that can be interpreted to be the Principal part of a Laurent expansion about $\rho=0$. Significantly, the expansion factors in such a way that the terms of the Principal part series vanish if $\epsilon=0$, leaving a simple pole and a regular term at the limit point. This is all consistent with the tentative interpretation that the limit point under consideration ($\rho=0, \epsilon=0$) is an essential singularity with the unusual property that approached from certain directions, it reverts to a simple pole, a property antithetical to the properties of a run-of-the-mill essential singularity. It will now be shown how to utilize this property to overcome the ambiguity discussed above.

\subsubsection{Two possibilities when $w=0$}
One way to deal with the ambiguity, employing the ordering limit used in \eqref{Ep5p} with $r=m_{i}+1/2,\,i=1,2$, is to add the two results with multipliers $(-1)^{m_{1}}\sqrt{m_{1}+1/2}$ and $(-1)^{m_{2}}\sqrt{m_{2}+1/2}$,  and follow by setting $\epsilon=0$, to eventually arrive at a more general, unambiguous identity

\begin{equation}
\int_{0}^{\infty}\Re \! \left(\left(\left(-1\right)^{m_{1}}\,\left({m_{1}+1/2}\right)^{\frac{1}{2}+i\,v}-\left(-1\right)^{m_{2}}\,\left({m_{2}+1/2}\right)^{\frac{1}{2}+i\,v}\right) \beta \! \left(1/2+i\,v \right)\right)d v
 = 0\,.
\label{AnsBig}
\end{equation}

thereby removing an unremovable (potentially) essential singularity -- a textbook impossibility -- suggesting that it may not be an essential singularity. See Figures (\ref{fig:Y12X}) and (\ref{fig:Ans1L}).

\begin{figure}[h] 
\centering
\subfloat[The right-hand side of \eqref{HR1} as a function of $\phi$, combined in the same way as utilized to obtain \eqref{AnsBig} with $m_{1}=0\,\,m_{2}=1$, demonstrating that the value at the point corresponding to $\phi=\pi$ is regular, determinate and equal to zero as in \eqref{AnsBig}.]
{
\includegraphics[width=.45\textwidth,height=.45\textwidth]{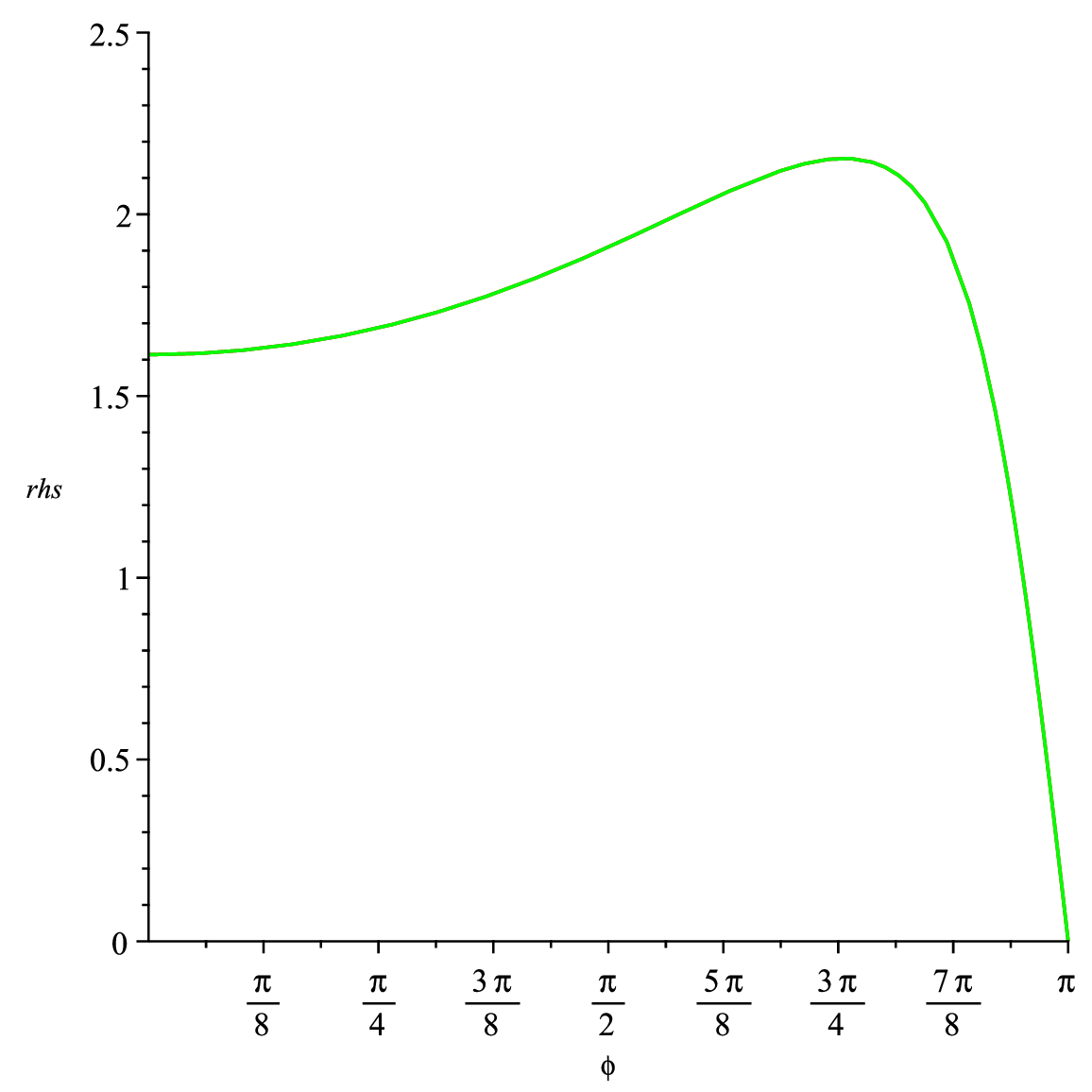}
\label{fig:Y12X}
}
\hfill
\subfloat[Values of the integrand defined by the left-hand side of \eqref{AnsBig} when $m_{1}=0\,,\,m_{2}=1$, as a function of $v$, over a randomly chosen range of $v$. It appears to be evenly balanced about the value zero, lending plausibility to the result.]
{
\includegraphics[width=0.45\textwidth,height=.45\textwidth]{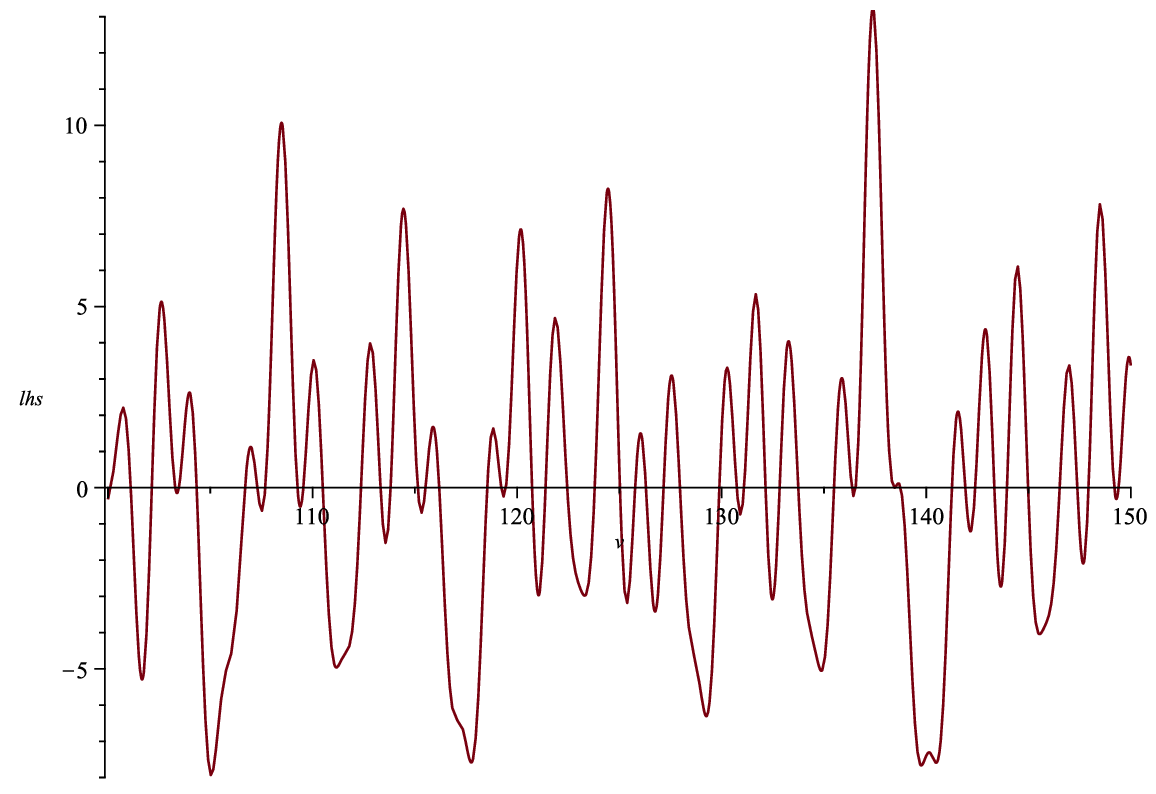}
\label{fig:Ans1L}
}
\caption{These figures investigate the consistency and plausibility of the identity developed in \eqref{AnsBig}.}
\end{figure}
Another, more general way of treating this ambiguous limit suggests that we follow Vladimirov \cite{Vlad} by rewriting \eqref{HR1} and \eqref{HI1} in terms of two new variables defined by
\begin{equation}
a\equiv M_{h}(-1+q\,e^{i\alpha})
\label{adef}
\end{equation}
where $M_{h}\equiv m+1/2$ and $q=0$ arrives at the exceptional limit point. Effectively, we are evaluating \eqref{HR1} and \eqref{HI1} along a radius (defined by the angle $\alpha$) inside a circle of ever diminishing radius $q\,M_{h}$ surrounding the exceptional complex point $a=(-M_{h},0i)$, corresponding to $\phi=\pm\pi$, $r=M_{h}$ in the previous co-ordinates. Again see Figure \ref{fig:Diagram}. After fairly lengthy calculations (Maple), we find the equivalents of \eqref{HR1} and \eqref{HI1} as follows


\begin{align} \nonumber
&\sin \! \left(\frac{\kappa \! \left(q , \alpha \right)}{2}\right) \int_{0}^{\infty}\Re \! \left((QM_{h}^2)^{-\frac{1}{4}+\frac{i\,v}{2}}\,\eta \! \left(\frac{1}{2}+i\,v , \frac{1}{2}\right)\right)\frac{ \cosh \! \left(v\,\kappa \! \left(q , \alpha \right)\right)}{\cosh \! \left(\pi \,v \right)}d v \\ \nonumber
&+\cos \! \left(\frac{\kappa \! \left(q , \alpha \right)}{2}\right) \int_{0}^{\infty}\Im \! \left((QM_{h}^2)^{-\frac{1}{4}+\frac{i\,v}{2}}\,\eta \! \left(\frac{1}{2}+i\,v , \frac{1}{2}\right)\right)\frac{ \sinh \! \left(v\,\kappa \! \left(q , \alpha \right)\right)}{\cosh \! \left(\pi \,v \right)}\,d v 
\\& = 
\frac{\left(-1\right)^{m}\,\sin \! \left(\alpha \right)}{q\,M_{h}}+\overset{\infty}{\underset{\substack{j =0\\j\neq m}}{\sum}}\! \left(\frac{\left(-1\right)^{j}\,M_{h}\,q\,\sin \! \left(\alpha \right)}{\left(j -m +M_{h}\,\cos \! \left(\alpha \right) q \right)^{2}+M_{h}^{2}\,q^{2} \sin^{2}\left(\alpha \right)}\right)\,,
\label{C1Icp}
\end{align}

and

\begin{align} \nonumber
&\cos \! \left(\frac{\kappa \! \left(q , \alpha \right)}{2}\right)\int_{0}^{\infty}\Re \! \left((QM_{h}^2)^{-\frac{1}{4}+\frac{i\,v}{2}}\,\eta \! \left(\frac{1}{2}+i\,v , \frac{1}{2}\right)\right)\frac{ \cosh \! \left(v\,\kappa \! \left(q , \alpha \right)\right)}{\cosh \! \left(\pi \,v \right)}d v \\ \nonumber
-&\sin \! \left(\frac{\kappa \! \left(q , \alpha \right)}{2}\right) \int_{0}^{\infty}\Im \! \left((QM_{h}^2)^{-\frac{1}{4}+\frac{i\,v}{2}}\,\eta \! \left(\frac{1}{2}+i\,v , \frac{1}{2}\right)\right)\frac{ \sinh \! \left(v\,\kappa \! \left(q , \alpha \right)\right)}{\cosh \! \left(\pi \,v \right)}d v\\
& = 
\frac{\left(-1\right)^{m}\,\cos \! \left(\alpha \right)}{q\,M_{h}}+\overset{\infty}{\underset{\substack{j =0\\j\neq m}}{\sum}}\! \frac{\left(-1\right)^{j} \left(j -m +M_{h}\,\cos \! \left(\alpha \right) q \right)}{\left(j -m +M_{h}\,\cos \! \left(\alpha \right) q \right)^{2}+M_{h}^{2}\,q^{2} \sin^{2}\left(\alpha \right)}\,,
\label{C1Rdp}
\end{align}
where 
\begin{equation}
Q\equiv 1-2\,q\,\cos \! \left(\alpha \right)+q^{2}
\label{Q}
\end{equation}
and 
\begin{equation}
\kappa \! \left(q , \alpha \right) \equiv \arg\left(-1+q\cos(\alpha)+i\,q\sin(\alpha)\right)\,.
\label{kappa}
\end{equation}

When $q=0$, $\kappa \! \left(0, \alpha \right)=\pm\pi$, creating a sign ambiguity in \eqref{C1Icp}, but leaving \eqref{C1Rdp} invariant. Both \eqref{C1Icp} and \eqref{C1Rdp} were verified numerically for $0<<q<1$ and random values of $\alpha$. However, in addition to the ambiguity attached to $\kappa(q,\alpha)$ discussed above, problems arise on the right-hand side when one attempts to determine the value of these identities {\bf at} the limit point $q=0$, which, of course, defines integrals of particular interest. Significantly, the numerator of the pole term on the right-hand side of both varies continuously and arbitrarily between $\pm\,1$ depending on $\alpha$, so the right-hand side function itself takes on all values in the neighbourhood of the point $q=0$ as required by Picard's Second theorem. In this respect, the limit point appears to satisfy the criterion that it is an essential singularity of the function on the right-hand side of both identities, and therefore creates a value ambiguity of the left-hand side by analytic continuation in the variable $q$. We now proceed to show how it can be unambiguously removed - a textbook impossibility for an essential singularity.

\subsection{Notation} \label{sec:Notate}
In the following, I use the symbol $ {\,\underset{q\rightarrow 0} \looparrowleft\,} $ to indicate that the {\bf value} of the left-hand side of an identity possesses an ambiguous {\bf directional limit} when the operation ${\underset{q\rightarrow 0} \lim}$ is applied to the right-hand side. This means that the limit point is approachable from {\bf any direction} in the complex plane, in contrast to a {\bf directed limit}, denoted by the commonly used superscripts ``$~{\bf^{\pm}}~$", utilized to indicate a limit where an ambiguity exists when all variables are real and the approach can only be made from either of two directions along the real axis. Thus, with reference to the left-hand side (LHS) of \eqref{C1Icp}, we could write
\begin{equation}
 \mbox{LHS} {\,\underset{q\rightarrow 0} \looparrowleft\,} \left\lbrace {\underset {\infty}  {^{O(q)}}} \right\rbrace
\end{equation}
to mean
\begin{equation}
{\,\underset{q\rightarrow 0} \lim} \;\mbox{LHS} =\,{\underset{q\rightarrow 0} \lim} \frac{\left(-1\right)^{m}\,\sin \! \left(\alpha \right)}{q\,M_{h}}+\cdots\,.
\end{equation}
For  any choice of $\alpha$, this means that the value is either finite (consider the iconoclastic choice $\alpha=\sin^{-1}(q)$) or indefinite, determined by the leading order of $q$. As shall be seen shortly in all the cases considered, a self-consistent family of integrals can be found by choosing $\alpha$ to yield finite values, in which case we write ``${\,\underset{\looparrowleft} =\,}$" simply to remind us that a choice has been made by means of a directional limit. In terms of the discontinuity discussed in Section \ref{sec:RealAn}, it simply means that the limit point is open-ended and the chosen value corresponds to the value arbitrarily close to, but never equal to the limit point. See \eqref{BetInTanX} below.

\subsubsection{Summary} \label{sec:Summary1}
To summarize the properties of such identities, consider \eqref{C1Rdp} when $q=0$. On the left-hand side, it is sufficient to simply substitute $q=0$ and allow the (symbolic) notation to convey its meaning, but, significantly, not its value. Insofar as the integrals are concerned, the variables $q$ and $\alpha$ are ``hidden variables" and the value of the integrals is determined by holomorphic continuation (regularization) in these variables. On the right-hand side, which defines the value of the left-hand side by analytic continuation in the variable $q$, there remains an ambiguity in the ordering of limits.

Specifically, to evaluate the limit $q\rightarrow 0$ {\bf at a constant value of $\alpha$}, one finds, from \eqref{C1Rdp}

\begin{equation}
\int_{0}^{\infty}\Re \! \left(\left(2\,M_{h} \right)^{-\frac{1}{2}+i\,v}\,\beta \! \left(1/2+i\,v \right)\right)d v
  {\,\underset{q\rightarrow 0} \looparrowleft\,}  \frac{\left(-1\right)^{m}\,\sin \! \left(\alpha \right)}{2\,q\,M_{h}}+\mathrm{O}(q),
\label{BetaReId}
\end{equation}

so that, if $m=0$ we have

\begin{equation}
\int_{0}^{\infty}\beta_{R} \! \left(1/2+i\,v \right)d v   {\,\underset{q\rightarrow 0} \looparrowleft}  \frac{\sin \! \left(\alpha \right)}{q} +\mathrm{O}(q).
\label{BetaReAv}
\end{equation}
To reiterate, if one chooses to approach the limit point $q=0$ along one of the two rays defined by $\sin(\alpha)=0$, (i.e. $\alpha=\{0,\pm\pi\}$ then the value of the right-hand side is governed by the leading (finite) term in the expansion in the variable $q$, otherwise it is indefinite.  From \eqref{C1Icp} we similarly have
\begin{equation}
\int_{0}^{\infty}\Im \! \left(\left(2\,M_{h} \right)^{-\frac{1}{2}+i\,v}\,\beta \! \left(1/2+i\,v \right)\right) \tanh \! \left(\pi \,v \right)d v
 {\,\underset{q\rightarrow 0} \looparrowleft} 
\frac{\eta \! \left(1, m +1\right)}{2}-\frac{\left(-1\right)^{m}\,\cos \! \left(\alpha \right)}{2\,q\,M_{h}},
\label{Ansx0}
\end{equation}
which, for the case $m=0$ yields
\begin{equation}
\int_{0}^{\infty}\beta_{I} \! \left(1/2+i\,v \right) \tanh \! \left(\pi \,v \right)d v
 {\,\underset{q\rightarrow 0} \looparrowleft\,} 
\frac{\ln \! \left(2\right)}{2}-\frac{\cos \! \left(\alpha \right)}{q}\,,
\label{BetInTan}
\end{equation}

and an ambiguity arises depending on whether one chooses to approach the limit point along the ray(s) defined by $\cos(\alpha)=0$ or not. That is, if we choose $\alpha=\pm \pi/2$, i.e. the imaginary axis, the integral \eqref{BetInTan} is finite as is \eqref{Ansx0}. Thus \eqref{BetInTan} with the choice $\alpha=\pi/2$ would be written

\begin{equation}
\int_{0}^{\infty}\beta_{I} \! \left(1/2+i\,v \right) \tanh \! \left(\pi \,v \right)d v
{\,\underset{\looparrowleft} =\,}
\frac{\ln \! \left(2\right)}{2}\,.
\label{BetInTanX}
\end{equation}

Finally, and significantly, in order to obtain a finite result in \eqref{C1Rdp} with $\alpha$ constant (a constraint added here for simplicity), we require the choice $\sin(\alpha)                                                                 =0$, i.e. $\alpha=\{0,\pm\pi\}$, whereas to obtain a finite result in \eqref{C1Icp} with $\alpha$ constant, we require the choice $\cos(\alpha)=0$, i.e. $\alpha={\pm\pi/2}$, suggesting an inconsistency between the choice of limits for the real and imaginary parts of the two identities, unless the ambiguity is banished. It will shortly be shown how this can be achieved, but first, for improved insight, we consider an example from the literature (see \cite{MilHughMMM2025}) that involves only real variables and directed limits.

\subsubsection{An Analytic Analogue} \label{sec:RealAn}


In the following, we reiterate that $r>0$ is real and $m$ is a non-negative integer. ``$\lfloor...\rfloor$" is the \textbf{open-ended} floor function (see below).
\begin{theorem} \label{sec:BetaThm}
\begin{equation}
\Re\int_{0}^{\infty}\frac{\eta \left(\frac{1}{2}+i\,v ,\frac{1}{2}\right) r^{\frac{1}{2}+i\,v}}{\frac{1}{2}+i\,v}d v
 = 
\left\{\begin{array}{cc}
\left(-\frac{\left(-1\right)^{{\lfloor r +\frac{1}{2}\rfloor}}}{2}+\frac{1}{2}\right) \pi  & r \neq m +\frac{1}{2} 
\\
 \frac{\pi}{2} & r =m +\frac{1}{2} 
\end{array}\right.
\label{BetaInt}
\end{equation}
\begin{proof}
From its fundamental definition \eqref{HzetaAlt}, write
\begin{equation}
\int_{0}^{\infty}\frac{\eta \left(\frac{1}{2}+i\,v ,\frac{1}{2}\right) r^{\frac{1}{2}+i\,v}}{\frac{1}{2}+i\,v}d v
 = 
\sqrt{r}\;\overset{\infty}{\underset{k =0}{\sum}}\frac{ \left(-1\right)^{k} }{\sqrt{k +\frac{1}{2}}}\;\int_{0}^{\infty}\frac{r^{i\,v}}{\left(k +\frac{1}{2}\right)^{i\,v} \left(\frac{1}{2}+i\,v \right)}d v 
\label{L2}
\end{equation}

after transposing the summation and integration because both are convergent. We now consider the real part of the integrand
\begin{equation}
\Re \left(\frac{r^{i\,v}}{\left(k +\frac{1}{2}\right)^{i\,v} \left(\frac{1}{2}+i\,v \right)}\right)
 = 
\frac{-v\,\sin \left(v\,\ln \left(\frac{2\,k +1}{2\,r}\right)\right)+\frac{1}{2}\cos \left(v\,\ln \left(\frac{2\,k +1}{2\,r}\right)\right)}{v^{2}+\frac{1}{4}}
\label{GintR}
\end{equation}
and, with some help\footnote {For a more thorough proof, see \cite[Section (2.2)]{MilHughMMM2025}.} from Maple \cite{Maple24}, if $k<r+1/2$, then

\begin{equation}
\int_{0}^{\infty}\frac{-v\,\sin \left(v\,\ln \left(\frac{2\,k +1}{2\,r}\right)\right)+\frac{1}{2}\cos \left(v\,\ln \left(\frac{2\,k +1}{2\,r}\right)\right)}{v^{2}+\frac{1}{4}}
 = \pi\frac{\sqrt{k + \frac{1}{2}}}{\sqrt{r}}\,.
\label{GintRL}
\end{equation}
Further, if $r=m+1/2$ and $k=m$, then
\begin{equation}
\int_{0}^{\infty}\frac{-v\,\sin \left(v\,\ln \left(\frac{2\,k +1}{2\,r}\right)\right)+\frac{1}{2}\cos \left(v\,\ln \left(\frac{2\,k +1}{2\,r}\right)\right)}{v^{2}+\frac{1}{4}}
 = \frac{\pi}{2}\,,
\label{GintRE}
\end{equation} 
and, if $k>r+1/2$, then
\begin{equation}
\int_{0}^{\infty}\frac{-v\,\sin \left(v\,\ln \left(\frac{2\,k +1}{2\,r}\right)\right)+\frac{1}{2}\cos \left(v\,\ln \left(\frac{2\,k +1}{2\,r}\right)\right)}{v^{2}+\frac{1}{4}}
 =0\,.
 \label{GintRM}
 \end{equation}
Summing, as in \eqref{L2} where we use the floor function to define the various summation limits, immediately yields \eqref{BetaInt}. See Figure \ref{fig:Bsum}.

\end{proof}
\end{theorem}

\begin{figure}[h] 
\centering
\includegraphics[width=.9\textwidth,height=.5\textwidth]{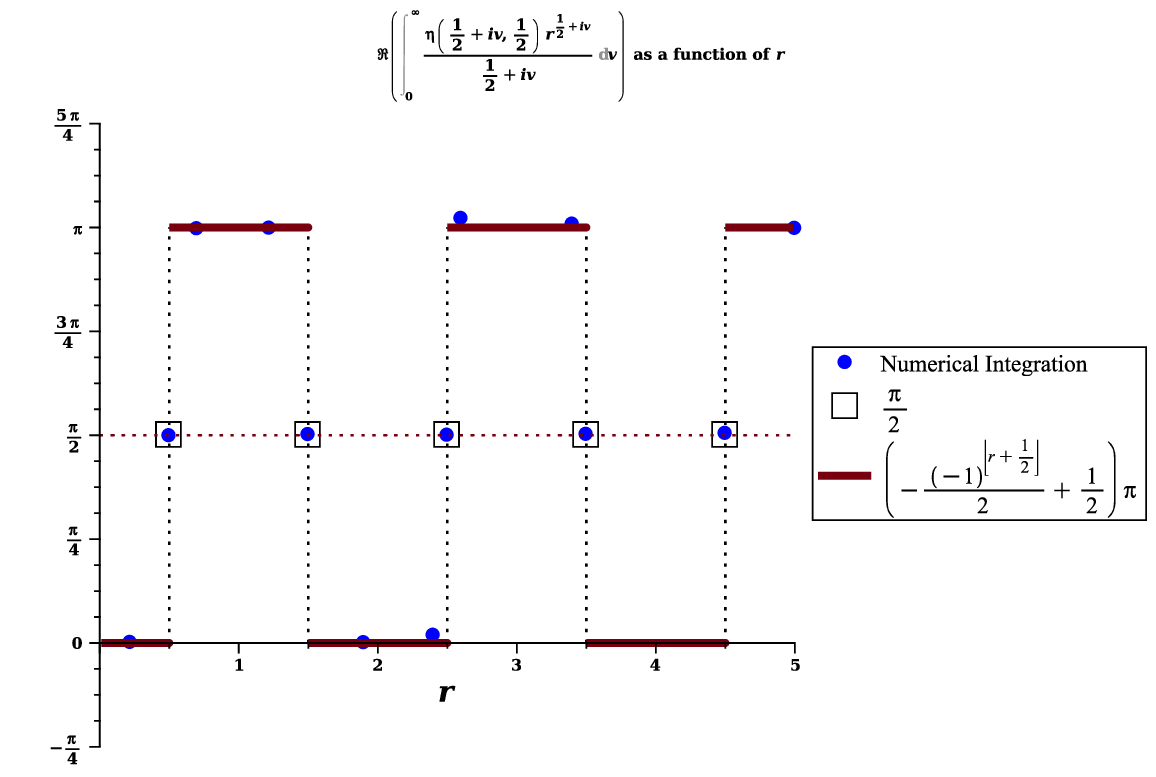}
\caption{This Figure shows the right-hand side of Theorem \ref{sec:BetaThm} calculated as a function of $r$. It is sometimes known as the boxcar function. The left-hand side of Theorem \ref{sec:BetaThm} can be shown to agree using numerical integration \cite{Math23} and selected values of $r$. As usual, numerical integration becomes less stable as $r$ increases, but as seen in \cite{MilHughMMM2025}, it is surprisngly stable at the half-integers.}
\label{fig:Bsum}
\end{figure}

In Figure \ref{fig:Bsum}, values at the discontinuities are calculated and happen to coincide, both analytically and numerically, with the mean of the values as $r\rightarrow (m+1/2)^{\pm}$, a property that is usually imposed by fiat because values that occur precisely at a discontinuity are not always known. For this reason, we need to relinquish the classic definition of limit and require that the floor function $\lfloor r \rfloor $ means the greatest integer less than, {\bf but not equal} to $r$, as indicated in \eqref{BetaInt}. If we now differentiate \eqref{BetaInt} with respect to $r$, and although the integral may not converge, it is clear that the integral in its regularized form becomes
\begin{equation}
\Re\left(\frac{1}{r}\int_{0}^{\infty}\eta \left(1/2+i\,v,1/2 \right) r^{\frac{1}{2}+i\,v}d v\right)=0\hspace{30pt} r\neq m+1/2\,,
\label{EtaInt}
\end{equation}
because $\frac{d}{d r}{\lfloor r +\frac{1}{2}\rfloor} = 0$ if $r\neq m+1/2$.

If $r=m+1/2$, the concept of differentiating fails at a point of discontinuity and the right-hand side of \eqref{EtaInt} becomes indeterminate, being either zero, if interpreted as  the limit of the derivative of the right-hand side of \eqref{BetaInt} as $r\rightarrow (m+1/2)^{\pm}$, or infinite using the standard definition of derivative (see \cite{MilHughMMM2025} for more details). Therefore we can finally write
\begin{Coroly}
\begin{equation}
\Re\left(\int_{0}^{\infty}\eta \left(1/2+i\,v,1/2 \right) r^{\frac{1}{2}+i\,v}d v\right)=(-1)^{r-1/2}\,\pi r\delta(r-(m+1/2)) ,
\label{EtaIntd}
\end{equation}
where $\delta(x)$ is the Dirac delta function. 
\end{Coroly}
\begin{proof}
This  follows exactly the proof of \cite[Theorem 3.1]{MilHughMMM2025} with the replacements $\zeta(\sigma+it):=(-1)^{r-1/2}\eta(1/2+it,1/2)$ and $n:=m+1/2$, as well as use of the elementary identity
$$\Re\left( \int_{0}^{\infty}\eta \left(\sigma+i\,v ,1/2 \right) r^{\sigma+i\,v }d v \right)
 = 
\frac{1}{2}\int_{-\infty}^{\infty}\eta \left(\sigma+i\,v  , 1/2\right) r^{\sigma+i\,v }d v .$$
\end{proof}
\begin{rem}
 Because the function $\eta \left(\sigma+i\,v  , \frac{1}{2}\right)$ is entire, there is no singularity comparable to that which occurs in the function $\zeta(1+iv)$, thereby justifying the replacement $\sigma:=1/2$ in the proof, because in this case, the contour of integration can be shifted from $\sigma>1$ used in \cite[Theorem 3.1]{MilHughMMM2025}, to $\sigma=1/2$, without corrections.\newline
\end{rem}

\begin{rem} 
With this analogue, it has been shown how a limiting case of the integral under discussion is conventionally analyzed when entities are limited to the realm of real variables. It contrasts with how that very same discontinuity presents when it is encountered in the complex realm. Recalling \eqref{Beta(s)}, we can compare \eqref{EtaIntd} with \eqref{BetaReId}. That is
\begin{equation}
\Re\left(\int_{0}^{\infty}\eta \left(1/2+i\,v , 1/2\right)\,\left(m +1/2\right)^{\frac{1}{2}+i\,v}\,d v\right)
 {\,\underset{q\rightarrow 0} \looparrowleft\,}  \frac{\left(-1\right)^{m}\,\sin \left(\alpha \right)}{q}\,,
\label{EtaBeta}
\end{equation}
showing how a complex singularity of the form discussed here, projects in the form of a Dirac delta function when limited to the real line. In both cases, if $r\neq m+1/2$, the integral vanishes because $\sin(\alpha)=0$, consistent with \eqref{IntREta}.

In summary, this comparison suggests that, when evaluating the pathology studied here, if finite values of an ambiguous limit are chosen, then the symbol ``${\,\underset{\looparrowleft} =\,}$" labels an equivalent variant of the classical directed $limit$ operation, except that the limit point is arbitrarily close, but never reached.
\end{rem}

\end{rem}


\subsection{Removal of the ambiguity} \label{sec:Remove}

\subsubsection{Method A: Weighted Difference} \label{sec:MremA}

Here, we investigate one possible way to remove the ambiguity embedded in the above in order to obtain consistency and limit the analysis to that of finite quantities, by applying the temporary condition that $\alpha\neq \pm  n\pi/2,\, n\geq0$. We then consider a weighted difference of \eqref{C1Icp} and \eqref{C1Rdp}, the former multiplied by $\cos(\alpha)$, the latter by $\sin(\alpha)$, yielding a single combined identity
\begin{align} \nonumber
\int_{0}^{\infty}
 & -\Im \! \left(\eta \! \left(1/2+i\,v ,1/2\right) \left(Q\,M_{h}^{2}\right)^{{i\,v}/{2}-1/4}\right) \cos \! \left(\alpha -\frac{\kappa(q,\alpha)}{2}\right) \frac{\sinh \! \left(v\,\kappa(q,\alpha)\right)}{\cosh \! \left(\pi \,v \right)}  
\\ \nonumber
& +\Re \! \left(\eta \! \left(1/2+i\,v , 1/2\right) \left(Q\,M_{h}^{2}\right)^{{i\,v}/{2}-1/4}\right) \sin \! \left(\alpha -\frac{\kappa(q,\alpha)}{2}\right)\frac{ \cosh \! \left(v\,\kappa(q,\alpha)\right)}{\cosh \! \left(\pi \,v \right)} d v
  \\
&= \sin \! \left(\alpha \right) \overset{\infty}{\underset{\substack{{j =0}\\j\neq m}}{\sum}}\,\,\, \frac{\left(-1\right)^{j} \left(j -m \right)}{2\,M_{h}\,q \left(j -m \right) \cos \! \left(\alpha \right)+M_{h}^{2}\,q^{2}+\left(j -m \right)^{2}}\,.
\label{Ans12Lb2}
\end{align}


Since the divergent term no longer appears on its right-hand side, \eqref{Ans12Lb2} is a simple, straightforward identity that can be (and was) tested numerically for reasonable values of $0\ll q<1$ and any value(s) of $\alpha$ or $M_{h}$ one desires. However, as $q\rightarrow 0$, the left-hand side becomes less and less amenable to numerical evaluation, and eventually confounds both Mathematica and Maple as $q\approx 0$, since it becomes heavily oscillatory and numerically intractable at the limit point $q=0$. But, since the right-hand side of \eqref{Ans12Lb2} is an entire function of $\alpha$ for $q\geq 0$, and, by simply setting $q=0$, the right-hand side yields the identity 
\begin{equation}
\overset{\infty}{\underset{\substack{{j =0}\\j\neq m}}{\sum}}\,\, \frac{\left(-1\right)^{j} \left(j -m \right)}{2\,M_{h}\,q \left(j -m \right) \cos \! \left(\alpha \right)+M_{h}^{2}\,q^{2}+\left(j -m \right)^{2}}\left|_{q=0}\right.=\,-\eta(1,m+1)\,.
\label{check}
\end{equation}

So, the singularity, but not the ambiguity,  has been removed for the entire range of both variables $q$ and $\alpha$, and with no remaining singularities, the range $q\geq0$ and $|\alpha|\le\pi$ is a holomorphic domain of \eqref{Ans12Lb2}. Thus, by the laws of analytic and holomorphic continuation, \eqref{Ans12Lb2} is a perfectly well-defined identity for all $\alpha$ including the points $\alpha=n\pi/2$, the validity of which was numerically verified when $q\gg0$.

Therefore, when $\sin(\alpha)=0$, the left-hand side of \eqref{Ans12Lb2} remains represented by its right-hand side, which vanishes for all values of $q\geq 0$. By simple substitution in \eqref{Ans12Lb2}, special cases of interest with $q=0$ can be obtained by choosing special values of $\alpha$; first set $\alpha=\{0,\pm \pi\}$, to obtain

\begin{equation}
\int_{0}^{\infty}\Re \! \left(\eta \! \left(1/2+i\,v , 1/2\right) M_{h}^{i\,v}\right)d v
 = 0\,,
\label{Sc1}
\end{equation}
which, with $m=0$ yields
\begin{equation}
\int_{0}^{\infty} \beta_{R} \! \left(1/2+i\,v \right)d v = 0\,,
\label{BRId}
\end{equation}
and second, let $\alpha=\pm \pi/2$, to find

\begin{equation}
\int_{0}^{\infty}\Im \! \left(M_{h}^{i\,v}\,\eta \! \left(1/2+i\,v , 1/2\right)\right) \tanh \! \left(\pi \,v \right)d v
 = \sqrt{M_{h}}\,\,\eta \, \left(1, m +1\right)\,,
\label{Sc2}
\end{equation}
which, with $m=0$ yields
\begin{equation}
\int_{0}^{\infty} \beta_{I} \! \left(1/2+i\,v \right) \tanh \! \left(\pi \,v \right)d v
 = \frac{\ln \! \left(2\right)}{2}\,,
\label{BImId}
\end{equation}
both of which are consistent with \eqref{BetaReAv} and \eqref{BetInTan}, if the finite variation of both had been chosen. 

\begin{rem} The identity \eqref{Sc1} is also consistent with the analogue case  \eqref{EtaInt}. If we  set $r=M_{h}$ in \eqref{EtaInt}, then \eqref{Sc1} and \eqref{EtaInt} agree, because the choice $\alpha=\{0,\pm\pi\}$ in Figure \ref{fig:Diagram}, corresponds exactly to approaching the discontinuity in Figure \ref{fig:Bsum}, where $r$ is limited to the realm of real variables, from above and below along the real line.
\end{rem}

We also note that \eqref{BImId} can be recast by making use of the simple identity
\begin{equation}
\sinh(\pi v)=\cosh(\pi v)-\exp(-\pi v)
\label{SinCoshId}
\end{equation}
in which case \eqref{BImId} becomes
\begin{equation}
\int_{0}^{\infty} \beta_{I} \! \left(1/2+i\,v \right) d v
 = \frac{\ln \! \left(2\right)}{2}-B\,,
\label{BImId2}
\end{equation}
where
\begin{equation}
B\equiv\int_{0}^{\infty}\frac{ \!\beta_{I} \! \left(\frac{1}{2}+i\,v \right) {\mathrm e}^{-\pi \,v}}{\cosh \! \left(\pi \,v \right)}d v\,.
\label{BId}
\end{equation}
The numerical value of \eqref{BId} is given in Table \ref{Tab:Cdefs} of Appendix \ref{sec:Consts}; the analytic value of \eqref{BImId2} defies the same analysis as that employed in the previous Section (\ref{sec:RealAn}).


\subsubsection{Method B: Choose $\alpha=\pi/4$ } \label{sec:MremB}

Here, we consider a second possibility leading to an unequivocal result in a different way. In  \eqref{C1Icp}  and \eqref{C1Rdp} set $\alpha=\pi/4$, and add, to eliminate the divergent term, yielding a simpler form of \eqref{Ans12Lb2}. Let $q=0$, to obtain
\begin{equation}
\int_{0}^{\infty}{\Re \! \left(M_{h}^{i\,v}\,\eta \! \left(1/2+i\,v ,1/2\right)\right) +\Im \! \left(M_{h}^{i\,v}\,\eta \! \left(1/2+i\,v , 1/2\right)\right) \tanh \! \left(\pi \,v \right)}{}d v
 = \sqrt{M_{h}}\,\eta \! \left(1, m +1\right)\,.
\label{T1A}
\end{equation}
Then repeat the same sequence with $\alpha=3\pi/4$, to obtain
\begin{equation}
\int_{0}^{\infty}{\Im \! \left(M^{i\,v}\,\eta \! \left(1/2+i\,v , 1/2\right)\right) \tanh \! \left(\pi \,v \right)-\Re \! \left(M_{h}^{i\,v}\,\eta \! \left(1/2+i\,v , 1/2\right)\right) }{}d v
 = \sqrt{M_{h}}\,\eta \! \left(1, m +1\right)\,.
\label{T2A}
\end{equation}

Adding and subtracting \eqref{T1A} and \eqref{T2A} will reproduce \eqref{Sc1} and \eqref{Sc2}. Because of the consistency among the various means of evaluation, one must conclude that removal of the singularity is the ``preferred", but not necessarily ``unique", choice when analyzing \eqref{BetaReId} and \eqref{BetaReAv}.

\subsubsection{Method C: The case $q=1$}


Instead of evaluating the various identities at the limit point defined by $q=0$, we can consider the same identity at the point $q=1,~ \sin(\alpha)=0$, so the (red) circle in Figure \ref{fig:Diagram} passes through the two nearest neighbour singular points (where $q=0$ and $m:=m\pm1$), the centre of expansion has shifted from the origin to $(-m-1/2,i0)$ and no series term is singular. With $m=0$, from \eqref{C1Icp} we find the particularly simple result
 \begin{equation}
 \int_{0}^{\infty}\eta_{R} \! \left(1/2+i\,v , 1/2\right) d v =0\,
\label{QA2}
 \end{equation}
 
and from \eqref{C1Rdp} we find

\begin{equation}
\int_{0}^{\infty}\eta_{I} \! \left(1/2+i\,v , 1/2\right) \tanh \! \left(\pi \,v \right)d v
 =\,2+ \pi/2\,.
 \label{QB2}
\end{equation}
Compare with \eqref{Sc1} and \eqref{Sc2}, corresponding to $m=1$.

\subsection{Eq. \eqref{Ans12Lb2} and the exceptional case(s) $w=1/2$} \label{sec:PhiPiwh}

In the case $w=1/2$, the exceptional point lies at $a=(-m-1,i0)$, in which case \eqref{C1Icp} and \eqref{C1Rdp} are virtually unchanged except that we employ $M_{u}\equiv m+1$ in place of $M_{h}$ and we are dealing with $\eta(1/2+iv)$, the alternating Zeta function, and indirectly Riemann's function $\zeta(1/2+iv)$. In the following, having studied the intricacies of indeterminacy in the previous section, only results will be given. The reader is encouraged to repeat these variations as an exercise.

To summarize, in what  follows, I use $\alpha=\pi/4$ and $\alpha=3\pi/4$ throughout, with $q\rightarrow 0$ in the appropriate combination of identities equivalent to \eqref{Ans12Lb2} with the substitutions given above, in which case indeterminism has been banished. Further, identities involving $\eta(s)$ are usually transcribed in terms of $\zeta(s)$ (see \eqref{etLim1}), since integrals written in the latter form are of greater general interest. The simplest identity found is

\begin{equation}
\int_{0}^{\infty}\Re \! \left(\left(m +1\right)^{-1/2+i\,v}\,\eta \! \left(1/2+i\,v \right)\right)d v
 = 0
\label{AzPlus}
\end{equation}
and
\begin{equation}
\int_{0}^{\infty}\Im \! \left(\left(m +1\right)^{-1/2+i\,v}\,\eta \! \left(1/2+i\,v \right)\right) \tanh \! \left(\pi \,v \right)d v
 = \eta \! \left(1, m +1\right)\,.
\label{AzMinus}
\end{equation}

\begin{rem} Prior to setting $q=0$, the combined identities were tested numerically for reasonable values of $0.5\lesssim q$.
\end{rem}
Special cases now follow:
\begin{itemize} 

\item{In \eqref{AzPlus}, put $m=0$, then $m=1$ along with \eqref{etLim1}, to obtain}

\begin{equation}
\int_{0}^{\infty} \left(\left(1-\sqrt{2}\,\cos \! \left(v\,\ln \! \left(2\right)\right)\right)\zeta_{R} \! \left(1/2+i\,v \right)-\sqrt{2}\,\sin \! \left(v\,\ln \! \left(2\right)\right) \zeta_{I} \! \left(1/2+i\,v \right)\right) d v=0
\label{S0}
\end{equation}
and
\begin{equation}
\int_{0}^{\infty}\left( \left(\sqrt{2}\,\cos \! \left(v\,\ln \! \left(2\right)\right)-2\right) \zeta_{R} \! \left(1/2+i\,v \right)-\sqrt{2}\,\sin \! \left(v\,\ln \! \left(2\right)\right) \zeta_{I} \! \left(1/2+i\,v \right)\right) d v=0\,.
\label{S1}
\end{equation}
After adding and subtracting we get several identities relating the components of $\zeta(1/2+iv)$:

\begin{equation}
\int_{0}^{\infty}\left(3
-
{2  \sqrt{2}} \cos \! \left(v\ln \! \left(2\right)\right)\right)\,\zeta_{R} \! \left(1/2+i\,v \right)d v =0
\label{s01ma}
\end{equation}
and
\begin{equation}
\int_{0}^{\infty}\left(\zeta_{R} \! \left(1/2+i\,v \right)
 +2\,\sqrt{2} \sin \! \left(v\,\ln \! \left(2\right)\right)\zeta_{I} \! \left(1/2+i\,v \right)\right) d v=0\,, 
\label{s01pa}
\end{equation}
or
\begin{equation}
\int_{0}^{\infty}\left(\cos \! \left(v\ln \! \left(2\right)\right)\,\zeta_{R} \! \left(1/2+i\,v \right) 
+
3 \sin \! \left(v\ln \! \left(2\right)\right) \zeta_{I} \! \left(1/2+i\,v \right)\right) d v =0\,.
\label{spboth}
\end{equation}

\item{From \eqref{AzMinus} we similarly find the general form}

\begin{align} \nonumber
\int_{0}^{\infty}&\left(\left(-\sqrt{2}\,\sin \! \left(v\,\ln \! \left({M_{u}}/{2}\right)\right)+\sin \! \left(v\,\ln \! \left(M_{u}\right)\right)\right) \zeta_{R} \! \left(1/2+i\,v \right) \right. \\ & \left. \nonumber
+\left(-\sqrt{2}\,\cos \! \left(v\,\ln \! \left(M_{u}/2\right)\right)+\cos \! \left(v\,\ln \! \left(m +1\right)\right)\right) \zeta_{I} \! \left(1/2+i\,v \right)\right) \tanh \! \left(\pi \,v \right)d v \\
 &= \sqrt{M_{u}}\,\eta \! \left(1, M_{u}\right)
\label{Azm}
\end{align}

so, for $m=0$ and $m=1$ respectively,
\begin{align} \nonumber
\int_{0}^{\infty}\left(\sqrt{2}\,\sin \! \left(v\,\ln \! \left(2\right)\right) \zeta_{R} \! \left(1/2+i\,v \right)\right.&\left.
+\left(1-\sqrt{2}\,\cos \! \left(v\,\ln \! \left(2\right)\right)\right) \zeta_{I} \! \left(1/2+i\,v \right)\right) \tanh \! \left(\pi \,v \right)d v\\&
 =\ln \! \left(2\right)\,,
\label{Wm0}
\end{align}
and
\begin{align} \nonumber
\int_{0}^{\infty}&\left(\sqrt{2}\,\sin \! \left(v\,\ln \! \left(2\right)\right) \zeta_{R} \! \left(1/2+i\,v \right)+ \left(\sqrt{2}\,\cos \! \left(v\,\ln \! \left(2\right)\right)-2\right)\zeta_{I} \! \left(1/2+i\,v \right)\right)\,\tanh \! \left(\pi \,v \right) \,d v \\
&\hspace{3cm} = 2-2\,\ln \! \left(2\right)\,.
\label{Wm1}
\end{align}

\item{Adding and subtracting then gives}
\begin{equation}
\int_{0}^{\infty}\left(2\,\sqrt{2}\,\sin \! \left(v\,\ln \! \left(2\right)\right) \zeta_{R} \! \left(1/2+i\,v \right)-\zeta_{I} \! \left(1/2+i\,v \right)\right) \tanh \! \left(\pi \,v \right)d v
 = 2-\ln \! \left(2\right)
\label{Az01p}
\end{equation}
and
\begin{equation}
\int_{0}^{\infty}\left(3-2\,\sqrt{2}\,\cos \! \left(v\,\ln \! \left(2\right)\right)\right) \zeta_{I} \! \left(1/2+i\,v \right)\,\tanh \! \left(\pi \,v \right)d v
 = 3\,\ln \! \left(2\right)-2\,.
\label{Az01m}
\end{equation}

\item{With reference to \eqref{SinCoshId}, \eqref{Wm0} and \eqref{Wm1} can be rewritten in the form}
\begin{align} \nonumber
\int_{0}^{\infty}\left(\sqrt{2} \sin \right.&\left. \! \left(v\,\ln \! \left(2\right)\right) \zeta_{R} \! \left(1/2+i\,v \right)
 + 
\left(1-\sqrt{2}\,\cos \! \left(v\,\ln \! \left(2\right)\right)\right) \zeta_{I} \! \left(1/2+i\,v \right)\right)\,d v \\=&
\ln \! \left(2\right)+Z_{0}
\label{ZA0}
\end{align}
and
\begin{align} \nonumber
\sqrt{2}\int_{0}^{\infty}( \sin \!& \left(v\,\ln \! \left(2\right)\right) \zeta_{R} \! \left(1/2+i\,v \right)
  -  \left(\sqrt{2}-\cos \! \left(v\,\ln \! \left(2\right)\right)\right) \zeta_{I} \! \left(1/2+i\,v \right))\,d v\\&= 2-2\,\ln \! \left(2\right)+Z_{1}+Z_{2}
\label{ZA1}
\end{align}
where
\begin{equation}
Z_{0}\equiv\int_{0}^{\infty}\frac{{\mathrm e}^{-\pi \,v}}{\cosh \! \left(\pi \,v \right)}\,\,\eta_{I} \! \left(1/2+i\,v\right)d v\,,
\label{Z0def}
\end{equation}
\begin{equation}
Z_{1}\equiv\sqrt{2} \int_{0}^{\infty}\frac{{\mathrm e}^{-\pi \,v} \left(\cos \! \left(v\,\ln \! \left(2\right)\right)-\sqrt{2}\right) }{\cosh \! \left(\pi \,v \right)}\zeta_{I} \! \left(1/2+i\,v \right)d v
\label{Z1def}
\end{equation}
and 
\begin{equation}
Z_{2}\equiv \sqrt{2} \int_{0}^{\infty}\frac{{\mathrm e}^{-\pi \,v}\,\sin \! \left(v\,\ln \! \left(2\right)\right) }{\cosh \! \left(\pi \,v \right)}\zeta_{R} \! \left(1/2+i\,v \right)d v \,.
\label{Z2def}
\end{equation}

Numerical values for the constants $Z_{j},j=0..2$ are listed in Appendix \ref{sec:Consts}, Table \ref{Tab:Cdefs}.
 
\begin{rem}Except for the latter three, it is not clear whether any of these integrals are convergent, or obtain their value by analytic continuation (see Appendix \ref{sec:Defs} , Lindel\"of's Hypothesis), since the conjectured asymptotic nature of $\zeta(1/2+i\,v)$ remains unproven. 
\end{rem}

\item{Adding and subtracting \eqref{s01ma} and \eqref{s01pa} yields the identities}
\begin{equation}
\int_{0}^{\infty}\left(1
 -
\frac{2\,\sqrt{2}}{3}\cos \! \left(v\,\ln \! \left(2\right)\right)\right) \zeta_{I} \! \left(1/2+i\,v \right)d v= \, -\frac{2}{3}+\ln \! \left(2\right)-\frac{Z_{2}}{3}-\frac{Z_{1}}{3}+\frac{Z_{0}}{3}
\label{Zans1p}
\end{equation}
and
\begin{equation}
\int_{0}^{\infty}\left(\zeta_{I} \! \left(1/2+i\,v \right) - 
2\,\sqrt{2}\sin \! \left(v\,\ln \! \left(2\right)\right) \zeta_{R} \! \left(1/2+i\,v \right)\right)d v=\,-2+\ln \! \left(2\right)-Z_{2}-Z_{1}-Z_{0}\,,
\label{Zans2p}
\end{equation}
from which we obtain
\begin{align} \nonumber
\int_{0}^{\infty}&\left(\sin \! \left(v\,\ln \! \left(2\right)\right) \zeta_{R} \! \left(1/2+i\,v \right)
- 
\frac{1}{3}\cos \! \left(v\,\ln \! \left(2\right)\right) \zeta_{I} \! \left(1/2+i\,v \right)\right) d v\\
&  =\frac{\sqrt{2}\,Z_{0}}{3}+\frac{\sqrt{2}\,Z_{1}}{6}+\frac{\sqrt{2}\,Z_{2}}{6}+\frac{\sqrt{2}}{3}\,.
\label{Zans12p}
\end{align}
\item{As a numerical check, \eqref{Zans1p} predicts that} 
\begin{equation}
\int_{0}^{\infty}\left(1-\frac{2\,\sqrt{2}\,\cos \! \left(v\,\ln \! \left(2\right)\right)}{3}\right) \zeta_{I} \! \left(1/2+i\,v \right)d v
 = 0.017324632\cdots.
\label{V12p}
\end{equation}
\end{itemize}

\begin{figure}[h] 
\centering
\includegraphics[width=0.7\textwidth,height=.5\textwidth]{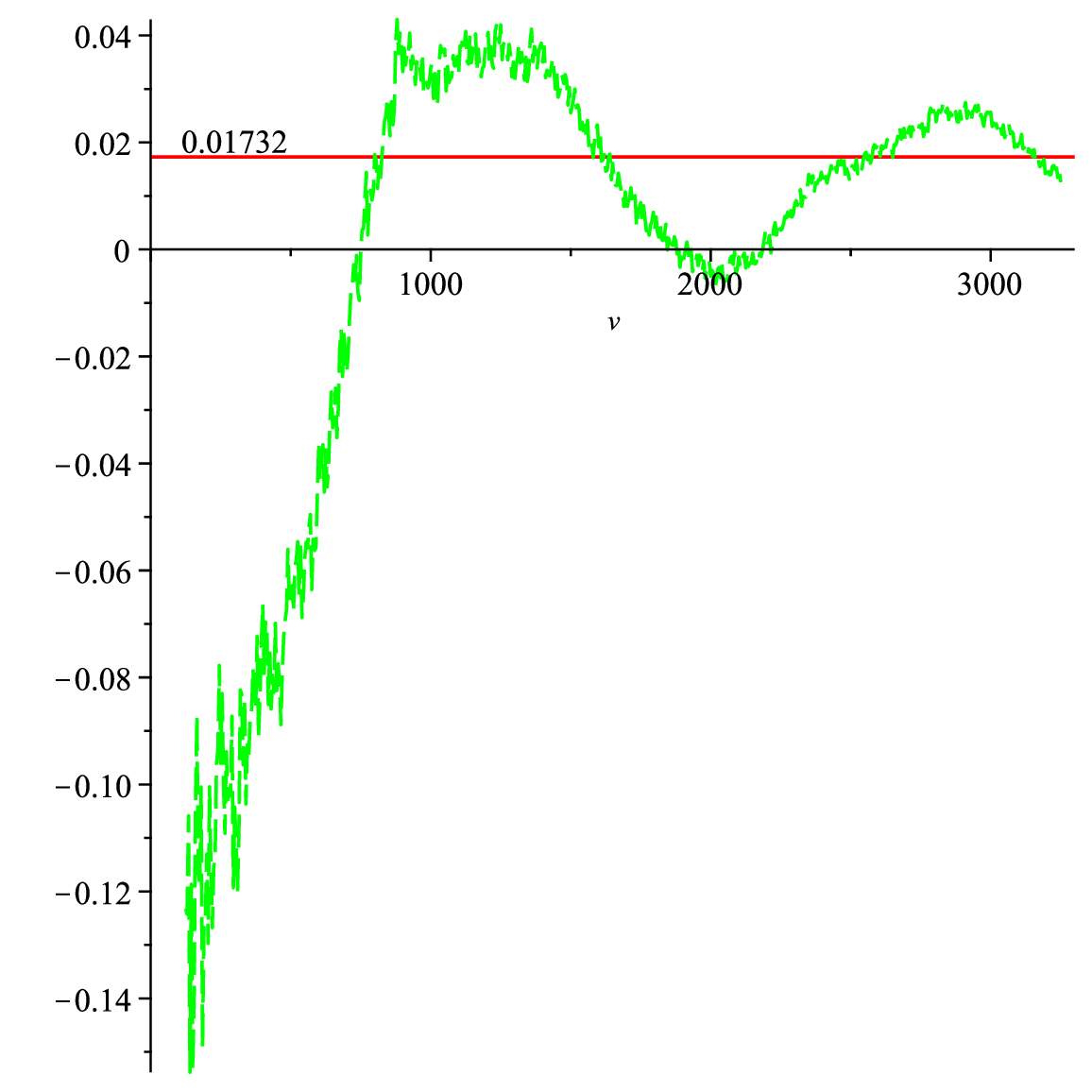}
\caption{Values of the average of the running partial sums of pieces of \eqref{V12p} in comparison to its predicted value. Numerical integration for $v\gtrsim 3000$ confounds both Maple and Mathematica.}
\label{fig:Cesaro2}
\end{figure}

It is interesting to attempt to verify this identity numerically, notwithstanding the fact that it possesses a highly oscillatory integrand and may not converge, by partitioning the integral into small pieces, numerically integrating each small piece, and averaging the running partial sums. The result is shown in Figure \ref{fig:Cesaro2}. This is an example of either Cesàro summation or Cesàro regularization as the case may be \cite[page 194]{Widder}, depending on whether the integral converges numerically or not. The outcome of the experiment shown in the figure, suggests that \eqref{V12p} is plausible and numerical convergence is a true possibility. Elsewhere, it has been shown that  Cesàro summation \cite[Appendix A]{MilHughMMM2025} is a very effective way to numerically approximate (or regularize) integrals such as \eqref{V12p}. 

\section{The exceptional cases Corollary \ref{sec:Cor2}, $\phi=\pi$} \label{sec:Except_zeta}
For the remainder of this Section, we shall be considering the case $\phi\rightarrow \pi$ and $r\rightarrow w+1/2+ m$ based on Corollary \ref{sec:Cor2}, i.e. \eqref{M1ng2}. 

\subsection{ The Real part}

\begin{itemize}
\item{First, let $w=0$ so we are interested in the case that $r\rightarrow M_{h}\equiv m+1/2$.}

Consider the right-hand side of \eqref{Ming2ca} labelled $R_{1}$, where we again find that an ambiguity resides. This time, for full generality, let $\phi=\pi-\epsilon,~r=M_{h}+\rho \mbox{ and } \sigma=1/2+\delta$. The leading terms, as each of $\delta\rightarrow 0,~ \rho\rightarrow 0 \mbox{ and } \epsilon\rightarrow 0$, independently is


\begin{equation}
R_{1}\sim \sqrt{M_{h}}\left(\frac{\left(1+\frac{M_{h}}{\rho}\right) \epsilon}{\rho}-\pi  \left(1+\frac{\delta}{\rho}\right)\right)\,.
\label{ALL}
\end{equation}
This demonstrates the indeterminism that we have come to expect, and, for the first time, also shows that the previous first level choice $\sigma=1/2$ ($\delta=0$) has spared us from a higher order of indeterminism, even in the case that $\epsilon=0$. To begin, choose $\phi=\pi$ (i.e. $\epsilon=0$) followed by $\sigma=1/2$, (i.e. $\delta=0$), leaving a finite result, that being


\begin{equation}
R_{1}{ \,\underset{\rho,\epsilon,\delta\rightarrow 0} \lim}\,\underset{\looparrowleft}=\,-\pi\sqrt{M_{h}}
\label{varLim1}
\end{equation}
where the ordering of limits is to be read from right to left. From \eqref{LhsRb} and \eqref{varLim1} with $m=0$, we obtain the identity
\begin{equation}
\int_{0}^{\infty}\left(\left(\sqrt{2}-\cos \! \left(v\,\ln \! \left(2\right)\right)\right)\zeta_{R} \! \left(1/2+i\,v \right) -\sin \! \left(v\,\ln \! \left(2\right)\right)\zeta_{I} \! \left(1/2+i\,v \right) \right)d v
 \underset{\looparrowleft\,}=\, -\frac{\pi \,\sqrt{2}}{2}\,.
\label{Car1}
\end{equation}
Combining \eqref{Car1} and \eqref{S1} yields
\begin{equation}
\int_{0}^{\infty} \left(\sqrt{2}\,\cos \! \left(v\,\ln \! \left(2\right)\right)-2\right)\zeta_{R} \! \left(1/2+i\,v \right)d v
 \underset{\looparrowleft\,}=\, \frac{\pi}{2},
\label{CbR1}
\end{equation}
which, together with \eqref{CR2} reproduces \eqref{CR5}, a non-exceptional result, demonstrating both consistency and indicating that the ordering in \eqref{varLim1} was the preferred choice.
 
The case for general values of $M_{h}$ and general values of $\phi$ is fairly lengthy and is left as an exercise, but the result is again as expected. In the case that we set $\phi=\pi $ and $ r=m+1/2\equiv M_{h}$, \eqref{Rparta} becomes

\begin{align} \nonumber
&\int_{0}^{\infty}\left(\frac{\left(-C \! \left(v , \sigma , r\right)\,\zeta_{R} \! \left(i\,v +\sigma \right)+S \! \left(v , \sigma , r\right)\,\zeta_{I} \! \left(i\,v +\sigma \right)\right) \sin \! \left(\pi \,\sigma \right) \cosh^{2}\left(\pi \,v \right)}{\cos \! \left(2\,\pi \,\sigma \right)-\cosh \! \left(2\,\pi \,v \right)}\right. \\ & \left.  \nonumber
+\frac{\left(-C \! \left(v , \sigma , r\right)\,\zeta_{I} \! \left(i\,v +\sigma \right)-S \! \left(v , \sigma , r\right)\,\zeta_{R} \! \left(i\,v +\sigma \right)\right) \cos \! \left(\pi \,\sigma \right) \sinh \! \left(\pi \,v \right) \cosh \! \left(\pi \,v \right)}{\cos \! \left(2\,\pi \,\sigma \right)-\cosh \! \left(2\,\pi \,v \right)}\right)d v \\ \nonumber
& =\frac{1}{2} 
\left(\left(-{\pi \,\sin \! \left(\pi \,\sigma \right)}{}+\left({\psi \! \left(m +1\right)}{}-{\ln \! \left(M_{h} \right)}{}\right) \cos \! \left(\pi \,\sigma \right)\right) M_{h} +\left({1}{}-{\sigma}{}\right) \cos \! \left(\pi \,\sigma \right)\right) M_{h}^{-\sigma}\\
&+\frac{M_{h}^{1-\sigma}\,\cos \! \left(\pi \,\sigma \right)}{2(\,r -M_{h})}.
\label{Way4}
\end{align}

Effectively, the right-hand side is of the form anticipated from \eqref{ALL}, replete with a singularity that vanishes if $\sigma=1/2$. So, letting $\sigma=1/2$, yields
\begin{align} \nonumber
\int_{0}^{\infty}&\left(\left(\sqrt{2}\,\cos \! \left(v\,\ln \! \left(2\,M_{h} \right)\right)-\cos \! \left(v\,\ln \! \left(M_{h} \right)\right)\right) \zeta_{R} \! \left(1/2+i\,v \right) \right. \\ & \left. 
-\left(\sqrt{2}\,\sin \! \left(v\,\ln \! \left(2\,M_{h} \right)\right)-\sin \! \left(v\,\ln \! \left(M_{h} \right)\right)\right) \zeta_{I} \! \left(1/2+i\,v \right)\right)d v
 = -\pi \,\sqrt{M_{h}}
\label{Way4a}
\end{align}
which reduces to \eqref{Car1} if $M_{h}=1/2$. The identity \eqref{Way4} suggests that the singularity can be removed by combining two instances of that identity. For example set $M_{h}=m_{1}+1/2$ then $M_{h}=m_{2}+1/2$ and subtract with suitable multipliers. Alternatively use $\sigma$ and $1-\sigma$ and add; $\sigma=1/4$ and $\sigma=3/4$ is also a propitious choice. However, these variations lead to very long identities that show no promise of reducing to something simple, so are left as exercises for the reader.

It is also interesting to digress slightly by evaluating the same identity \eqref{LhsRb} in a different order. If $w=0,~\sigma=1/2,~ m=0 \mbox{ and }r=1/2$, the equivalent identity becomes

\begin{align} \nonumber
&\int_{0}^{\infty}\left(-{\left(\cos \! \left(v\,\ln \! \left(2\right)\right)-\sqrt{2}\right)  \zeta_{R} \! \left(1/2+i\,v \right)}{}-{\sin \! \left(v\,\ln \! \left(2\right)\right)\zeta_{I} \! \left(1/2+i\,v \right)}{}\right)\frac{\cosh \! \left(v\,\phi \right)}{\cosh \! \left(\pi \,v \right)}d v\\ \nonumber
& = 
\frac{\sqrt{2}\,\sin \! \left(\frac{\phi}{2}\right) \sin \! \left(\phi \right) }{2}\overset{\infty}{\underset{j =1}{\sum}}\, \frac{1}{\left(1+2\,j \right) \cos \left(\phi \right)+2\,j^{2}+2\,j +1}\\ \nonumber
&-\frac{ \sqrt{2} \cos \! \left(\frac{\phi}{2}\right)}{2}\overset{\infty}{\underset{j =1}{\sum}}\, \frac{\left(j +1\right) \left(\cos \left(\phi \right)+1\right)}{\left(2\,j^{2}+j \right) \cos \left(\phi \right)+2\,j^{3}+2\,j^{2}+j}\\ 
&+\frac{\sqrt{2}\,\cos \! \left(\frac{\phi}{2}\right) \gamma}{2}-\frac{\sqrt{2}\,\cos \! \left(\frac{\phi}{2}\right) \ln \! \left(2\right)}{2}-\frac{\sqrt{2}\,\cos \! \left(\frac{\phi}{2}\right) \left(\phi \,\sin \! \left(\frac{\phi}{2}\right) \cos \! \left(\frac{\phi}{2}\right)-1\right)}{\cos \! \left(\phi \right)+1}
\label{Way2a}
\end{align}
which was easily numerically verified for small values of $\phi$. However, in the limit $\phi\rightarrow \pi$, the right-hand side of \eqref{Way2a} diverges as predicted by \eqref{ALL}; specifically we find 
\begin{equation}
\int_{0}^{\infty}\left(\sqrt{2}-\cos \! \left(v\,\ln \! \left(2\right)\right)\right) \zeta_{R} \! \left(1/2+i\,v \right)-\sin \! \left(v\,\ln \! \left(2\right)\right) \zeta_{I} \! \left(1/2+i\,v \right)d v
\rightarrow\,- \frac{\sqrt{2}}{\phi -\pi}-\frac{\pi \,\sqrt{2}}{2}
\label{Wex1}
\end{equation}
which should be compared to \eqref{Car1}.


\item{If $w=1/2,~\sigma=1/2$, indeterminism again arises in \eqref{Ming2e}:}

Let $\sigma=1/2,~w=1/2,~\mbox{and }\phi=\pi$ to obtain
\begin{equation}
\int_{0}^{\infty}\sin \! \left(v\,\ln \! \left(r \right)\right) \zeta_{I} \! \left(1/2+i\,v \right)-\cos \! \left(v\,\ln \! \left(r \right)\right) \zeta_{R} \! \left(1/2+i\,v \right)d v  =\,\pi\, \sqrt{r}\,,
\label{MingH2}
\end{equation}
an identity that is valid for all values of $r$, particularly $r=m+1$ and does not appear to be indeterminate. However, considering the oppositely reordered  directional limit, let $\sigma=1/2,~w=1/2,~r=m+1$, followed by $\phi= \pi-\epsilon$ and, let $\epsilon\rightarrow 0$ to obtain 
\begin{align} \nonumber
\int_{0}^{\infty}\sin \! \left(v\,\ln \! \left(m+1 \right)\right) \zeta_{I} \! \left(1/2+i\,v \right)&-\cos \! \left(v\,\ln \! \left(m+1 \right)\right) \zeta_{R} \! \left(1/2+i\,v \right)d v  \\
&=\,- \frac{1}{\sqrt{m+1}}\,\epsilon^{-1}+\sqrt{m+1}\,\pi
\label{MingH3}
\end{align}
another indeterminism. However, choosing \eqref{MingH2} produces
\begin{equation}
\int_{0}^{\infty}\sin \! \left(v\,\ln \! \left(m +1\right)\right) \zeta_{I} \! \left(1/2+i\,v \right)-\cos \! \left(v\,\ln \! \left(m +1\right)\right) \zeta_{R} \! \left(1/2+i\,v \right)d v
\underset{\looparrowleft}=\, \pi\,\sqrt{m +1}\,.
\label{MingWha}
\end{equation}
By choosing \eqref{MingH2} it is possible to obtain an identity for the left-hand side at a constant value of $r$, valid for all $\phi$. Being rather lengthy, (but numerically verifiable for $|\phi|<<\pi$), this identity is recorded as \eqref{MingH1bX} in Appendix \ref{sec:LongEqs}.  If $m=0$, \eqref{MingWha} reduces to \eqref{CR5}. If $m=1$, \eqref{MingWha} becomes
\begin{equation}
\int_{0}^{\infty}\left(-\sin \! \left(v\,\ln \! \left(2\right)\right) \zeta_{I} \! \left(1/2+i\,v \right)+\cos \! \left(v\,\ln \! \left(2\right)\right) \zeta_{R} \! \left(1/2+i\,v \right)\right)d v
\underset{\looparrowleft}=\, -\pi \,\sqrt{2}\,,
\label{MgM1}
\end{equation}
which can be combined with several of the results previously obtained to reproduce others already obtained. So, \eqref{MgM1}, is not an independent identity, although it was independently derived, and bears the tag of indeterminism because an alternative evaluation exists.
\end{itemize}

\subsection{The imaginary part}
\begin{itemize}

\item{In the case that $w=0$ and $\sigma=1/2$, let $R_2$ denote the right-hand side of \eqref{LhsIa}}. 

If we first set $\sigma=1/2$, a choice emerges. For the case $m=0$, with the limits to be read from right to left, and $\rho,~\epsilon$ being as defined previously (see \eqref{ALL}), we have
\begin{equation}
\underset{\rho \rightarrow 0,\epsilon\rightarrow 0}{\mathrm{lim}}\;R_{2}=-\frac{1+\left(\ln \! \left(2\right)-\gamma \right) \sqrt{2}}{2}-\frac{\sqrt{2}}{2\,\rho}
\end{equation}
whereas
\begin{equation}
\underset{\epsilon\rightarrow 0,\rho \rightarrow 0}{\mathrm{lim}}\;R_{2}=-\frac{\left(\ln \! \left(2\right)-\gamma \right) \sqrt{2}}{2} +O(\rho/\epsilon^2)
\label{LimUsed}
\end{equation}
a far more desirable, but not unique, outcome . From the preferential choice \eqref{LimUsed}  we then have


\begin{align} \nonumber
\int_{0}^{\infty}&\left(-\sin \! \left(v\,\ln \! \left(2\right)\right) \zeta_{R} \! \left(1/2+i\,v \right)+\left(\cos \! \left(v\,\ln \! \left(2\right)\right)-\sqrt{2}\right) \zeta_{I} \! \left(1/2+i\,v \right)\right) \tanh \! \left(\pi\, v\,\right)\,d v
 \\&\hspace{30pt} \underset{\looparrowleft}=
\frac{\left(-\ln \! \left(2\right)+\gamma \right) \sqrt{2}}{2}\,.
\label{Limw0}
\end{align}
From previous results, it is possible to check this choice for consistency. Combine \eqref{Limw0} with \eqref{Wdm1} to obtain
\begin{align}
\int_{0}^{\infty}&\left(\cos \! \left(v\,\ln \! \left(2\right)\right) \zeta_{I} \! \left(1/2+i\,v \right)-\sin \! \left(v\,\ln \! \left(2\right)\right) \zeta_{R} \! \left(1/2+i\,v \right)\right)\tanh \! \left(\pi \,v \right) d v\\&
\underset{\looparrowleft}=
-\frac{\sqrt{2}\,\gamma}{2}-\frac{\sqrt{2}\,\ln \! \left(2\right)}{2}\,.
\label{Limw0A}
\end{align}
Adding \eqref{Limw0A} and \eqref{Wdp} yields
\begin{equation}
\int_{0}^{\infty} \cos \! \left(v\,\ln \! \left(2\right)\right) \zeta_{I} \! \left(1/2+i\,v \right)\,\tanh \! \left(\pi \,v \right)d v
 \underset{\looparrowleft}=\, 
-\frac{3\,\sqrt{2}\,\gamma}{4}-\frac{3\,\sqrt{2}\,\ln \! \left(2\right)}{4}+\frac{\sqrt{2}}{2}
\label{Limw0B}
\end{equation}
and combining \eqref{Limw0B} with \eqref{Wdm} reproduces \eqref{Wdm1}, demonstrating consistency. It is also possible to show that the right-hand side of \eqref{Limw0} is the only possible value that will reproduce this sequence of equalities, proving both that the above argument is not circular, and, at the same time, that \eqref{Limw0} is not independent of  the other identities, although it was derived independently. Finally, on the same ordering basis as \eqref{LimUsed}, for $m>0$ and $M_{h}=m+1/2$, we find


\begin{align} \nonumber
\int_{0}^{\infty}&\left(\left(-\sqrt{2}\,\sin \! \left(v\,\ln \! \left(2\,M_{h} \right)\right)+\sin \! \left(v\,\ln \! \left(M_{h} \right)\right)\right) \zeta_{R} \! \left(\frac{1}{2}+i\,v \right) \right. \\ & \left. \nonumber
+\left(-\sqrt{2}\,\cos \! \left(v\,\ln \! \left(2\,M_{h} \right)\right)+\cos \! \left(v\,\ln \! \left(M_{h} \right)\right)\right) \zeta_{I} \! \left(\frac{1}{2}+i\,v \right)\right) \tanh \! \left(\pi \,v \right)d v\\
& \underset{\looparrowleft}=
\left(-\psi \! \left(m \right)+\frac{1}{2\,M_{h}\,m}+\ln \! \left(M_{h}\right)\right) \sqrt{M_{h}}\,,\hspace{1cm} m>0.
\label{Limw0M}
\end{align}


If $m=1$, \eqref{Limw0M} becomes
\begin{align} \nonumber
\int_{0}^{\infty}&\left(\left(\sqrt{2}\,\cos \! \left(v\,\ln \! \left(3\right)\right)-\cos \! \left(v\,\ln \! \left({2}/{3}\right)\right)\right) \zeta_{I} \! \left({1}/{2}+i\,v \right) \right. \\ \nonumber &  \left.
+\left(\sqrt{2}\,\sin \! \left(v\,\ln \! \left(3\right)\right)+\sin \! \left(v\,\ln \! \left({2}/{3}\right)\right)\right) \zeta_{R} \! \left({1}/{2}+i\,v \right)\right) \tanh \! \left(\pi \,v \right)d v
\\& \underset{\looparrowleft}= 
\frac{\left(3\,\ln \! \left(\frac{2}{3}\right)-3\,\gamma -1\right) \sqrt{6}}{6}
\label{Limw1}
\end{align}


\item{In the case $w=1/2$, define $R_{3}$ to represent the right-hand side of \eqref{C2I}.} With $w=1/2~,M_{u}\equiv~ m+1$ and $ \sigma=1/2$ we have two possible limits, similar to those given in \eqref{LimUsed}:

\begin{equation}
\underset{r \rightarrow M_{u}+\rho,\phi\rightarrow \pi-\epsilon}{\mathrm{lim}}R_{3}=\,-\frac{\sqrt{M_{u}}}{\rho}-\frac{1}{2\,\sqrt{M_{u}}}+\left(-\psi \! \left(M_{u}\right)+\ln \! \left(M_{u}\right)\right) \sqrt{M_{u}}+O\! \left(\epsilon\right)+O\! \left(\rho \right),
\label{O1}
\end{equation}
whereas an alternative is
\begin{equation}
\underset{\phi\rightarrow \pi-\epsilon,r \rightarrow M_{u}+\rho}{\mathrm{lim}}R_{3}=-\sqrt{M_u}\,\psi \! \left(M_u\right)+\sqrt{M_u}\,\ln \! \left(M_u\right)+O(\epsilon)+O(\rho).
\label{O2}
\end{equation}
Notice that the non-singular part of \eqref{O1} differs from the leading terms of \eqref{O2}. In keeping with accumulating experience, we choose to employ \eqref{O2} rather than attempting to eliminate the singular term in \eqref{O1}, in which case, \eqref{C2I} becomes
\begin{align} \nonumber
\int_{0}^{\infty}&\left(\cos \! \left(v\,\ln \! \left(M_u \right)\right) \zeta_{I} \! \left(1/2+i\,v \right)+\sin \! \left(v\,\ln \! \left(M_u \right)\right) \zeta_{R} \! \left(1/2+i\,v \right)\right) \tanh \! \left(\pi \,v \right)d v\\
 &\hspace{3cm}\underset{\looparrowleft\,}= 
\sqrt{M_u} \left(\psi \! \left(M_u \right)-\ln \! \left(M_u \right)\right)\,,
\label{Oa3}
\end{align}
reducing to \eqref{Wdm1} when $m=0$, suggesting that \eqref{O2} is indeed the preferred choice. If $m=1$, we find
\begin{align} \nonumber
\int_{0}^{\infty}&\left(\cos \! \left(v\,\ln \! \left(2\right)\right) \zeta_{I} \! \left(1/2+i\,v \right)+\sin \! \left(v\,\ln \! \left(2\right)\right) \zeta_{R} \! \left(1/2+i\,v \right)\right) \tanh \! \left(\pi \,v \right)d v \\
&\hspace{3cm} \underset{\looparrowleft\,}= \sqrt{2} \left(1-\gamma -\ln \! \left(2\right)\right)\,.
\label{Oa4}
\end{align}
Finally, adding  \eqref{Limw0} and \eqref{Oa4} yields the identity
\begin{equation}
\int_{0}^{\infty}\left(2\,\cos \! \left(v\,\ln \! \left(2\right)\right)-\sqrt{2}\,\right) \zeta_{I} \! \left(1/2+i\,v \right)\,\tanh \! \left(\pi \,v \right)d v
 \underset{\looparrowleft\,}= 
\sqrt{2} \left(1-\frac{3\,\sqrt{2}\,\ln \! \left(2\right)}{2}-\frac{\gamma}{2}\right)\,.
\label{Oa4b}
\end{equation}
Combining \eqref{Oa4b} and \eqref{Wdm2} reproduces \eqref{Wdm1}, again demonstrating consistency in the choice \eqref{O2}.
\end{itemize}

\section{A Principal Value Analogue} \label{sec:Pv}
Many of the identities quoted so far started as functions of the single complex variable $a=re^{i\phi}$ and morphed into  meromorphic  functions of the real variables $r$ and $\phi$ considered independently. So, by the ``Edge of the Wedge" theorem \cite[pages 251 ff. ]{Vlad}, which states that, with sufficient conditions on analyticity, analytic continuation is valid along a line corresponding to the real part of one of several complex variables, this is sufficient to allow both $r$ and $\phi$ to be treated as independent and real variables as has been done  throughout. Now consider an extension such that $\phi$ becomes complex.

If one writes $\phi:=\phi+i\theta$ and simultaneously $r:=re^{\theta}$, the variable  $a=re^{i\phi}$ is unchanged. However, now let $\phi:=\phi+i\theta$, with no change in $r$. This is equivalent to adding an imaginary ($\theta$) axis perpendicular to the plane of Figure \ref{fig:Diagram}. Since each of the original identities develops a new real and imaginary part, this generates a new set of identities, the left-hand sides of which are easily written, the right-hand sides, lengthy and not so easily reported (this is  left as an exercise for the reader with access to a computer algebra programme). With $w=1/2$, the four integrals that arise are the following:

\begin{align}
J_{1}(r,\phi,\theta)\equiv \int_{0}^{\infty}\frac{\left(\sin \! \left(v\,\ln \! \left(r \right)\right) \eta_{I} \! \left(\frac{1}{2}+i\,v\right)-\cos \! \left(v\,\ln \! \left(r \right)\right) \eta_{R} \! \left(\frac{1}{2}+i\,v \right)\right)}{\cosh \! \left(\pi \,v \right)} \cosh \! \left(\phi \,v \right) \cos \! \left(v\,\theta \right) d v\,,
\label{J1Int}\\
J_{2}(r,\phi,\theta)\equiv \int_{0}^{\infty}\frac{\left(\sin \! \left(v\,\ln \! \left(r \right)\right) \eta_{I} \! \left(\frac{1}{2}+i\,v\right)-\cos \! \left(v\,\ln \! \left(r \right)\right) \eta_{R} \! \left(\frac{1}{2}+i\,v\right)\right) }{\cosh \! \left(\pi \,v \right)}\sinh \! \left(\phi \,v \right) \sin \! \left(v\,\theta \right) d v\,,
\label{J2Int}\\
J_{3}(r,\phi,\theta)\equiv  \int_{0}^{\infty}\frac{\left(\cos \! \left(v\,\ln \! \left(r \right)\right) \eta_{I} \! \left(\frac{1}{2}+i\,v \right)+\sin \! \left(v\,\ln \! \left(r \right)\right) \eta_{R} \! \left(\frac{1}{2}+i\,v\right)\right)}{\cosh \! \left(\pi \,v \right)}  \sinh \! \left(\phi \,v \right) \cos \! \left(v\,\theta \right) d v \,,
\label{J3Int}\\
J_{4}(r,\phi,\theta)\equiv \int_{0}^{\infty}\frac{\left(\cos \! \left(v\,\ln \! \left(r \right)\right) \eta_{I} \! \left(\frac{1}{2}+i\,v \right)+\sin \! \left(v\,\ln \! \left(r \right)\right) \eta_{R} \! \left(\frac{1}{2}+i\,v\right)\right)}{\cosh \! \left(\pi \,v \right)} \cosh \! \left(\phi \,v \right) \sin \! \left(v\,\theta \right) d v \,.
\label{J4Int}
\end{align}
In special cases, the right-hand sides simplify considerably, so, for the record the following presents the corresponding non-zero identities when $\phi=0$:
\begin{align} 
	&\begin{aligned}
J{_1}(r,0,\theta)=\frac{ r^{\frac{3}{2}}}{2}\left(\cosh \! \left(\frac{3\,\theta}{2}\right)-\cosh \! \left(\frac{\theta}{2}\right)\right)& \overset{\infty}{\underset{j =0}{\sum}}\, \frac{\left(-1\right)^{j}}{\left(2\,j +2\right) r\,\cosh \left(\theta \right)+r^{2}+\left(j +1\right)^{2}}\\
 - \sqrt{r}\,\cosh \! \left(\frac{\theta}{2}\right)&\overset{\infty}{\underset{j =0}{\sum}}\, \frac{\left(-1\right)^{j} \left(r\,\cosh \! \left(\theta \right)+j +1\right)}{\left(2\,j +2\right) r\,\cosh \! \left(\theta \right)+r^{2}+\left(j +1\right)^{2}}\,,
\label{J1_phi0}
	\end{aligned}
\\
&\begin{aligned}
 J_{4}(r,0,\theta)=r^{3/2}\cosh \! \left(\frac{\theta}{2}\right) \overset{\infty}{\underset{j =0}{\sum}}&\,\, \frac{ \left(-1\right)^{j}\,\sinh \! \left(\theta \right)}{r^{2}+\left(j +1\right)^{2}+2\,r\,\cosh \! \left(\theta \right) \left(j +1\right)}\\
 -\sqrt{r}\sinh \! \left(\frac{\theta}{2}\right) \overset{\infty}{\underset{j =0}{\sum}}&\,\, \frac{\left(-1\right)^{j} \left(r\,\cosh \! \left(\theta \right)+j +1\right)}{2\,r\,\cosh \! \left(\theta \right) \left(j +1\right)+j^{2}+r^{2}+2\,j +1}\,.
\label{J4_phi0}
\end{aligned}
\end{align}
For $r,\theta>0$, there are no singularities in these cases, all of which have been verified numerically. 

\subsection{The case $\phi=\pm\pi$} \label{sec:PvalExcpt}

We now consider the case $\phi=\pi$ with $\theta\neq0$, with the eventual goal of studying the limit $\theta\rightarrow 0$, thereby approaching the pathological singularity from a different direction. If $r\neq(1+m)\exp({\pm \theta})$ -- this defines exceptional cases where one indexed term in the denominator vanishes -- we find the non-zero identities
\begin{align} \nonumber
\\
&J_{2}(r,\pi,\theta)=\sqrt{r}\,\sinh \! \left(\frac{\theta}{2}\right)\overset{\infty}{\underset{j =0}{\sum}}\, \frac{\left(-1\right)^{j} \left(j +1+r \right)}{r^{2}+\left(j +1\right)^{2}-2\,r\,\left(j +1\right)\cosh \! \left(\theta \right) }\,,
\label{HR0I}\\
&J_{3}(r,\pi,\theta)=\,-\sqrt{r}\,\cosh \! \left(\frac{\theta}{2}\right)\overset{\infty}{\underset{j =0}{\sum}}\, \frac{\left(-1\right)^{j} \left(j +1-r \right)}{r^{2}+\left(j +1\right)^{2}-2\,r \left(j +1\right)\cosh \! \left(\theta \right)}\,,
\label{HI0R}
\end{align}

From these exceptional cases, arises another possibility to deal with the ambiguous singularity discussed throughout. To demonstrate the possibilities, here I consider only the case $J_{2}((m+1)\exp(\pm\theta),\pi,\theta)$ in detail. In \eqref{J2Int} set $r=(m+1)\,\exp(\pm\theta)$ followed by $\phi=\pm\pi\mp\epsilon$; then evaluate the limit $\epsilon\rightarrow 0$, to find, in either case
\begin{align} \nonumber
&\int_{0}^{\infty}\left(\sin \! \left(v\,\ln \! \left({\mathrm e}^{\theta} \left(1+m \right)\right)\right) \eta_{I} \! \left(1/2+i\,v \right)\right.\\&\left.\nonumber \hspace{50pt}-\cos \! \left(v\,\ln \! \left({\mathrm e}^{\theta} \left(1+m \right)\right)\right) \eta_{R} \! \left(1/2+i\,v \right)\right) \tanh \! \left(\pi v \right) \sin \! \left(v\,\theta \right)d v\\
& = 
\frac{\sqrt{1+m} \left({\mathrm e}^{\theta}-1\right) }{2}\overset{\infty}{\underset{\substack{j =0\\j\neq m}}{\sum}}\! \frac{\left(-1\right)^{j} \left(j +1+{\mathrm e}^{\theta} \left(1+m \right)\right)}{\left(j -m \right) \left(\left(-m -1\right) {\mathrm e}^{2\,\theta}+j +1\right)}-\frac{\left(-1\right)^{m}\,{\mathrm e}^{\theta}}{2\,\sqrt{1+m} \left({\mathrm e}^{2\,\theta}-1\right)}\,.
\label{H2A}
\end{align}
This result is valid for all $\theta$, and diverges with a simple pole if $\theta= 0$, where it reduces to the results obtained previously, but becomes considerably more complicated if $\theta=\ln(k)/2$ where $k$ is a positive integer. Alternatively, repeat the exercise but in the opposite order. First, in \eqref{J2Int}, let $\phi=\pi$ followed by $r=(m+1+\rho)\,\exp(\pm\theta)$, and in the limit $\rho\rightarrow 0$, we find
\begin{align} \nonumber
\int_{0}^{\infty}&\left(\sin \! \left(v\,\ln \! \left({\mathrm e}^{\pm\theta} \left(1+m \right)\right)\right) \eta_{I} \! \left(1/2+iv \right)\right.\\\left.\nonumber -\right.&\left.\cos \! \left(v\,\ln \! \left({\mathrm e}^{\pm\theta} \left(1+m \right)\right)\right) \eta_{R} \! \left(1/2+iv \right)\right) \tanh \! \left(\pi v \right) \sin \! \left(v \theta \right)d v\\
 &\sim\mp\frac{\left(-1\right)^{m}\,\sqrt{1+m}}{2\,\rho}+O\! \left(\rho^{0}\right)
\label{HR1IRB}
\end{align}
reigniting the indeterminism that has pervaded this work. Since the singular term in \eqref{HR1IRB} changes sign depending on which of the two polarities of $\theta$ one chooses, it becomes possible to remove the singularity by adding the two versions of \eqref{HR1IRB} corresponding to $+$ and $-$ signs. This is analogous to a Principal Value regularization, where a contour along a real axis is deformed above and below that axis into the imaginary part of the complex plane to avoid a singularity (e.g. Appendix \ref{sec:Reg}). After performing the addition, and following considerable simplification (Maple) we find
\begin{align} \nonumber
&\int_{0}^{\infty}\left(\eta_{I} \! \left(1/2+i\,v \right)\,\sin \! \left(v\,\ln \! \left(1+m \right)\right)-\eta_{R} \! \left(1/2+i\,v \right)\,\cos \! \left(v\,\ln \! \left(1+m \right)\right)\right)\,\tanh \! \left(\pi \,v \right) \sin \! \left(2\,v\,\theta \right) \,d v\\
& = 
\frac{1}{2}\sqrt{1+m} \left({\mathrm e}^{2\,\theta}-1\right) {\mathrm e}^{\theta} \overset{\infty}{\underset{\substack{j =0\\j\neq m}}{\sum}}\! \frac{\left(-1\right)^{j} \left(j +m +2\right)}{\left(\left(j +1\right) {\mathrm e}^{2\,\theta}-m -1\right) \left(j +1-\left(m +1\right) {\mathrm e}^{2\,\theta}\right)}-\frac{\left(-1\right)^{m}\,{\mathrm e}^{\theta}}{\sqrt{1+m} \left({\mathrm e}^{2\,\theta}-1\right)}\,.
\label{H2B}
\end{align} 
In this variation, no limits were applied, so the result is deterministically exact, although numerically unverifiable. Note that the left-hand sides of \eqref{H2A} and \eqref{H2B} differ, but they share pole terms at $\theta=0$ that differ in magnitude by a factor of two. Thus subtracting \eqref{H2A} and \eqref{H2B} with a suitable multiplier, yields a new result, equivalent to one of the identities obtained previously using the same variables. In fact for the case $m=0$ we find,
\begin{align} \nonumber
&\int_{0}^{\infty}\sin^{2}\left(v\,\theta \right) \eta_{I} \! \left({1}/{2}+i\,v \right)\,\tanh \! \left(\pi \,v \right)d v
\\& = 
-\frac{{\mathrm e}^{2\,\theta}\,\sinh \! \left(\theta \right)}{2}\overset{\infty}{\underset{j =0}{\sum}}\! \frac{\left(-1\right)^{j} \left(j +3\right)}{\left(\left(j +2\right) {\mathrm e}^{2\,\theta}-1\right) \left(2-{\mathrm e}^{2\,\theta}+j \right)}
+
\frac{1}{2}(e^{\theta}-1)
 \overset{\infty}{\underset{j =0}{\sum}}\! \frac{\left(-1\right)^{j} \left(j +2+{\mathrm e}^{\theta}\right)}{\left(j +1\right) \left(2-{\mathrm e}^{2\,\theta}+j \right)}\,.
\label{HAB2}
\end{align}
In Appendix \ref{sec:LongEqs}, the general case ($m>0$) is listed (see \eqref{H2AB1}). Because the singularity at $\theta=0$ has now been removed, by taking the limit $\theta\rightarrow 0$ on the right-hand side of \eqref{HAB2} and expanding the left-hand side about $\theta=0$, we recover \eqref{Diff2}; also, since both sides can be expressed in even powers of $\theta$, the next moment corresponding to the coefficients of $\theta^4$ gives the identity

\begin{equation}
\int_{0}^{\infty}v^{4}\,\eta_{I} \! \left(1/2+i\,v \right)\,\tanh \! \left(\pi \,v \right)d v
 = 
-\frac{\ln \! \left(2\right)}{16}-\frac{5\,\pi^{2}}{12}-\frac{87\,\zeta \! \left(3\right)}{4}-\frac{7\,\pi^{4}}{15}-\frac{45\,\zeta \! \left(5\right)}{2}\,,
\label{HAB4}
\end{equation}
in agreement with a higher order differentiation of \eqref{Diff2}. 

\begin{rem}
The expansion of the left-hand side of \eqref{HAB2} says nothing about the odd power moments of \eqref{Diff2} that were conjectured to vanish.
\end{rem}



\section{Summary} \label{sec:Summary}

In a fairly lengthy exposition, I have evaluated a number of improper integrals containing Riemann's zeta and related functions and identified a pathology which, in many respects, emulates an essential singularity. Many of these integral identities are new; others reproduce known identities derived elsewhere. The pathology arises because the analysis is based on the use of two independent real variables, which, as noted, do not obey the usual rules of complex analysis. However (see Section (\ref{sec:RealAn})), by examining the related case of a real function of one real variable, we come to realize that the pathology is in fact the generalization of a commonly encountered discontinuity in one real variable, extended into the complex plane, or, in reverse, the pathology studied here in in the complex plane, degenerates into a run-of-the-mill discontinuity when projected onto the real line. This requires that we distinguish a \textbf{directed} limit in real analysis, where a limit point can only be approached from one of two directions, from a \textbf{directional} limit, where a limit point can be approached from a variety of directions.

In order to arrive at a useful calculus embodying the directional limit, several different methods were explored, each of which succeeded in arriving at a meaningful and self-consistent result. As noted, and to reiterate, in some cases it was possible to compare two identities, one of which was calculated according to the accepted methods of analysis, the other by means of the method in which the pathology was circumvented. In all cases the comparison yielded consistency. In one instance, a numerical experiment was attempted in order to check the plausibility of a regularized result. In order to focus on the identified pathology, many threads were left hanging. A few possibilities that merit future study include the application of the results reported here to:

\begin{itemize}

\item{\bf Indeterminism: }
By removing the pathology, a self-consistent calculus arose involving finite quantities, but the possibility remains that a different choice could have been made resulting in singular quantities. Is it possible that this second choice will yield a self-consistent calculus and, if so, what does it mean? Notice the difference between the regular parts of \eqref{O1} and \eqref{O2}.

\item{\bf Off the Critical Line: }
Throughout, for brevity, most examples were chosen to lie on the critical line $\sigma=1/2$, although the methodology allows for a more general choice. In a few places - see \eqref{ALL} and \eqref{Way4} -  there were suggestions that choosing to evaluate an integral along the critical line creates a new level of indeterminism that may distinguish the critical line from elsewhere in the critical strip. Is this significant?

\item{\bf Integral transforms: }
Many of the results quoted here effectively define the coefficients of the cosine and sine transforms of the function $\zeta(\sigma+it/ln(M))$ - see \eqref{WhI1} for example. Is it possible to express $\zeta(\sigma+it)$ as a Fourier or Laplace transform in closed form - see \cite[Section 8.2]{IvicZ}, other than via \eqref{Cry1A} and known transformations among Mellin and other transforms - see \cite{Widder}?

\item{\bf Moments: }
Several samples were given of the higher moments for the integrals being studied - see \eqref{Hiworh} for example. Is there a general form? Does it lead to new identities?

\item{\bf Mean Value Formulae: } In \cite[Section 11.4]{BoChRoWe}, Ivi{\'c} studies the asymptotic limit of the integral

\begin{equation}
\int_{0}^{^T}|\zeta(1/2+it)|^{2k}dt\,.
\label{Iv1}
\end{equation}
Can \eqref{Iv1} be related to integrals studied here, perhaps by extending $T:=T+i\tau$?

\item{\bf By Parts: } Does integration by parts applied to any of the identities presented here lead to new identities?

\item{\bf Extension: } There is no requirement that the variable $b$ in Section \ref{sec:Cor3a} be real. What happens if it becomes complex?
\item{\bf Principal Value and Odd Moments:}
In Section \ref{sec:Pv}, in the interests of brevity, many interesting possibilities were omitted. These might form the basis of a future study, including the conjecture that all odd moments vanish (see \eqref{Diff2} and \eqref{HAB4}).

\item{\bf Application to the Riemann Hypothesis:} Many identities are known and labelled as ``equivalent to RH" -- e.g. \cite{BoChRoWe} and \cite{Spigler} -- but few of them involve the evaluation of integrals. The notable exceptions are the Volchkov equivalence, (which has been shown to be both incorrect and only-half equivalent even if it were correct - see \cite[Sections 8 and 9]{Milgram_Exploring}), and Lindel{\"o}f's Hypothesis (LH) (see Appendix \ref{sec:Defs}) which governs the convergence of integrals such as those being considered here (see also \cite{2012arXiv1201.2633F}). Thus  any successful evaluation of an integral such as 
\begin{equation}
\label{Zint}
\int_{0}^{\infty}|\zeta(1/2+it)| dt
\end{equation}
will yield information about its convergence and hence the asymptotic nature of its integrand, i.e. LH. It is generally agreed that LH is a consequence of RH, but not the converse (see \cite[page 328]{Titch2}, \cite[page 65]{Widder}, although Ivi{\'c} is not so sure \cite[page 138]{BoChRoWe}) and Fokas \cite{Fokas}, citing earlier work, points out that ``the validity of Lindelöf’s hypothesis reduces significantly the number of possible zeros that disobey the Riemann hypothesis". However, if LH is proven false, then RH must be false. 

Let us now speculate that the integral \eqref{Zint}, which bears a strong similarity to integrals discussed here, shares the property of being indeterminate. That means that an analyst might someday prove, not necessarily using the methods discussed here, that the value of that integral is infinite, consistent with the truth of LH and therefore likely RH. However, a second analyst could similarly prove that the value of the integral \eqref{Zint} is finite, reflecting the indeterminism discussed here, and thereby disprove both LH and RH. 

A resolution of these opposing viewpoints resides in  the simple possibility that both camps are correct if \eqref{Zint} is indeterminate as conjectured, in which case one could conclude that RH is ambiguously both true and false and therefore unprovable, according to the usual rules of mathematical logic.  
%

\end{itemize}

\section{Acknowledgements}

\section{Declarations}
{\bf: Funding:} All expenses related to this work have been borne by the author.

{\bf: Conflicts of interest} None.

{\bf: Availability of data and material} All relevant material is contained herein.

{\bf: Code availability}  N/A

{\bf: Authors' contributions} This is the work of the sole author.

{\bf: Ethics approval} N/A

{\bf: Consent to participate} N/A

{\bf: Consent for publication} N/A

\bibliographystyle{unsrt}

\bibliography{Intrep_12b.bib}

@book{Paris&Kaminski,
title={Asymptotics and Mellin-Barnes Integrals},
author={R.B. Paris and D. Kaminski},
series={Encyclopedia of Mathematics and its Applications},
volume=85,
year=2001,
publisher={Cambridge University Press},
address={Cambridge, U.K.}
}

@book{IvicZ,
author={Aleksandar Ivi{\'c}},
title="The theory of {H}ardy's Z-Function",
series={Cambridge tracts in Mathematics},
year=2013,
publisher ={Cambridge University Press},
address={New York},
volume=196
}

@book{Titch2,
author="Titchmarsh,E.C.
and Heath-Brown, D.R",
title={The Theory of the {R}iemann {Z}eta-{F}unction},
edition={{S}econd},
publisher={Oxford Science Publications},
address={Oxford},
year={1986},
}

@book{Nielsen,
title={Handbuch der Theorie der Gammafunktion},
author={Niels Nielsen},
year={1906},
publisher={Druck und Verlag von B.G. Teubner},
address={Leipzig},
note={reprinted from the original by Wentworth Press},
}

@article{HuKimKim,
title="Special values and Integral representations for the {H}urwitz-type {E}uler {Z}eta Functions",
author={Su Hu, Daeyeoul Kim and Min-Soo Kim},
journal={J. Korean Math. Soc.,},
volume=55,
number=1,
pages="185-210",
year=2018,
note={https://doi.org/10.4134/JKMS.j170110}
}

@article{Esp&Moll2002,
author={Espinosa, Olivier and Moll, Victor H.},
year=2002,
title={{O}n Some Integrals Involving the {H}urwitz {Z}eta Function: Part 1},
journal= "The Ramanujan Journal",
volume=6,
pages="159-188",
url="https://doi.org/10.1023/A:1015706300169",
note="also available from https://doi.org/10.1023/A:1015706300169 or https://arxiv.org/0012078v1"
}

@article{Shpot&Paris,
author={{Shpot}, M.~A. and {Paris}, R.~B.},
title={Integrals of products of {H}urwitz zeta functions via {F}eynman parametrization and two double sums of {R}iemann zeta functions},
year=2020,
month=may,
day=24,
note={available from https://arxiv.org/pdf/1609.05658},
}

@book{NIST,
      editor = "F.~W.~J. Olver and D.~W. Lozier and R.~F. Boisvert and C.~W. Clark",
       title = "{NIST Handbook of Mathematical Functions}",
   publisher = "Cambridge University Press",
     address = "New York, NY",
        year = "2010",
        note = "Print companion to \cite{NIST:DLMF}"}

@manual{Math,
title={Mathematica, version 12},
organization={Wolfram Research},
url={http://wolfram.com},
year="2020",
address="Champaign, Illinois",}

@manual{Maple,
title="Maple.",
organization={Maplesoft, a division of Waterloo Maple Inc., version 2021},
}

@article{MilHughMMM2025,
author={Michael Milgram and Roy J. Hughes},
title={On a generalized moment integral containing {R}iemann's {Z}eta function: {A}nalysis and experiment},
journal={Modern Mathematical Methods},
pages={14-41},
volume=3,
number=1,
year={2025},
note={ https://doi.org/10.64700/mmm.48}
}

@book{G&R,
author= {I.S. Gradshteyn and I.M. Ryzhik},
title={Tables of {I}ntegrals, {S}eries and {P}roducts, corrected and enlarged Edition},
publisher= "Academic Press", 
year=1980,
}

@book{prudnikov,
  title={Integrals and Series: More Special functions},
  author={A.P. Prudnikov and {\relax Yu.} A. Brychkov and O.I. Marichev},
  isbn={9782881240904},
  lccn={85027065},
  url={https://books.google.ca/books?id=2t2cNs00aTgC},
  year={1986},
  volume={3},
  publisher={Gordon and Breach Science Publishers},
  address={New York},
}

@misc{connon2012integral,
      title={On an integral involving the digamma function}, 
      author={Donal F. Connon},
      year={2012},
      eprint={1212.1432},
      archivePrefix={arXiv},
      primaryClass={math.GM},
      note={available from https://arxiv.org/abs/1212.1432}
}

@article{CoppoAndCandel,
title="A note on some formulae related to {E}uler sums",
author="Marc-Antoine Coppo and Bernard Candelpergher",
year=2021,
journal="hal-03170892v4",
note="available from https://hal.univ-cotedazur.fr/hal-03170892v4",
}

@book{Vlad,
author={{V}asiliy {S}ergeyevich {V}ladimirov},
title="Methods of the Theory of Functions of Many Complex Variables",
publisher={Dover Publications Inc},
address={Mineola, New York},
year={2007},
note={Translated from the Russian by MIT Press, 1966},
}

@book{CounterX,
title={Counterexamples in Analysis},
author={Gelbaum Bernard R. and Olmsted John M.H.},
publisher={Holden-Day Inc},
address={San Francisco},
year=1964,
}

@book{Hardy,
title={Pure Mathematics},
author={Hardy,G.H.},
publisher={Cambridge at the University Press},
year=1945,
}

@book{WW,
title={A Course of Modern Analysis},
author={Whittaker E.T. and Watson G.N.},
publisher={Cambridge University Press},
year=1950,
}

@book{M&F,
title={Methods of Theoretical Physics, Part 1},
author={Morse P.M. and Feshbach H.},
publisher={McGraw-Hill Book Company Inc.},
year=1953,
address={Toronto},
}

@manual{Maple24,
title="Maple.",
organization={Maplesoft, a division of Waterloo Maple Inc., version 2024},
}

@manual{Math23,
title={Mathematica, version 13.2},
organization={Wolfram Research},
url={http://wolfram.com},
year="2023",
address="Champaign, Illinois",}

@book{Widder,
title={An Introduction to Transform Theory},
author={Widder D.V.},
publisher={Academic press},
address={New York},
year=1971,
}

@book{BoChRoWe,
   title = "The {R}iemann {H}ypothesis : a resource for the afficionado and virtuoso alike",
   author = "{Borwein P.,} {Choi S.,} {Rooney B.,} and Weirathmueller A. ",
   series = "CMS books in mathematics",
   publisher = "Springer",
   address = "New York",
   isbn = "978-0-387-72125-5",
   year = 2008
}

@article{Spigler,
author = {Spigler, Renato},
title = {A Brief Survey on the {R}iemann Hypothesis and Some Attempts to Prove It},
journal = {Symmetry},
volume = {17},
year = {2025},
number = {2},
page = {225},
URL = {https://www.mdpi.com/2073-8994/17/2/225},
note={https://doi.org/10.3390/sym17020225},
}

@article{Milgram_Exploring,
author = {Michael Milgram},
title = {Exploring {R}iemann's functional equation},
journal = {Cogent Mathematics},
volume = {3},
number = {1},
pages = {1179246},
year = {2016},
doi = {10.1080/23311835.2016.1179246},

note = { 
        http://dx.doi.org/10.1080/23311835.2016.1179246
    
},
eprint = { 
        http://dx.doi.org/10.1080/23311835.2016.1179246
    
}
}

@ARTICLE{2012arXiv1201.2633F,
       author = {{Fokas}, A.~S. and {Lenells}, J.},
        title = "{On the Asymptotics to all Orders of the Riemann Zeta Function and of a
        Two-Parameter Generalization of the Riemann Zeta Function}",
      journal = {Memoires of the American Mathematical Society to be published,https://www.ams.org/cgi-bin/mstrack/accepted},
         year = 2015,
        month = Dec,
          eid = {arXiv:1201.2633},
        note = {Available  from ArXiv.org as arXiv:1201.2633v2},
archivePrefix = {arXiv},
       eprint = {1201.2633},
 primaryClass = {math.NT},
       adsurl = {https://ui.adsabs.harvard.edu/#abs/2012arXiv1201.2633F}
}

@article{Fokas,
author={A.S.Fokas},
title={A novel approach to the {L}indel{\"o}f hypothesis},
journal ="Trans. Math. and Appl.",
volume=3,
number=1,
month=Feb,
pages={1-49},
note={https://doi.org/10.1093/imatrm/tnz006},
year=2019,
}

@book{TitchF,
title={The Theory of Functions, Second Edition},
author={Titchmarsh, E.C.},
publisher={Oxford University Press},
year={1939, reprinted and corrected 1949},
}

@misc{Analytic,
author={Weisstein, E.W.},
year={1999-2001},
title={Analytic Continuation. MathWorld -- A Wolfram Web Resource.},
note={Retrieved from https://mathworld.wolfram.com/AnalyticContinuation.html},
}

@article{HuKim2024,
title = {Asymptotic expansions for the alternating {H}urwitz {Z}eta function and its derivatives},
journal = {J. Math. Anal. and Appl.},
author = {Su Hu and Min-Soo Kim},
volume = {537},
number = {1},
pages = {128306},
year = {2024},
note = {https://doi.org/10.1016/j.jmaa.2024.128306},
}

@book{EulerIV,
title={Institutiones Calculi Integralia},
author={Euler L.},
volume={IV},
year=1794,
}

\begin{appendices}

\section{Definitions and Identities} \label{sec:Defs}
\subsection{Definitions}
 By way of review:

\begin{Def}
{\bf Holomorphic function:} ``A function $f(z)$ is said to be holomorphic at a point $z^{0}\in\mathfrak{C}^{n}$ if, in some neighbourhood of $z^{0}$, it is the sum of an absolutely convergent power series $f(z)={\underset{n\ge0}\sum}\,a_{n}\left(z-z_{0}\right)^{n}$. \cite[Section I.4.1]{Vlad} - Weierstrass' defintion".
\end{Def}

\begin{Def}
{\bf Domain of holomorphy:} ``To every holomorphic function there corresponds a unique domain of holomorphy, namely, the domain of existence of the function" \cite[Section I.8.6]{Vlad}".
\end{Def}

\begin{Def}
{\bf Singular points:} ``The boundary points of a domain of holomorphy of a function are called the singular points of that function" \cite[Section I.8.7]{Vlad}".
\end{Def}

\begin{Def}
{\bf Essential Singularity:} ``In the neighbourhood of an isolated essential singularity, a one-valued function takes every value, with one possible exception, an infinity of times" \cite[Picard's Second Theorem, page 283] {TitchF}. \label{Picard} \newline
\end{Def}

\begin{Def} \label{sec:Lind}
{\bf Lindel{\"o}f's Hypothesis} ``The Lindel{\"o}f Hypothesis is that $\zeta(\sigma+it)=O(t^\epsilon)$ for every positive $\epsilon$ and every $\sigma\geq \frac{1}{2}$" \cite[Chapter XIII]{Titch2}.
\end{Def}

\begin{Def} \label{sec:AnalCont}
{\bf Analytic Continuation:} ``Let $f_{1}$ and $f_{2}$ be analytic functions on domains $\Omega_{1}$ and  $\Omega_{2}$, respectively, and suppose that the intersection  $\Omega_{1}\cap \Omega_{2}$ is not empty and that $f_{1}=f_{2}$ on $\Omega_{1} \cap \Omega_{2}$. Then $f_{2}$ is called an analytic continuation of $f_{1}$ to $\Omega_{2}$, and vice versa. Moreover, if it exists, the analytic continuation of $f_{1}$ to $\Omega_{2}$ is unique." \cite{Analytic}
\end{Def}

\subsection{Identities and Lemmas}

\begin{itemize}

\item{By dissecting the sum \eqref{HzetaAlt} into its even and odd parts it is easily shown that}

\begin{equation}
\eta(s,a)=\left( \zeta(s,a/2)-\zeta(s,a/2+1/2)\right)/2^s.
\label{EtaZeta}
\end{equation}

\item{From \cite[Eq. 25.11.43] {NIST} and \cite{HuKim2024}, we quote the asymptotic limit}
\begin{equation}
\zeta(s,a)\sim a^{1-s}/(s-1)+\frac{1}{2}a^{-s}+\dots\,
\label{Zasyv}
\end{equation}
as $a\rightarrow\infty$, from which, using \eqref{EtaZeta} we also find
\begin{equation}
\eta \! \left(s ,a \right) \sim 
\frac{a^{-s}}{2}
\label{EtaAsy}
\end{equation}
in the same limit.

\item{From \cite[Section 13.21]{WW}}
\begin{equation}
\underset{b \rightarrow 0}{\mathrm{lim}}\! \left(\zeta \left(b +1, v\right)-\frac{1}{b}\right)
 = -\psi \! \left(v \right)\,.
\label{Zlim0}
\end{equation}

\item

From \cite[Eqs. (5.4.3), (5.4.4) and (5.7.6) respectively]{NIST} we have
\begin{equation}
{| \Gamma \! \left(i\,v \right)|}^{2} = 
\frac{\pi}{v\,\sinh \! \left(\pi \,v \right)}\,,
\label{Nist5p4p3}
\end{equation}

\begin{equation}
{| \Gamma \! \left(1/2+i\,v \right)|}^{2} = 
\frac{\pi}{\cosh \! \left(\pi \,v \right)}\,
\label{Nist5p4p4}
\end{equation}
and
\begin{equation}
\psi \! \left(z \right) = 
-\gamma -\frac{1}{z}+\overset{\infty}{\underset{j =1}{\sum}}\! \frac{z}{j \left(j +z \right)}\hspace{12pt} z\neq 0,-1\dots\,.
\label{PsiZ}
\end{equation}
\item
Differentiating \eqref{Cry1a} with respect to $z$ yields
\begin{align} \nonumber
`\int_{0}^{\infty}v^{z -1}\,\ln \! \left(v \right) \eta \! \left(1, v +1\right)d v
 &= 
\frac{\pi  \left(\left(2^{z}-1\right) \zeta^{\,\prime}\! \left(1-z \right)-2^{z}\,\ln \! \left(2\right) \zeta \! \left(1-z \right)\right)}{\sin \! \left(\pi \,z \right)}\\
&+\frac{\pi^{2}\,\cos \! \left(\pi \,z \right) \left(2^{z}-1\right) \zeta \! \left(1-z \right)}{\sin \! \left(\pi \,z \right)^{2}}\,.
\label{T2p2a}
\end{align}

\item
From Maple \cite{Maple} and verified numerically

\begin{equation}
\int_{1/2}^{1}\zeta\, \left(b +1, v\right)d v
 = \frac{\zeta \! \left(b \right) \left(2^{b}-2\right)}{b}\,.
\label{Zident}
\end{equation}

\item 
Based on \cite[Eq. (5.7.6)]{NIST} we separate the real and imaginary parts of a special case of the digamma function:

\begin{align} \nonumber
\psi_{R} \! \left(r\,\cos \! \left(\phi \right)+w +\frac{1}{2}+i\,r\,\sin \! \left(\phi \right)\right)
 = 
-\gamma -\frac{r\,\cos \! \left(\phi \right)+w +\frac{1}{2}}{2 \left(w +\frac{1}{2}\right) r\,\cos \! \left(\phi \right)+r^{2}+\left(w +\frac{1}{2}\right)^{2}}\\
+\,\overset{\infty}{\underset{j =1}{\sum}}\,\, \frac{r \left(j +2\,w +1\right) \cos \! \left(\phi \right)+r^{2}+\left(w +\frac{1}{2}\right) \left(j +w +\frac{1}{2}\right)}{j \left(2\,r \left(j +w +\frac{1}{2}\right) \cos \! \left(\phi \right)+j^{2}+r^{2}+\left(w +\frac{1}{2}\right) \left(2\,j +w +\frac{1}{2}\right)\right)}
\label{PsR}
\end{align}
and
\begin{equation}
 \psi_{I} \! \left(r\,\cos \! \left(\phi \right)+w +1/2+i\,r\,\sin \! \left(\phi \right)\right)
 = 
 \overset{\infty}{\underset{j =0}{\sum}}\,\,\frac{r\,\sin \! \left(\phi \right)}{\,2 \left(w +j +\frac{1}{2}\right) r\,\cos \! \left(\phi \right)+r^{2}+\left(w +j +\frac{1}{2}\right)^{2}}\,.
\label{PsI}
\end{equation}

\end{itemize}

\section{Long Equations} \label{sec:LongEqs}

Here we collect many of the most general results used throughout to obtain special case identities.
\begin{itemize}

\item{ The following two identities present the general form of \eqref{HR} and \eqref{HI} for arbitrary values of $\sigma$ as discussed in Section \ref{sec:GenCon1}.}
\begin{align} \nonumber
 &{{\int}}_{\!\!\!0}^{\infty}\left( \frac{\left(-\sin \! \left(v \ln \! \left(r \right)\right) \sin \! \left(\pi  \sigma \right) \cosh \! \left(\pi  v \right)+\cos \! \left(v \ln \! \left(r \right)\right) \sinh \! \left(\pi  v \right) \cos \! \left(\pi  \sigma \right)\right) \eta_{\mathrm{I}}\! \left(\sigma +\mathit{iv} ,w +\frac{1}{2}\right)}{\cos \! \left(\pi  \sigma \right)^{2}-\cosh \! \left(\pi  v \right)^{2}} \right. \\ \nonumber
&\left. +\frac{\left(\sin \! \left(v \ln \! \left(r \right)\right) \cos \! \left(\pi  \sigma \right) \sinh \! \left(\pi  v \right)+\cosh \! \left(\pi  v \right) \sin \! \left(\pi  \sigma \right) \cos \! \left(v \ln \! \left(r \right)\right)\right) \eta_{R}\! \left(\sigma +\mathit{iv} ,w +\frac{1}{2}\right)}{\cos \! \left(\pi  \sigma \right)^{2}-\cosh \! \left(\pi  v \right)^{2}}\right) \cosh \! \left(\phi\,  v\right) {d}v\\ \nonumber
 &=r^{1-\sigma} \left(r \sin \! \left(\phi \right)\sin \! \left(\phi  \left(\sigma-1 \right)\right) {{{\sum_{j=0}^{\infty}}}}\, \frac{\left(-1\right)^{j}}{2 \left(j +w +\frac{1}{2}\right) r \cos \! \left(\phi \right)+r^{2}+\left(j +w +\frac{1}{2}\right)^{2}}
 \right. \\
&\left. -{{\cos \! \left(\phi  \left(\sigma-1 \right)\right){\sum_{j=0}^{\infty}}}}\, \frac{\left(-1\right)^{j} \left(r \cos \! \left(\phi \right)+j +w +\frac{1}{2}\right)}{2 \left(j +w +\frac{1}{2}\right) r \cos \! \left(\phi \right)+r^{2}+\left(j +w +\frac{1}{2}\right)^{2}}\right)
\label{C2R}
\end{align}
and
\begin{align} \nonumber
&{{\int}}_{\!\!\!0}^{\infty}\left(\frac{\left(-\sin \! \left(v \ln \! \left(r \right)\right) \cos \! \left(\pi  \sigma \right) \sinh \! \left(\pi  v \right)-\cosh \! \left(\pi  v \right) \sin \! \left(\pi  \sigma \right) \cos \! \left(v \ln \! \left(r \right)\right)\right) \eta_{\mathrm{I}}\! \left(\sigma +\mathit{iv} ,\frac{1}{2}+w \right)}{\cos \! \left(\pi  \sigma \right)^{2}-\cosh \! \left(\pi  v \right)^{2}}\right. \\ \nonumber
& \left.+\frac{\left(-\sin \! \left(v \ln \! \left(r \right)\right) \sin \! \left(\pi  \sigma \right) \cosh \! \left(\pi  v \right)+\cos \! \left(v \ln \! \left(r \right)\right) \sinh \! \left(\pi  v \right) \cos \! \left(\pi  \sigma \right)\right) \eta_{R}\! \left(\sigma +\mathit{iv} ,\frac{1}{2}+w \right)}{\cos \! \left(\pi  \sigma \right)^{2}-\cosh \! \left(\pi  v \right)^{2}}\right) \sinh \! \left(\phi\,  v \right){d}v\\ \nonumber
& =r^{1-\sigma} \left(r\,\cos \! \left(\phi  \left(\sigma -1\right)\right)\sin \! \left(\phi \right) {{{\sum_{j=0}^{\infty}}}}\,\, \frac{\left(-1\right)^{j}}{2 \left(j +w +\frac{1}{2}\right) r \cos \! \left(\phi \right)+r^{2}+\left(j +w +\frac{1}{2}\right)^{2}} \right. \\
&\left . +\sin \! \left(\phi  \left(\sigma -1\right)\right) {{{\sum_{j=0}^{\infty}}}}\,\, \frac{\left(-1\right)^{j} \left(r \cos \! \left(\phi \right)+j +w +\frac{1}{2}\right)}{2 \left(j +w +\frac{1}{2}\right) r \cos \! \left(\phi \right)+r^{2}+\left(j +w +\frac{1}{2}\right)^{2}}\right) \,.
\label{C2I}
\end{align}

\item{ The following is the general case corresponding to \eqref{HAB2} with $M_{u}=m+1$:}

\begin{align} \nonumber
&\int_{0}^{\infty}\left(\left(\sin \! \left(v\,\ln \! \left(M_{u} \right)\right) \sin \! \left(2\,v\,\theta \right)-2\,\sin \! \left(v\,\ln \! \left({\mathrm e}^{\theta}\,M_{u} \right)\right) \sin \! \left(v\,\theta \right)\right) \eta_{I} \! \left(\frac{1}{2}+i\,v \right) 
\right. \\ \nonumber  &\left.
+\left(-\sin \! \left(2\,v\,\theta \right) \cos \! \left(v\,\ln \! \left(M_{u} \right)\right)+2\,\cos \! \left(v\,\ln \! \left({\mathrm e}^{\theta}\,M_{u} \right)\right) \sin \! \left(v\,\theta \right)\right) \eta_{R} \! \left(\frac{1}{2}+i\,v \right)\right) \tanh \! \left(\pi \,v \right)d v
\\ \nonumber & = 
\sqrt{M_{u}} \left(-1\right)^{M_{u}}\overset{\infty}{\underset{j =0}{\sum}}\! \frac{\displaystyle -\frac{{\mathrm e}^{2\,\theta} \left({\mathrm e}^{\theta}-{\mathrm e}^{-\theta}\right) \left(-1\right)^{j} \left(j +2\,M_{u} +1\right)}{2 \left(j +M_{u} +1\right) {\mathrm e}^{2\,\theta}-2\,M_{u}}+\frac{\left({\mathrm e}^{\theta}-1\right) \left(-1\right)^{j} \left(j +M_{u} +1+{\mathrm e}^{\theta}\,M_{u} \right)}{j +1}}{1-M_{u}\,{\mathrm e}^{2\,\theta}+j +M_{u}}\\
&+ \sqrt{M_{u}}\overset{M_{u} -2}{\underset{j =0}{\sum}}\! \frac{\left({\mathrm e}^{\theta}-1\right) \left(-1\right)^{j} \left(j +1+{\mathrm e}^{\theta}\,M_{u} \right)}{\left(j -M_{u} +1\right) \left(-M_{u}\,{\mathrm e}^{2\,\theta}+1+j \right)}\,.
\label{H2AB1}
\end{align}

\item{ The following are the most general forms of the real and imaginary parts of \eqref{Ming2ca} and \eqref{Ming2e}:}


\begin{align} \nonumber
&2\int_{0}^{\infty}\left(\frac{\Im \! \left(r^{i\,v}\,\zeta\left(\sigma+i\,v  , w_{h} \right)\right) \cos \! \left(\pi \,\sigma \right) \sinh \! \left(\pi \,v \right)}{\cos \! \left(2\,\pi \,\sigma \right)-\cosh \! \left(2\,\pi \,v \right)}+\frac{\Re \! \left(r^{i\,v}\,\zeta\left(\sigma+i\,v  , w_{h} \right)\right) \sin \! \left(\pi \,\sigma \right) \cosh \! \left(\pi \,v \right)}{\cos \! \left(2\,\pi \,\sigma \right)-\cosh \! \left(2\,\pi \,v \right)}\right) \cosh \! \left(v\,\phi \right)d v  \\ \nonumber
& = 
r^{1-\sigma} \left(\left(-r\,\sin \! \left(\phi \right) \overset{\infty}{\underset{j =0}{\sum}}\! \frac{1}{\left(j +r\,\cos \! \left(\phi \right)+w_{h}\right)^{2}+r^{2} \left(\sin^{2}\left(\phi \right)\right)}+\phi \right) \sin \! \left(\phi  \left(1-\sigma \right)\right)\right. \\ & \left.
-\left(\gamma+\ln \! \left(r \right) +\frac{r\,\cos \! \left(\phi \right)+w_{h}}{2\,w_{h}\,r\,\cos \! \left(\phi \right)+r^{2}+w_{h}^{2}} -\overset{\infty}{\underset{j =1}{\sum}}\! \frac{r \left(2\,w_{h}+j \right) \cos \! \left(\phi \right)+r^{2}+w_{h} \left(w_{h}+j \right)}{\,j \left(\left(j +r\,\cos \! \left(\phi \right)+w_{h}\right)^{2}+r^{2} \left(\sin^{2}\left(\phi \right)\right)\right)}\right) \cos \! \left(\phi  \left(\sigma-1 \right)\right)\right)
\label{LhsRb}
\end{align}
and

\begin{align} \nonumber
&2 \int_{0}^{\infty}\left(\frac{\Re \! \left(r^{i\,v}\,\zeta\left(\sigma+i\,v , w_{h}\right)\right) \cos \! \left(\pi \,\sigma \right) \sinh \! \left(\pi \,v \right)}{-\cos \! \left(2\,\pi \,\sigma \right)+\cosh \! \left(2\,\pi \,v \right)}-\frac{\Im \! \left(r^{i\,v}\,\zeta \left(\sigma+i\,v , w_{h}\right)\right) \cosh \! \left(\pi \,v \right) \sin \! \left(\pi \,\sigma \right)}{-\cos \! \left(2\,\pi \,\sigma \right)+\cosh \! \left(2\,\pi \,v \right)}\right) \sinh \! \left(v\,\phi \right)d v  \\ & \nonumber
 = 
r^{1-\sigma} \left(\left(-r\,\sin \! \left(\phi \right) \overset{\infty}{\underset{j =0}{\sum}}\! \frac{1}{\left(j +r\,\cos \! \left(\phi \right)+w_{h}\right)^{2}+r^{2} \left(\sin^{2}\left(\phi \right)\right)}+\phi \right) \cos \! \left(\phi  \left(\sigma-1 \right)\right) \right. \\ & \left. 
-\left(\gamma+\ln \! \left(r \right) +\frac{r\,\cos \! \left(\phi \right)+w_{h}}{2\,w_{h}\,r\,\cos \! \left(\phi \right)+r^{2}+w_{h}^{2}}-\overset{\infty}{\underset{j =1}{\sum}}\! \frac{r \left(2\,w_{h}+j \right) \cos \! \left(\phi \right)+r^{2}+w_{h} \left(w_{h}+j \right)}{j \left(\left(j +r\,\cos \! \left(\phi \right)+w_{h}\right)^{2}+r^{2} \left(\sin^{2}\left(\phi \right)\right)\right)}\right) \sin \! \left(\phi  \left(\sigma-1\right)\right)\right)
\label{LhsIa}
\end{align}

where $w_{h}\equiv w+1/2$. \eqref{LhsRb} and \eqref{LhsIa} were easily tested numerically for $\phi<<\pi$.\newline

\item{ The identity corresponding to \eqref{Rparta} with $w=1/2$ and $r\neq m$, follows:}
\begin{align} \nonumber
2\int_{0}^{\infty}& \frac{\left(\Im \! \left(\zeta \! \left(i\,v +\sigma \right) r^{i\,v}\right) \cos \! \left(\pi \,\sigma \right) \sinh \! \left(\pi \,v \right)+\Re \! \left(\zeta \! \left(i\,v +\sigma \right) r^{i\,v}\right) \sin \! \left(\pi \,\sigma \right) \cosh \! \left(\pi \,v \right)\right) \cosh \! \left(\pi \,v \right)}{\cos \! \left(2\,\pi \,\sigma \right)-\cosh \! \left(2\,\pi \,v \right)}d v \\
& = 
r^{1-\sigma}\left(\pi \,\sin \! \left(\pi \,\sigma \right)+\left(-\psi \! \left(1-r \right)+\ln \! \left(r \right)\right) \cos \! \left(\pi \,\sigma \right)\right) 
\label{Rpartb}
\end{align}

\item{ The following is the identity \eqref{MingWha} previous to evaluating the limit $\phi\rightarrow\pi$ and therefore gives the relationship between the left-hand side of \eqref{MingWha} at a constant value of $r$ for all $-\pi\leq\phi\leq\pi$. This identity was numerically verified for $|\phi|<<\pi$.} 
\begin{align} \nonumber
\int_{0}^{\infty}&\frac{-\left(\sin \! \left(v\,\ln \! \left(m +1\right)\right) \zeta_{I} \! \left(\frac{1}{2}+i\,v \right)+\cos \! \left(v\,\ln \! \left(m +1\right)\right) \zeta_{R} \! \left(\frac{1}{2}+i\,v \right)\right) \cosh \! \left(\phi \,v \right)}{\cosh \! \left(\pi \,v \right)}d v\\ \nonumber
& = 
\left(\gamma -\overset{\infty}{\underset{j\neq m}{\underset{j =1}{\sum}}}\! \frac{1+j +\left(m +1\right) \cos \! \left(\phi \right) j +2 \left(m +1\right) \cos \! \left(\phi \right)+\left(m +1\right)^{2}}{j \left(j^{2}+2\,j +1+2 \left(m +1\right) \cos \! \left(\phi \right) j +2 \left(m +1\right) \cos \! \left(\phi \right)+\left(m +1\right)^{2}\right)}\right. \\ & \left.  \nonumber
+\frac{\left(2\,m^{3}+2\,m^{2}-4\,m -4\right) \cos \! \left(\phi \right)-m^{3}-2\,m^{2}-4\,m -4}{2\,m \left(\left(2\,m +2\right) \cos \! \left(\phi \right)+m^{2}+2\,m +2\right) \left(m +1\right)}+\ln \! \left(m +1\right)\right) \sqrt{m +1}\,\cos \! \left(\frac{\phi}{2}\right) \\ & \nonumber
-\left(-\frac{\sin \! \left(\phi \right) \left(\left(m^{2}+3\,m +2\right) \cos \! \left(\phi \right)+\frac{3\,m^{2}}{2}+3\,m +2\right)}{\left(\cos \! \left(\phi \right)+1\right) \left(\left(2\,m +2\right) \cos \! \left(\phi \right)+m^{2}+2\,m +2\right) \left(m +1\right)}+\phi \right. \\ &\left. 
-\left(m +1\right) \sin \! \left(\phi \right) \overset{\infty}{\underset{j\neq m}{\underset{j =1}{\sum}}}\! \frac{1}{2 \left(m +1\right) \left(1+j \right) \cos \! \left(\phi \right)+\left(m +1\right)^{2}+\left(1+j \right)^{2}}\right) \sqrt{m +1}\,\sin \! \left(\frac{\phi}{2}\right)
\label{MingH1bX}
\end{align}
\end{itemize}

\subsection{Constants} \label{sec:Consts}
\begin{center}
\captionof{table}{  The following give numerical values for integrals defined in the text. All were calculated by both Mathematica and Maple using 50 significant digits. They define integrals specified in \eqref{Hdef}, \eqref{BId}, \eqref{Z0def}, \eqref{Z1def} and \eqref{Z2def} respectively.\newline
\label{Tab:Cdefs}}
\begin{tabular}{||l||l||} 
\hline
\rule[-1ex]{0pt}{2.5ex} H & 0.73128770329627095609488699144402008790399467173125\dots \\ 
\hline  
\rule[-1ex]{0pt}{2.5ex} B & 0.011631879936622239647933640848662374256023406252918\dots \\ 
\hline 
\rule[-1ex]{0pt}{2.5ex} $Z_{0}$ & 0.0080562886088479438618780062072017192398728472183716\dots \\ 
\hline 
\rule[-1ex]{0pt}{2.5ex} $Z_{1}$ & 0.070086041002631260196371300230060136931855307780993\dots \\ 
\hline 
\rule[-1ex]{0pt}{2.5ex} $Z_{2}$ & -0.034562106682472096758925963508715396071235328556772\dots \\ 
\hline 
\end{tabular}\newline

\flushleft

\section{Analytic Continuation and Regularization - A review}

\subsection{Analytic Continuation: An Example from the pages of history.} \label{sec:Analogue}

In 1794, Euler \cite[page 342, as cited by Nielsen \cite{Nielsen}]{EulerIV} studied the following integral, valid for $a>0$ and $b\geq0$:

\begin{equation}
J(a,b,s)\equiv\int_{0}^{\infty}{\mathrm e}^{-a\,t}\,\cos \! \left(t\,b \right) t^{s -1}d t
 = 
\frac{\Gamma \! \left(s \right) \cos \! \left( s\,\arctan \! \left(b/a\right) \right)}{\left(a^{2}+b^{2}\right)^{s/2}}
\label{Euler1}
\end{equation}
and considered the value at $a= 0$. Since the integral is defined for $a>0$, and both sides of \eqref{Euler1} are valid analytic representations of $J(a,b,s)$ over a continuous domain ($a>0$), then by the principle of analytic continuation, the left-hand side obtains its value from the right-hand side in the limit $a\rightarrow 0^{+}$. So, Euler, after considerable analysis, which {\bf did not} apply the principle of analytic continuation, arrived at the now well-known \cite[Eq. (3.761.9)]{G&R} result:

\begin{equation}
\int_{0}^{\infty}\cos \! \left(t\,b \right) t^{s -1}d t = 
\cos \! \left(\frac{\pi \,s}{2}\right) b^{-s}\,\Gamma \! \left(s \right)\,\hspace{1cm}0<s<1,~b>0\,.
\label{Euler2}
\end{equation}

Hence, if $s=1/2$, 
\begin{equation}
\int_{0}^{\infty}\frac{\cos \! \left(v \right)}{\sqrt{v}}d v =\sqrt{\frac{\pi}{2}}\,, 
\label{Euler3}
\end{equation}
after a simple change of variables. Since it is not immediately clear that the left-hand side of \eqref{Euler3} is convergent, this identity must have been a considerable tour-de-force at the time. For our purposes, it simply demonstrates how a symbolic entity -  the left-hand side- can obtain a value when it (possibly) lies outside its range of validity.\newline

Now, let us ask, what happens if $s=1$? Again, by the principle of analytic continuation, this time in the variable $s$, since both sides of \eqref{Euler2} are analytic and equal over a continuous domain ($0<s<1$) as originally defined, they both represent the same function when $s=1$, giving
\begin{equation}
\int_{0}^{\infty}\cos(t)\,dt=0\,.
\label{Euler4}
\end{equation}
Although it is difficult to determine if the left-hand side of \eqref{Euler4} converges to zero numerically, its meaning is symbolically clear, and its value is obtained from the right-hand side. A similar argument would apply to the case $s>1$, where the right-hand side of \eqref{Euler2} assigns a value to the left-hand side, which does not converge. Compare to \eqref{CR5}:
\item{}
\begin{equation} \nonumber
\int_{0}^{\infty}\zeta_{R} \! \left(1/2+iv \right)d v =\, -\pi\,,
\end{equation}
obtained in a similar fashion and recall that, although $\zeta_{R} \! \left(1/2+iv \right)$ is oscillatory, it is mostly positive, as has been widely observed.

\subsection{Regularization - A Quick Review} \label{sec:Reg}
Anyone perusing the symbolic entity -- the left-hand side of \eqref{Euler4} -- who was unfamiliar with \eqref{Euler1}, would not be aware that the integral had obtained its value as a limit point of a hidden variable, and might consider it a challenge to obtain that result rigorously, based solely on the (largely symbolic) integral as written. The somewhat ad hoc method of regularization provides a way out of this impasse. Formally, one could write

\begin{equation}
\int_{0}^{\infty}\cos(t)\,dt={\underset{s\rightarrow 1} {\lim}}\int_{0}^{\infty}\cos \! \left(t \right) t^{s -1}d t = 
{\underset{s\rightarrow 1} {\lim}}\cos \! \left(\frac{\pi \,s}{2}\right) \,\Gamma \! \left(s \right)=0
\label{Euler2a}
\end{equation}
and obtain an answer, arguing that analytic continuation in the variable $s$ is valid if \eqref{Euler2a} were to be read from right to left. However, there are an infinite number of ways of doing this. In the example under consideration, one could alternatively write
\begin{equation}
\int_{0}^{\infty}\cos(t)\,dt={\underset{a\rightarrow 0} {\lim}}\int_{0}^{\infty}{\mathrm e}^{-t\,a}\,\cos \! \left(t\right)\,d t=
{\underset{a\rightarrow 0} {\lim}}\frac{\cos \! \left( \arctan \! \left(1/a\right) \right)}{\left(a^{2}+1\right)^{1/2}}=0\,,
\label{Euler2b}
\end{equation}
and again justify the result reading from right to left and citing the principle of analytic continuation, this time in the variable $a$.\newline

The two examples given above represent two of the many possible implementations of regularization. In addition to Cesàro regularization (see Figure \ref{fig:Cesaro2}), other popular methods, chosen according to the properties of the function $f(x)$ under consideration, are as follows:
\begin{itemize}
\item{from Quantum Electrodynamics (see also \eqref{Iv1}):}
\begin{equation}
\int_{0}^{\infty}f(x)dx = {\underset{T\rightarrow \infty} {\lim}}\int_{0}^{T}f(x)dx\,;
\end{equation}
\item{Principal Value Regularization with $a\in\Re,~a>0$ (see also Section \ref{sec:Pv}):}

\begin{equation}
\int_{0}^{\infty}f(x)/(x-a)dx = {\underset{\epsilon\rightarrow 0} {\lim}}\int_{0}^{\infty}f(x)/(x-a\pm i\epsilon)dx\,;
\end{equation}

\item{Analytic Regularization where $\Re(a)>0,~\nu\in\Re$:}
\begin{equation}
\int_{0}^{\infty}f(x)/(x+a)^n\,dx = {\underset{\nu\rightarrow n} {\lim}}\int_{0}^{\infty}f(x)/(x+a)^{\nu}dx\,;
\end{equation}

\item{and Dimensional Regularization, from Quantum Chromodynamics:}
 
\begin{equation}
\int_{0}^{\infty}f(x)\,d^{\,n}x = {\underset{\epsilon\rightarrow n} {\lim}}\int_{0}^{\infty}f(x)d^{\epsilon}x\,.
\end{equation}
\end{itemize}

\end{center}
\end{appendices}

\end{flushleft}
\end{document}